\documentclass[12pt]{article} 

\usepackage{amsfonts}
\usepackage{xr}
\usepackage{graphicx}
\usepackage{amsmath}
\usepackage{epsfig}
\usepackage{color}

\def\d{{\partial}}

\def\CC{{\mathbb C}}
\def\C{{\mathbb C}}
\def\A{{\mathbb A}}
\def\AA{{\mathbb A}}

\def\DD{{\mathbb D}}
\def\D{{\mathbb D}}

\def\NN{{\mathbb N}}
\def\S{{\mathbb S}}

\def\RR{{\mathbb R}}
\def\R{{\RR}}

\def\TT{{\mathbb T}}
\def\T0{{\mathbb T}_{x_0}}
\def\T{{\mathbb T}}

\def\cF{{\cal F}}

\def\trait (#1) (#2) (#3){\vrule width #1pt height #2pt depth #3pt}
\def\fin{\hfill\trait (0.1) (5) (0) \trait (5) (0.1) (0) \kern-5pt 
\trait (5) (5) (-4.9) \trait (0.1) (5) (0)}
\def\tr{\mbox{tr}}
\def\supess{\mathop{\mbox{ess sup}\,}}

\newtheorem{property}{Property}

\newtheorem{remarque}{Remark}

\newtheorem{corollary}{Corollary}

\newtheorem{proposition}{Proposition}

\newtheorem{lemma}{Lemma}

\newtheorem{theorem}{Theorem}

\newcommand{\Bk}{\color{black}}%
\newcommand{\boite}{\mbox{} \hfill \mbox{\rule{2mm}{2mm}}}

\begin{document}

\author{\begin{tabular}{ccc}
Laurent Baratchart\footnote{INRIA Sophia-Antipolis, Team APICS; Laurent.Baratchart@inria.fr.}, &  
Yannick Fischer\footnote{INRIA Bordeaux, Team Magique3D; Yannick.Fischer@inria.fr;
corresponding author: INRIA-Bordeaux Sud Ouest, UFR Sciences, B\^{a}timent B1, Avenue de l'Universit\'e, BP 1155, 64013 Pau cedex, France; tel: +33 5 40 17 75 40, fax: +33 5 40 17 51 50.}, &  
Juliette Leblond\footnote{INRIA Sophia-Antipolis, Team APICS; Juliette.Leblond@inria.fr.}
\end{tabular}}

\title{Dirichlet/Neumann problems and Hardy classes for the planar conductivity equation}

\maketitle

\noindent{\small{{\bf Abstract.} }} 
We study Hardy spaces $H^p_\nu$ of the conjugate Beltrami equation 
$\overline{\partial} f=\nu\overline{\partial f}$ over Dini-smooth
finitely connected domains, for real contractive
$\nu\in W^{1,r}$ with $r>2$, in the range $r/(r-1)<p<\infty$. 
We develop a theory of conjugate functions and apply it to solve 
Dirichlet and Neumann problems for the conductivity equation
$\nabla.(\sigma \nabla u)=0$ where $\sigma=(1-\nu)/(1+\nu)$. 
In particular situations, we also consider some density properties of traces of solutions together with boundary approximation issues.


\section{Introduction} 

\label{sec:pb}
The foundations  of pseudoanalytic 
function theory, that generalizes some key features of classical holomorphic
function theory,  go back to \cite{bn, vekua} and were historically 
applied to
boundary value problems for partial differential equations. 
This  has recently  been a topic of renewed interest \cite{ap,EfWend, imsurv, krav,Wen}, and reference \cite{BLRR} was apparently first to investigate the 
connections between generalized Hardy spaces on simply connected domains
and Dirichlet problems for the planar conductivity equation with 
$L^p$ boundary data. In \cite{EfendievRuss},
part of this material was carried over to annuli under  more general 
regularity assumptions on the coefficients, and used there  to approach 
certain mixed Dirichlet-Neumann problems. 
The present paper expands and generalizes the results of 
\cite{BLRR} to finitely connected domains under still weaker,
conjectured optimal regularity assumptions. We also take up Neumann problems 
and discuss some density issues for traces of solutions.

We shall consider  a simple class of  pseudoanalytic functions 
(also called generalized analytic functions), namely those satisfying the
conjugate Beltrami equation:

\begin{equation}
\label{dbar} \tag{CB}
\overline\partial f  = \nu  \, \overline{\partial f} \mbox{ a.e. in } \Omega
\,, 
\end{equation}
where $\Omega\subset \R^2 \simeq \C$ is a Dini-smooth domain.
The  dilation coefficient $\nu$ is
{\it real valued} and lies in a Sobolev class
$W^{1,r}(\Omega)$, $2 < r \leq \infty$, while satisfying a 
uniform bound of the type 
\begin{equation} \label{kappa} \tag{$\kappa$}
\left\Vert \nu\right\Vert_{L^\infty(\Omega)}\leq \kappa \mbox{ for some }
\kappa\in (0,1) \, .
\end{equation}
If one writes $f=u+iv$ with real $u$, $v$,  
then (\ref{dbar}) becomes a pair of equations generalizing the
Cauchy-Riemann system:
\begin{equation} \label{system}
\left\{
\begin{array}{l}
\partial_xv=-\sigma\partial_y u,\\
\partial_yv=\sigma\partial_xu,
\end{array}
\right.
\end{equation}
with 
\begin{equation} \label{dfnsigma}
\sigma=\frac{1-\nu}{1+\nu} \, .
\end{equation} 
Now, a compatibility condition for \eqref{system} is the planar 
conductivity equation:
\begin{equation} \label{div}
\mbox{div} (\sigma\nabla u)=0\mbox{ in }\Omega \, ,
\end{equation}
and this fact is the main motivation to study \eqref{dbar}. An interesting 
connection, this time  to a Schr\"odinger equation, was also 
pointed out in \cite{SylUhl}.
Observe that \eqref{kappa} is equivalent {\it via} \eqref{dfnsigma} to the 
ellipticity condition
\begin{equation} \label{ellipticsigma}
0<c \leq \sigma \leq C<\infty \ \mbox{ a.e. in } \Omega
\end{equation}
for some constants $c,C$.
Note also, if $u$ meets (\ref{div}), that the $\sigma$-conjugate function
$v$ satisfies
\begin{equation} \label{system2}
\mbox{div} \, \left(\frac 1{\sigma}\nabla v\right)=0\mbox{ in }\Omega \, .
\end{equation}
We shall study Hardy classes of solutions to \eqref{dbar} 
(see definition in Section \ref{dfnhpnu}), analyze their boundary 
behaviour and give a complete description 
of $\sigma$ conjugate functions in this context. We shall apply the results to 
the Dirichlet problem of equation \eqref{div}, and it will turn out that 
data in $L^p(\partial\Omega)$ for $p> r/(r-1)$ are exactly boundary 
values of solutions satisfying a Hardy condition. 
Trading smoothness of the boundary for smoothness of the 
coefficients, we also give an application to the Dirichlet problem 
with Lipschitz coefficients on piecewise 
$C^{1,\lambda}$ domains without outward pointing cusp, 
when the boundary data are integrable with respect to harmonic measure;
inward pointing cusps are allowed, so that the domain may not be Lipschitz.
In addition, we solve the Neumann 
problem with data $\sigma\partial u/\partial n\in W^{-1,p}(\partial\Omega)$.

From the point of view of regularity
theory, and though we deal with two dimensions and scalar conductivity only,
it is noteworthy that our assumptions are
not covered by the Carleson condition set up in \cite{DPP,KP}. As we rely 
rather extensively on complex methods,
higher dimensional analogs of our results, if true at all, require new ideas 
of proof.

The authors' motivation for such a study originates in certain free boundary 
problems of Bernoulli type for equation \eqref{div} that arise naturally 
when trying to locate  the boundary of a plasma at equilibrium in a tokamak 
\cite{blum}. These are genuine 2-D problem, due to
rotational symmetry. Their approach {\it via} extremal problems,
developed in \cite{TheseYannick,fl,flps}, raises some density
issues for traces of solutions to \eqref{dbar} on subsets of
$\partial\Omega$ which are interesting in their own right and deserves further studies \cite{vlongue}. 

The paper is organized as follows. After some preliminaries
on conformal mapping and Sobolev spaces in Section
\ref{sec:notations_generales}, we introduce in Section
\ref{dfnhpnu} Hardy classes $H^p_\nu(\Omega)$ of equation
\eqref{dbar}, along with their companion Hardy classes $G^p_\alpha(\Omega)$ of
equation \eqref{eq:w} which are of great technical importance for our 
approach. Dwelling on classical works \cite{bn, vekua} to make connection with 
holomorphic Hardy spaces, we then derive  the main properties of  
$H^p_\nu(\Omega)$. Section \ref{secGA} is devoted to a decomposition theorem 
which identifies Hardy classes over finitely connected domains  with sums of 
Hardy classes over simply connected domains, much like in the holomorphic case.
In Section \ref{secdir}, we deal with analytical and topological conditions
for the existence of $\sigma$-conjugate functions, and we apply our results to
the Dirichlet and Neumann problems for equation \eqref{div}. 
Finally, we discuss in Section \ref{sec:density}
some density properties of traces of $H^p_\nu(\Omega)$-functions
on $\partial\Omega$ which are relevant to inverse boundary value problems.
Concluding remarks are given in Section \ref{sec:conclu}.
We append in Appendix some of the more technical results and proofs.

\section{Notations and preliminaries}
\label{sec:notations_generales}

We put $\overline{\C}:=\C\cup\{\infty\}$ for the extended complex plane, 
which identifies to the unit sphere $\S^2$ under stereographic projection.
We let $\D_r$ and $\T_r$ designate the open 
disk and the circle centered at 0 of radius $r$;
when $r=1$ we omit the subscript. More generally, 
$\T_{a,r}$ (resp. $\D_{a,r}$) indicates the circle (resp. open disk) of center 
$a$ and radius $r$.

For $\varrho \in (0,1)$, we set 
$\A_\varrho := \D \setminus \overline{\D}_\varrho$ to be the annulus lying
between $\TT$ and $\T_\varrho$. 
 For more general annuli  we use the notation
$\A_{\varrho_1,\varrho_2}:=\{z;\ \varrho_1<|z|<\varrho_2\}$.
If $\Omega$ is a doubly connected domain 
such that no component of
$\overline{\C}\setminus\Omega$ reduces to a single point, it is well known
there is a unique $\varrho$ making $\Omega$  conformally equivalent to
$\A_\varrho$ \cite[Thm VIII.6.1]{SafTot}. More generally, any finitely 
connected 
domain whose complement is infinite is conformally equivalent to $\D$  with 
finitely many closed disks removed (some of which may degenerate to 
points)\footnote{Indeed, any finitely connected domain is conformally 
equivalent to a domain whose boundary consists of
circles or points \cite[Sec. V.6, Thm 2]{Goluzin}; if the complement 
is infinite there is at least one circle
whose interior can be mapped onto $\overline{\C}\setminus\overline{\D}$ by a
M\"obius tranform.}.
Such a domain will be termed a normalized circular domain. 
Moreover, the above conformal map is unique up to a
M\"obius transformation.

Recall that a function $h$ is called {\emph {Dini-continuous}} if 
$\int_0^\varepsilon(\omega_h(t)/t)dt<+\infty$ for some, 
hence any $\varepsilon>0$, 
where $\omega_h$ is the modulus of continuity of $h$. 
A function is {\emph {Dini-smooth}} if it has Dini-continuous 
derivative.
A  domain $\Omega\subset\overline{\C}$ is said to be Dini-smooth 
if its boundary $\partial\Omega$ lies in $\C$ and
consists  of finitely many  Jordan 
curves with nonsingular Dini-smooth parametrization. Note that a Dini-smooth 
domain is finitely connected by definition, and it contains $\infty$ if 
it is unbounded.

Any conformal map between Dini-smooth domains extends to a homeomorphism 
of their closures, and the derivative also extends continuously to the closure
of the initial domain in such a way that it is never zero,
{\it cf.} Lemma \ref{conf} in Appendix \ref{app1}.

We orient the boundary of a Dini-smooth domain $\Omega$  
in a canonical way, 
i.e. $\Omega$ lies on the left side when moving along $\partial \Omega$, 
and the unit normal $\vec{n}$ points outward. 

We denote interchangeably (the differential of) planar Lebesgue measure by
\[
 dm(\xi) = dt_1dt_2 = (i/2) \ d\xi \wedge d \overline{\xi},\qquad \xi=t_1+it_2.
\]



Given a domain $\Omega\subset\C$, we put ${\cal D}(\Omega)$ for the space of 
$C^\infty$-smooth
complex valued 
functions with compact support in $\Omega$, equiped with the usual 
topology\footnote{{\it i.e.} the inductive topology of its subspaces  
$\mathcal{D}_K$ comprised of functions whose support lies in a compact set 
$K$, each $\mathcal{D}_K$ being topologized by uniform convergence of 
a function and all its partial derivatives \cite[Sec. I.2]{Schwartz}).}.
Its dual $\mathcal{D}'(\Omega)$ 
is the space of distributions on $\Omega$. 
For $1\leq p\leq \infty$ and $k\in\NN$,  
we let $L^p(\Omega)$ and $W^{k,p}(\Omega)$ 
be the familiar Lebesgue and Sobolev spaces  with respect to $dm$;
we sometimes write $L^p_\RR(\Omega)$ or $W_\RR^{k,p}(\Omega)$
to emphasize restriction to real-valued functions.
Recall that $W^{1,p}(\Omega)$ consists of functions in
$L^p(\Omega)$ whose distributional derivatives lie in
$L^p(\Omega)$ up to order $k$. Actually we only need $k=1,2$,
the norms on $W^{1,p}(\Omega)$, $W^{2,p}(\Omega)$ being defined as
\[
 \|f\|_{W^{1,p}(\Omega)} = \|f\|_{L^p(\Omega)}+\|\partial f\|_{L^p(\Omega)}+\|\overline{\partial} f\|_{L^p(\Omega)},
\]
\[\|f\|_{W^{2,p}(\Omega)}=\|\partial f\|_{W^{1,p}(\Omega)}+
\|\overline{\partial} f\|_{W^{1,p}(\Omega)}+\|f\|_{L^p(\Omega)},\]
where $\partial$ and $\overline{\partial}$ stand for the usual
(distributional) complex derivatives, i.e. 
\[
 \partial f := \partial_z f = \dfrac{1}{2}(\partial_x - i \partial_y)f~~~\mbox{and}~~~\overline{\partial} f := \partial_{\overline{z}} f = \dfrac{1}{2}(\partial_x + i \partial_y)f\qquad z=x+iy.
\]
Note the obvious identity: $\overline{\partial f}=\overline{\partial} 
\,\,\overline{f}$. The closure of $\mathcal{D}(\Omega)$ in $W^{1,p}(\Omega)$ is
denoted by $W^{1,p}_0(\Omega)$. Recall the notation 
$W^{-1,p}(\Omega) = (W^{1,q}_0(\Omega))^*$, $1/p+1/q=1$.
For basic properties of Sobolev spaces that we use, 
see {\it e.g.} \cite{brezis,ziemer}. 

When $\Omega\subset\C$ is a bounded Dini-smooth domain,  
$L^p(\partial\Omega)$ is understood with respect to normalized arclength
and $W^{1,p}(\partial\Omega)$ is naturally defined using local coordinates, 
since Lipschitz-continuous changes of variable preserve Sobolev classes.
Each 
$f \in W^{1,p}(\Omega)$ with $1<p\leq\infty$  has a trace on 
$\partial \Omega$, denoted 
by $\tr_{\partial \Omega} \ f$, lying in the so-called fractional Sobolev space
$W^{1-1/p,p}(\partial\Omega)$. The latter is a real interpolation space
between $L^p(\partial\Omega)$ and $W^{1,p}(\Omega)$ of exponent $1-1/p$,
an intrinsic definition of which can be found
in \cite[{ Thm} 7.48]{Adams}. The trace operator defines a continuous
surjection from
$W^{1,p}(\Omega)$ onto  $W^{1-1/p,p}(\partial\Omega)$. By the Sobolev 
embedding 
theorem, each $f \in W^{1,p}(\Omega)$ with $p>2$ is H\"older-smooth of 
exponent  $1-2/p$ on $\Omega$, hence $f$ extends continuously to
$\overline{\Omega}$ in this case. The space
$W^{1,\infty}(\Omega)$ identifies with 
Lipschitz-continuous functions on $\Omega$.

We also introduce the spaces $L^p_{loc}(\Omega)$ and $W^{1,p}_{loc}(\Omega)$ 
of distributions\footnote{These
are topologized by the family of 
semi-norms $\|f_{\Omega_n}\|_{L^p(\Omega_n)}$ and
$\|f_{\Omega_n}\|_{W^{1,p}(\Omega_n)}$ respectively, with $\{\Omega_n\}$ 
a nested family of relatively compact open subset exhausting $\Omega$.} 
whose restriction to any relatively compact open subset 
$\Omega_0$ of $\Omega$ lies in $L^p(\Omega_0)$ or 
$W^{1,p}(\Omega_0)$.
All classes of functions 
we will consider are embedded in $L_{loc}^p(\Omega)$ 
for some $p \in (1,+\infty)$, and solutions to differential
equations are understood in the distributional 
sense. For instance, to define distributions like 
$\nu \overline{\partial f}$ where $\nu \in W^{1,r}_\RR(\Omega)$ and
$f\in L_{loc}^p(\Omega)$ with $1/p+1/r\leq1$,
we use Leibniz's rule:
\[
 <\nu \overline{\partial f},\phi> = -\int_{\Omega} (\nu \overline{f} \,\overline{\partial}\phi + \overline{\partial}\nu \overline{f}\phi) \ dm, ~~~~~\forall \ \phi \in \mathcal{D}(\Omega).
\]
where $<,>$  denotes the duality product between $\mathcal{D}'(\Omega)$  and
$\mathcal{D}(\Omega)$.

If in addition $r>2$ and $\sigma$  meets \eqref{ellipticsigma}
while $u \in W^{1,p}_\RR(\Omega)$ 
solves (\ref{div}), 
the normal derivative $\partial_n u$ is 
the unique member of the dual space 
$W^{-1/p,p}_\RR(\partial \Omega) = (W^{1-1/q,q}_\RR(\partial \Omega))^*$,
$1/p+1/q=1$, such that
\begin{equation}
\label{Greenn}
 <\sigma \partial_n u, \phi >_{\partial\Omega} = 
\int_\Omega \sigma \nabla u \cdot \nabla g \ dm, ~~~ \  
\ g \in W^{1,q}(\Omega),\ \ \tr_{\partial \Omega} \ g = \phi.
\end{equation}
In fact, \eqref{Greenn} defines 
$\sigma\partial_n u\in W^{-1/p,p}(\partial\Omega)$ and, under the stated 
assumptions, multiplication by $\sigma$ is an isomorphism of the latter
because it is an isomorphism of $W^{1-1/q,q}(\partial\Omega)$ ({\it e.g.} by interpolation). Clearly then, $\|\partial_n u\|_{W^{-1/p,p}(\partial\Omega)}\leq 
C(\Omega,\sigma,p)\|u\|_{W^{1,p}(\Omega)}$.

Sobolev spaces are naturally defined on the
Riemann surface $\overline{\C}\sim\S^2$ 
\cite{Hebey}, integration being understood with respect to spherical measure.
We shall not be concerned with intrinsic notions:
to us it suffices to say that if $\Omega\subset\overline{\C}$ is a
(possibly unbounded)  
Dini-smooth domain distinct from $\overline{\C}$, then
it can be mapped onto a bounded Dini-smooth
domain $\Omega'$ by some conformal map $\varphi$ and then
$f\in W^{1,p}(\Omega)$ (resp. $f\in L^p(\Omega)$) if and only if 
$f\circ\varphi^{-1}\in W^{1,p}(\Omega')$ (resp. $L^p(\Omega')$),
with equivalence of norms (the equivalence depends on $\varphi$).
This is consistent with previous definitions, since by Lemma 
\ref{conf} (in Appendix \ref{app1}) conformal maps between
bounded Dini-smooth domains are Lipschitz continuous.
A function in $W^{1,p}(\overline{\C})$ is one whose restriction to any 
proper Dini-smooth subdomain $\Omega$ belongs to $W^{1,p}(\Omega)$.

On a bounded domain $\Omega$, relation \eqref{dbar} may be regarded as a 
differential equation for $1$-forms in local
coordinates on the Riemann surface $\S^2$, namely
$\overline{\partial}fd\bar{z}=\nu\overline{\partial f}d\bar{z}$. 
Subsequently, if $\Omega$, 
$\Omega'$, and $\varphi$ are as before, we say that $f$ solves \eqref{dbar} 
on $\Omega$ if and only if $f\circ\varphi^{-1}$ satisfies a similar equation 
on $\Omega'$ only with $\nu$ replaced by $\nu\circ\varphi^{-1}$; 
this agrees with the complex chain rule when $\Omega$ is bounded
\cite[Sec. 1.C]{ahlfors},
and allows us to make sense of \eqref{dbar} when $\Omega$ is unbounded.

If $f$ is a function defined on $\Omega$,
the symbol $f_{|\Omega_1}$ indicates the restriction of $f$ to 
$\Omega_1\subset\Omega$. 
Whenever $f$ is defined on $\Omega_1$ and $h$ on $\Omega_2 = \Omega \setminus \Omega_1$, the notation $f \vee h$ is used for the concatenated function 
defined on $\Omega$ which is equal to $f$ on $\Omega_1$  and to $h$ on $\Omega_2$.

We let $\partial_t$ and $\partial_n$ denote respectively the tangential and 
normal derivatives of a function at a smooth point on a rectifiable curve.
As became customary, the same symbol ({\it e.g.}
``$C$'') is used many times to mean
different constants.

\section{Generalized Hardy classes}
\label{dfnhpnu}

Hardy classes of equation \eqref{dbar}  over a
bounded Dini-smooth simply connected domain were introduced
in \cite{BLRR}. Their study is twined with that of Hardy classes 
of equation \eqref{eq:w} 
further below, whose connection to \eqref{dbar} was originally stressed in
\cite{bn}.  
Hardy classes of \eqref{dbar} over bounded annular domains 
with analytic boundary 
have  subsequently been defined in \cite{EfendievRuss}.
This section carries out their generalization to arbitrary 
Dini-smooth domains in $\overline{\C}$. 

Although \cite{BLRR} restricts to 
the case where $\nu\in W_\R^{1,\infty}(\Omega)$, 
it was 
observed in \cite{EfendievRuss} that many results
still hold when $\nu\in W_\R^{1,r}(\Omega)$ 
for some $r>2$, provided that $p>r/(r-2)$. We improve on this throughout 
by assuming $r>2$ and $p>r/(r-1)$, which we conjecture is
the optimal range of exponents for the validity of whatever follows.
When $r=\infty$, and only in 
this case, we  cover the whole range of exponents $1<p<\infty$.


To recap, our working assumptions 
will be that 
$\Omega\subset\overline{\C} \mbox{ is 
Dini-smooth}$ (in particular finitely connected) and that
\begin{equation}
\label{hypothesesnup}
\nu \mbox{ meets \eqref{kappa}}, \qquad \nu\in W_\R^{1,r}(\Omega)\ \mbox{for some}\  r\in(2,+\infty],\qquad r/(r-1)<p<+\infty.
\end{equation}

Note the assumptions on $\nu$ are equivalent to require that
$\sigma$ given by \eqref{dfnsigma} 
lies in  $W_\R^{1,r}(\Omega)$ and satisfies \eqref{ellipticsigma}.

\subsection{$H_\nu^p(\Omega)$}
\label{sec:defHpnu}

In analogy to classical holomorphic Hardy spaces, the Hardy space 
$H^p_\nu(\D)$ was defined in \cite{BLRR} to consist
of those functions $f$ in $L^p(\D)$  satisfying 
\eqref{dbar}  in the sense of distributions and such that 
 \begin{equation} \label{esssuplp}
 \left\Vert f\right\Vert_{H^p_{\nu}(\D)}:= \supess_{0<r<1}
 \left\Vert f\right\Vert_{L^p(\T_r)}=
\supess_{0<r < 1}
 \left(\int_0^{2\pi} \left\vert
 f(re^{i\theta})\right\vert^p\frac{d\theta}{2\pi}\right)^{1/p} <+\infty \,. 
 \end{equation}
Likewise, in \cite{EfendievRuss}, the Hardy space $H_\nu^p(\A_{\varrho})$
was set to be comprised of functions in $L^p(\A_\varrho)$ solving \eqref{dbar}
and such that
 \begin{equation}
\label{boundHpA}
\left\Vert f\right\Vert_{H^p_{\nu}(\A_\varrho)}:=  \supess_{\varrho< r < 1}
 \left\Vert f\right\Vert_{L^p(\T_r)}= \supess_{\varrho< r < 1}
 \left(\int_0^{2\pi} \left\vert
 f(re^{i\theta})\right\vert^p\frac{d\theta}{2\pi}\right)^{1/p}<+\infty.
\end{equation}
Now, for $\Omega$ a Dini-smooth proper
subdomain of $\overline{\C}$
and $\nu$, $p$ as in \eqref{hypothesesnup},
we define $H^p_\nu(\Omega)$ to be comprised of those $f\in L^p_{loc}(\Omega)$
solving \eqref{dbar} in the sense of distributions for which
there is a sequence of domains $\Delta_n$ with 
$\overline{\Delta}_n\subset\Omega$, whose boundary $\partial\Delta_n$ 
is a finite union of rectifiable Jordan curves of uniformly bounded length, 
such that
each compact subset of $\Omega$ is eventually contained in $\Delta_n$, and
having the property that 
\begin{equation}
\label{systema}
\sup_{n\in\NN}\|f\|_{L^p(\partial\Delta_n)}<\infty.
\end{equation}
When $\nu\equiv0$, condition \eqref{systema} defines the so-called 
{\emph Smirnov class} of index $p$ of holomorphic functions in $\Omega$, 
which coincides with the Hardy class on Dini-smooth domains\footnote{ The Hardy
class is defined by the condition that $|f|^p$ has a harmonic majorant; 
the two classes coincide as soon as harmonic measure and arclength are 
comparable  up to a multiplicative constant on $\partial\Omega$
\cite[Ch. 10]{duren}, \cite{TuHa2}, which is 
the case for Dini-smooth domains thanks to Lemma \ref{conf}, Appendix \ref{app1}.}. This class we
consistently denote by $H^p(\Omega)$ (no subscript).

It is true, although not immediately clear, that $H^p_\nu(\Omega)$ is a 
vector space and that there is a \emph{fixed} 
sequence $\{\Delta_n\}$ for which (\ref{systema}) yields a complete norm.
It is not obvious either that \eqref{systema} is equivalent to
\eqref{esssuplp} or \eqref{boundHpA} for the disk or the annulus. All
this is known to hold for holomorphic functions \cite{TuHa2}, 
\cite[Sec. 10.5]{duren}, but the proof when $\nu\neq0$
will await Section \ref{secGD}.
Note that $H^p_\nu$ is only a {\emph real} Banach space if $\nu\neq0$.

\vskip .3cm
The definition of $H^p_\nu(\Omega)$  just given is conformally invariant: 
if $\varphi$ conformally maps a Dini-smooth domain $\Omega'$ onto
a Dini-smooth domain $\Omega$, then  $\nu\in W_\R^{1,r}(\Omega)$
and $f\in H^p_\nu(\Omega)$ if and only if
$\nu\circ\varphi\in W_\R^{1,r}(\Omega')$ and
$f\circ\varphi\in H^p_{\nu\circ\varphi}(\Omega')$. Indeed the
$\varphi^{-1}(\Delta_n)$ form an admissible sequence of compact sets in $\Omega'$
since their boundary is eventually contained in a compact neighborhood of 
$\partial\Omega'$ where $|\varphi'|$ is bounded below by a strictly positive 
constant in view of Lemma \ref{conf}.
In \cite{BLRR,EfendievRuss}, conformal invariance was used to \emph{define}
$H^p_\nu(\Omega)$ on simply or doubly connected bounded Dini-smooth 
domains\footnote{These works do not mention the case of unbounded domains, 
but it requires no change as we just stressed. The paper
\cite{EfendievRuss} restricts to analytic boundaries, 
which is also unnecessary thanks to 
Lemma \ref{conf}.}. 

In connection with unbounded domains, the following 
{\emph reflexion principle} is useful:
for $f \in L^p_{loc}(\D)$, set 
\begin{equation}
\label{defcheck}
 \check{f}(z)=\overline{f\left(\frac 1{\overline{z}}\right)} \, , \ \ z\in \C\setminus
 \overline{\D} \, .
\end{equation}
Then 
\begin{equation}
\label{equiv}
f \in H_{\nu}^{p}(\D) \Longleftrightarrow \check{f} \in
 {H_{\check{\nu}}^{p}(\overline{\C} \setminus \overline{\D})} \, . 
\end{equation}
Indeed, if we put $g(z)=\overline{f(\overline{z})}$ and 
$\mu(z)=\nu(\overline{z})$ 
for $z\in\D$, we get by definition
upon  using the conformal map $z\mapsto 1/z$
that $\check{f}\in H^p_{\check{\nu}}(\overline{\C} \setminus \overline{\D})$
if and only if $\mu\in W^{1,r}_\R(\D)$ and $g\in H^p_{\mu}(\D)$. 
Clearly this is the case if and only if $\nu\in W^{1,r}_\R(\D)$ and
$f\in H^p_\nu(\D)$, as follows from \eqref{dbar} by conjugation,
which proves  (\ref{equiv}).

\begin{remarque}
\label{rmq:xcte}
We did not define Hardy spaces of $\overline{\C}$ (a Dini-smooth domain with empty 
boundary), but this case is of little interest since no non-constant 
distributional solution to \eqref{dbar} exists in 
$L^p_{loc}(\overline{\C})$. In fact, 
by Propositions \ref{trick-hardy},
\ref{expsf} further below,
a function $f$ with these properties must lie in
$W^{1,k}(\overline{\C})$ for some $k>2$. In particular it is bounded,
so by the extended Liouville theorem \cite[Cor. 3.4]{ap}, 
$f=Ce^g$, where $C$ is constant and $g$ is continuous on $\overline{\C}$. 
Applying this to $f-f(0)$ we conclude the latter is identically 
zero, as desired.
\end{remarque}



\subsection{$G_\alpha^p(\Omega)$}
\label{secGD}
When $\Omega=\D$ or $\Omega=\A_\varrho$ and $\alpha\in L^r(\Omega)$,
the  Hardy space 
$G^p_\alpha(\Omega)$ was defined in \cite{BLRR,EfendievRuss} to consist of 
those $w\in L^p_{loc}(\Omega)$
such that
\begin{equation}
\label{eq:w}
\bar\partial w=\alpha\overline{w}\ \ \ \mathrm{on}\  \Omega
\end{equation}
in the distributional sense,
and meeting condition \eqref{esssuplp} or \eqref{boundHpA}
(with $w$ instead of $f$).

On a bounded Dini-smooth domain $\Omega\subset\overline{\C}$, given 
$\alpha\in L^r(\Omega)$, we define
$G_\alpha^p(\Omega)$  to  consist of those $w\in L^p_{loc}(\Omega)$
meeting \eqref{systema}  for some admissible sequence 
$\{\Delta_n\}\subset\Omega$, and such that \eqref{eq:w} holds.
If we set further
\begin{equation}
\label{defA}
A(z)=\frac{1}{2i\pi}\int_\Omega \frac{\alpha(\xi)}{\xi-z}\,d\xi\wedge
d\overline{\xi},\qquad z\in\Omega,
\end{equation}
then $A\in W^{1,r}(\Omega)$ by the Sobolev embedding theorem together with
standard properties of the Cauchy and 
Beurling transforms \cite[[Ch. 1, (1.7)-(1.9)]{im2}; moreover
$\alpha=\overline{\partial}A$. Rewriting $G^p_\alpha(\Omega)$ as
$G^p_{\overline{\partial} A}(\Omega)$ is suggestive of a conformally 
invariant definition
valid over arbitrary Dini-smooth domains in $\overline{\C}$, namely
if $\varphi$ conformally
maps a Dini-smooth domain $\Omega$ onto a bounded Dini-smooth domain 
$\Omega'$ and if $A\in W^{1,r}(\Omega)$, then $w\in G^p_{\overline{\partial} A}(\Omega)$
if and only if  $w\circ\varphi^{-1}\in G^p_{\alpha}(\Omega')$
with $\alpha=\overline{\partial}(A\circ\varphi^{-1})$.
 


We will prove momentarily that $G^p_\alpha$ is 
a real Banach space but first we stress the motivation behind its
definition. If we let $A=\log\sigma^{1/2}$ so that 
$\alpha=\bar\partial \log \sigma^{1/2}$, an explicit connection 
between  $H^p_{\nu}(\Omega)$ and $G^p_{\alpha}(\Omega)$ 
stems from a transformation introduced in  \cite{bn}:
\begin{proposition} 
\label{trick-hardy}
Assume that $\Omega \subset\overline{\C} $ is a proper Dini-smooth domain and that 
$\nu$, $p$, $r$ satisfy \eqref{hypothesesnup}.  Let $\sigma$ be as in 
\eqref{dfnsigma} and define $\alpha\in L^r(\Omega)$ by 
\begin{equation}
\label{alpha}
\alpha=-\frac{\bar\partial\nu}{1-\nu^2} = \bar\partial \log \left[\frac{1-\nu}{1+\nu}\right]^{1/2}
= \bar\partial \log \sigma^{1/2}.
\end{equation}
\begin{equation}
\label{ftow} \mbox{Then: }~~~~~~~~~~ f  = u+ i \, v \in H_{\nu}^p(\Omega) \Longleftrightarrow
w =\frac{ f - \nu \overline{f}}{\sqrt{1-\nu^2}} 
= \sigma^{1/2}\, u+ i \, \sigma^{-1/2}\, v \Bk
\in G_\alpha^p(\Omega) \, .
\end{equation}
For $\Delta_n$ as in \eqref{systema}, there are constants 
$C_1,C_2>0$ independent of $f$ and $w$ such that
\begin{equation}
\label{equivHG}
\sup_{n\in\NN}\|f\|_{L^p(\partial\Delta_n)}\leq C_1
\sup_{n\in\NN}\|w\|_{L^p(\partial\Delta_n)}\leq C_2\sup_{n\in\NN}\|f\|_{L^p(\partial\Delta_n)}.
\end{equation}
\end{proposition}
The proof of Proposition \ref{trick-hardy} is a straightforward computation 
using Leibnitz's rule and the fact that $(1-\nu^2)^{-1/2}\in W^{1,r}(\Omega)$
since $r>2$ and $\nu$ satisfies  \eqref{kappa}.

This proposition entails that it is essentially equivalent
to work with $H_{\nu}^p$ or $G_\alpha^p$. However, equation \eqref{eq:w} is
technically simpler to handle because the derivative of 
the first order
({\it i.e.} $\overline{\partial}w$) is expressed in terms of the derivative
of zero-th order ({\it i.e.} $w$).

A primary example of such a simplification is the factorization principle
asserted in Proposition \ref{expsf}
below, which lies at the root of the connections between solutions to
\eqref{eq:w} and holomorphic functions. For slightly smoother classes of 
solutions, this principle goes back to \cite{bn} and
was later extended to accomodate more general planar elliptic equations
including those defining pseudo-analytic functions \cite[Thm 2.3.1]{Wen},
see also \cite{krav}, \cite[Thm 2.1]{EfWend}; it was adapted to Hardy classes
in \cite{BLRR,EfendievRuss}. We provide in the Appendix \ref{app2}
a proof which differs from \cite{BLRR}  in that
normalization  must be argued differently in the multiply 
connected case (compare \cite[Prop. 2.2.1]{EfendievRuss}).
\begin{proposition}
\label{expsf}
 Assume that $\Omega \subset\overline{\C} $ is a proper Dini-smooth domain and that 
$p$, $r$ satisfy \eqref{hypothesesnup}.
Let $A\in W^{1,r}(\Omega)$, $\alpha=\overline{\partial}A$,
and suppose $w\in L^p_{loc}(\Omega)$ is a distributional solution to
\eqref{eq:w}. Then $w$ admits a factorization of the form
\begin{equation}
\label{decompexpH}
 w(z) = \exp(s(z))F(z),~~~~z \in \Omega,
\end{equation}
where $F$ is holomorphic and $s \in W^{1,r}(\Omega)$ satisfies
\begin{equation}
\label{regs}
 \|s\|_{L^{\infty}(\Omega)} \leq C_0\|s\|_{W^{1,r}(\Omega)} 
\leq C_0'\|\alpha\|_{L^r(\Omega)} \leq C_0''\|A\|_{W^{1,r}(\Omega)} \, ,
\end{equation}
$C_0$, $C'_0$, $C''_0$ being strictly positive 
constants depending only on $r$ and $\Omega$.

In particular $w$ belongs to $W^{1,r}_{loc}(\Omega)$,
and $w \in G_\alpha^p(\Omega)$ if and only if $F \in H^p(\Omega)$.

If $\Omega$ is $n$-connected
and $\partial\Omega=\cup_{j=0}^{n} \Gamma_j$ where the $\Gamma_j$ are 
disjoint  Jordan curves, we may choose $s$ so that 
$\mbox{\rm Im}\,s_{|\Gamma_j}=c_j$ where the $c_j$ are constants such that
$\sum_{j=0}^{n}c_j=0$, one of which can be chosen arbitrarily.
\end{proposition}

\begin{remarque}
\label{remHolder}
 It follows from Proposition \ref{expsf} and Sobolev's embedding theorem 
that $s \in C^{0,\gamma}(\overline{\Omega})$, uniformly with 
respect to $w$, and 
$w \in C^{0,\gamma}_{loc}(\Omega)$ with $\gamma=1-2/r$.
\end{remarque}

Proposition \ref{expsf} quickly gives interior estimates for solutions to
\eqref{dbar}:
\begin{corollary}
\label{regint}
 Assume that $\Omega \subset\overline{\C} $ is a proper Dini-smooth
domain and that $\nu$, $p$, $r$ satisfy 
\eqref{hypothesesnup}. Let $f\in L^p_{loc}(\Omega)$ be a distributional
solution to \eqref{dbar}. Then $f\in W^{2,r}_{loc}(\Omega)$.
\end{corollary}
{\it Proof.} Defining $\alpha$ as in \eqref{alpha},
it is straightforward to check  that $w$ given by \eqref{ftow} 
lies in $L^p_{loc}(\Omega)$  and satisfies 
\eqref{eq:w}. By Proposition \ref{expsf}, $w\in W^{1,r}_{loc}(\Omega)$,
hence the same is true of $f$. Using this fact it is
easily verified that the distributional derivative of 
$\nu\overline{\partial f}$ can be computed 
according to Leibnitz's rule. Consequently,
setting $G:=\partial f$ and applying $\partial$ to \eqref{dbar}, 
we  obtain since $\partial$ and
$\overline \partial$ commute, that
$\overline{\partial} G=\nu \partial\overline{G}+
(\partial\nu)\overline{G}$.
As $\nu$ is real, conjugating this relation provides us with
another expression for $\partial\overline{G}$,
and solving for $\overline{\d} G$ after substitution yields 
\[
 \overline{\partial} G=\frac{\nu\overline{\partial}\nu}{1-\nu^2}G+
 \frac{\partial\nu}{1-\nu^2}\overline{G} \, 
\]
from which we deduce that $H = (1-\nu^2)^{1/2} \, G$ satisfies
\begin{equation}
\label{commew}
\overline{\partial} H = \frac{\partial\nu}{1-\nu^2}\overline{H}
\end{equation}
in the sense of distributions. As $H\in L^r_{loc}(\Omega)$, we deduce from 
Proposition \ref{expsf} (applied with $\alpha={\partial\nu}/(1-\nu^2)$)
that $H\in W^{1,r}_{loc}(\Omega)$, implying that the same is true 
of $G=\partial f$ since $W^{1,r}_{loc}(\Omega)$ is an algebra for $r>2$.
From \eqref{dbar} we see that also
$\overline{\partial} f \in W^{1,r}_{loc}(\Omega)$, thereby achieving the proof.
\hfill\boite\\
 
We are now in position to prove that $H^p_\nu(\Omega)$ and 
$G_\alpha^p(\Omega)$ are 
indeed Banach spaces. For this, let 
$\varphi$ conformally map 
$\Omega$ onto  a normalized circular domain $\Omega'$, and 
$\delta_{\Omega'}$ be 1/2 if $\Omega'=\D$ 
and half the minimal distance between 
two components of $\partial\Omega'$ otherwise. Set
\[K_\varepsilon=\{z\in\Omega';\  \mbox{\rm dist}(z,\overline{\C}\setminus\Omega')\geq\varepsilon\,\delta_{\Omega'}
\},\qquad\varepsilon\leq 1, \]
where ${\rm dist}(z,E)$ indicates the distance from $z$ to the set $E$.
Then, $K_\varepsilon$ is a compact subset of $\Omega'$ bounded by 
circles concentric with the components of $\partial\Omega'$. We put
\begin{equation}
\label{defDeltat}
\widetilde{\Delta}_\varepsilon=\varphi^{-1}(K_\varepsilon)\subset\Omega,
\end{equation}
and define
for $g\in H^p_\nu(\Omega)$ or $G^p_\alpha(\Omega)$:
\begin{equation}
\label{systemn}
\|g\|_p=\sup_{n\in\NN}\|g\|_{L^p(\partial\widetilde{\Delta}_{1/n})}.
\end{equation}
Although $\varphi$ is not uniquely defined, 
different $\varphi$ will give rise to equivalent $\|.\|_p$.

We use a lemma on holomorphic functions which is well known on the disk
and the annulus but that we could not ferret out in the literature in the
multiply connected case:
\begin{lemma}
\label{holoHp}
Let $\Omega$ be a Dini-smooth domain and $f$  holomorphic in $\Omega$.
Then $f\in H^p(\Omega)$ if and only if $\|f\|_p<\infty$.
\end{lemma}
The proof of Lemma \ref{holoHp} is given in Section \ref{app3}.
\begin{theorem}
\label{BS}
Assume that $\Omega$ is a proper Dini-smooth
domain in $\overline{\C}$
and that $\nu$, $p$, $r$ satisfy \eqref{hypothesesnup}. 
Let moreover $A\in W^{1,r}(\Omega)$ and $\alpha=\overline{\partial}A$.
\begin{itemize}
\item[(i)] Endowed with \eqref{systemn}, $H^p_\nu(\Omega)$ and $G_\alpha^p(\Omega)$ are real Banach spaces, that coincide with those defined by 
\eqref{esssuplp}and \eqref{boundHpA} when $\Omega=\D$ or $\Omega=\A_\varrho$.
\item[(ii)] It holds that $f\in H^p_\nu(\Omega)$ (resp. $w\in G_\alpha^p(\Omega)$) if, and only if $f$ satisfies \eqref{dbar} (resp. $w$ satisfies \eqref{eq:w})
and $|f|^p$ (resp. $|w|^p$) has a harmonic majorant 
in $\Omega$.
\end{itemize}
\end{theorem}
{\it Proof.} Since $\nu$ meets \eqref{kappa} and
$f=(w+\nu\overline{w})/\sqrt{1-\nu^2}$ by \eqref{ftow},
it is enough to prove the result for $G_\alpha^p(\Omega)$. 

By Proposition \ref{expsf} and Lemma \ref{holoHp}, it holds that 
$w\in G_\alpha^p$ if and only if $\|w\|_p<\infty$, in particular 
$G_\alpha^p$ is a real vector space on which $\|.\|_p$ defines a norm.
To see that it is complete, it is enough to check that a Cauchy sequence $w_k$
has a converging subsequence. Write $w_k=e^{s_k}F_k$ according to   
\eqref{decompexpH}. By Remark \ref{remHolder},
the sequence $s_k$ is equicontinuous,
on $\overline{\Omega}$, therefore some subsequence $s_{l_k}$
converges uniformly 
there. By \eqref{regs} again, 
$\|F_{l_k}\|_{L^p(\partial\widetilde{\Delta}_n)}$ 
is uniformly bounded, hence a normal family argument provides us with a 
subsequence $F_{m_k}$ converging locally uniformly in $\Omega$. 
Then $w_{m_k}$ converges locally uniformly to some $w\in L^p_{loc}(\Omega)$, 
and it
is clear from the definition of distributional derivatives that $w$ solves 
\eqref{eq:w}. Moreover, passing to the limit under the integral sign 
for fixed $n$ in the right hand side of \eqref{systemn} shows that
$w\in G_\alpha^p$ and that $\|w-w_{m_k}\|_p\to 0$ as $k$ goes to infinity. 
This proves
$(i)$.

If $g$ is holomorphic on $\Omega$ and $\varphi$ conformally maps the latter on 
a domain $\Omega'$ with analytic boundaries,
 it is known that \eqref{systema} holds for $g$
for some admissible sequence $\Delta_n$ if, and only if 
$|g\circ\varphi^{-1}|^p(\varphi^{-1})'$ has a harmonic majorant on 
$\Omega'$ \cite[Sec. 10.5]{duren}. In view of Lemma \ref{conf}, 
we conclude 
that \eqref{systema} holds if and only if $|g|^p$ has a harmonic 
majorant on $\Omega$ ({\it i.e.} the so-called Smirnov and Hardy 
classes coincide on Dini-smooth domains). Applying this to 
$F$ in \eqref{decompexpH}, point $(ii)$ now follows from 
Proposition \ref{expsf} since $s$ is bounded.\hfill\boite


\subsection{Basic properties of $G^p_\alpha$ and $H^p_\nu$ classes}
\label{sechpnud}

Below we enumerate some properties that  $G^p_\alpha(\Omega)$ and 
$H^p_\nu(\Omega)$ inherit from $H^p(\Omega)$ {\it via} Proposition 
\ref{expsf}.
These generalize results stressed in \cite{BLRR,EfendievRuss}
for the simply or doubly connected case (except the last two
which are not mentioned in \cite{EfendievRuss}). 

Recall $f$ defined on $\Omega$ has non tangential (``n.t.'') limit
$\ell$  at $\xi\in\partial\Omega$ if and only if, for every
$0<\beta<\pi/2$, $f(z)$ tends to $\ell$ as
$\Omega\cap S_{\xi\beta}\ni z\to \xi$, where $S_{\xi\beta}$
is the cone with vertex $\xi$ and aperture
$2\beta$ whose axis is normal to $\partial\Omega$ at $\xi$.

Also, the non-tangential maximal function of $f$ at $\xi\in\T$ 
is
\begin{equation}
\label{defntmax}
{\mathcal M}_f(\xi):=\sup_{z\in \Omega\cap S_{\xi\beta}}
|f(z)|,
\end{equation}
where we dropped the dependence of ${\mathcal M}_f$ on $\beta$.

Further ({\it cf.} \eqref{defDeltat}),
we define a map $P_{\partial\Omega,\varepsilon}:
\partial\widetilde{\Delta}_\varepsilon\to\partial\Omega$ as follows
(projection on the boundary). When $\Omega$ is normalized circular,
$P_{\partial\Omega,\varepsilon}(\xi)$ is the radial projection of 
$\xi$ on the boundary circle nearest to $\xi$. 
When $\Omega$ is a general Dini-smooth domain
in $\overline{\C}$, we pick a conformal map $\psi$ onto a normalized 
circular domain $\Omega'$ and we set $P_{\partial\Omega,\varepsilon}=\psi^{-1}\circ P_{\partial\Omega',\varepsilon}\circ\psi$. Clearly, 
$P_{\partial\Omega,\varepsilon}$ is a homeomorphism.
Different $\psi$ produce different
$P_{\partial\Omega,\varepsilon}$, but the results below hold for any of them.

Assumptions being as in Theorem \ref{BS}, the following properties holds.
\begin{property} 
\label{1}
Any $f$ in $H^p_{\nu}(\Omega)$ (resp. $w\in G_\alpha^p(\Omega)$)
has a non-tangential limit almost everywhere on $\partial \Omega$, thereby
defining a trace function $\tr_{\partial \Omega} f\in L^p(\partial \Omega)$\footnote{There is no discrepancy in the notation since the nontangential limit coindes with the Sobolev trace whenever it exists
\cite[Prop. 4.3.3]{BLRR}}.
It holds that 
\begin{equation}
\label{cvLpbord}
\lim_{\varepsilon\to0}\|\tr_{\partial\Omega}f-f\circ P^{-1}_{\partial\Omega,\varepsilon}
\|_{L^p(\partial\Omega)}=0\qquad\Bigl({\rm resp.}\ 
\lim_{\varepsilon\to0}\|\tr_{\partial\Omega}w-w\circ P^{-1}_{\partial\Omega,\varepsilon}
\|_{L^p(\partial\Omega)}=0\Bigr).
\end{equation}
\end{property}
\begin{property}
\label{2}
The quantity $\| \tr_{\partial \Omega} \, . \, \|_{L^p(\partial \Omega)}$ defines an 
equivalent norm on $H^p_{\nu}(\Omega)$ (resp. $G_\alpha^p(\Omega))$. 
As to the maximal function, it 
holds when $f\in H^p_\nu(\Omega)$ (resp. $w\in G^p_\alpha(\Omega)$) that
\begin{equation}
\label{ntestimates}
\|{\mathcal M}_f\|_{L^p(\partial\Omega)}\leq C
\|\tr_{\partial\Omega}f\|_{L^p(\partial\Omega)}\qquad
\Bigl(\mbox{resp.}\ 
\|{\mathcal{M}}_{w}\|_{L^p(\partial\Omega)}\leq C
\|\tr_{\partial\Omega} {w}\|_{L^p(\partial\Omega)}\Bigr)
\end{equation}
where $C$ depends only on $\Omega$, $\sigma$ (resp. $\alpha$), $p$ and the 
aperture $\beta$ used in the definition of the maximal function.
\end{property}
{\it Proof.} In view of Proposition \ref{trick-hardy} and the H\"older-continuity of 
$\nu$, Properties \ref{1} and \ref{2} for $f$ follow from their counterpart 
for $w$. The latter are consequences of \eqref{decompexpH}, 
Remark \ref{remHolder}, and
the corresponding property in $H^p(\Omega)$, {\it cf.} Lemma \ref{IntH}.

\hfill\boite
\begin{property}
\label{3}The space $\tr_{\partial \Omega} H^p_{\nu}(\Omega)$ 
(resp. $\tr_{\partial \Omega} G^p_{\alpha}(\Omega)$)
is closed in $L^p(\partial \Omega)$.  If $f\in H_{\nu}^p(\Omega)$
(resp. $w\in G^p_{\alpha}(\Omega)$) is not identically zero, 
then $\tr_{\partial \Omega} f$ 
(resp. $\tr_{\partial \Omega} f$) cannot vanish on a subset of 
$\partial \Omega$ with positive measure.
\end{property}
{\it Proof.} As before it is enough to prove it for $G^p_\alpha$.
That  $\tr_{\partial \Omega} G^p_{\alpha}(\Omega)$
is closed in $L^p(\partial \Omega)$ follows from Property \ref{2}
and Theorem \ref{BS} point $(i)$. That $w\neq0$ cannot vanish on
a subset of  $\partial \Omega$ with positive measure is immediate from 
\eqref{decompexpH} and the corresponding result in $H^p(\Omega)$
\cite[Thm 2.2]{duren}. \hfill\boite
\begin{property}
\label{4}
If $f$ (resp.$w$) is a nonzero member of 
$H^p_\nu(\Omega)$ (resp. $G_\alpha^p(\Omega)$), then
$\log|\tr_{\partial\Omega}f|$ (resp. 
$\log|\tr_{\partial\Omega}w|$) lies in $L^1_\RR(\partial\Omega)$. Moreover the zeros of 
$f$ 
(resp. $w$) are isolated, and if we enumerate them as $\xi_j$, $j\in\NN$, it 
holds for any $z_0\in\Omega$, $z_0\neq\xi_j$ for all $j$, that
\begin{equation}
\label{Blaschkec}
\sum_{j=1}^\infty g_\Omega(\xi_j,z_0)<\infty,
\end{equation}
where $g_\Omega(.,z_0)$ is the Green function of $\Omega$ with pole at 
$z_0$\footnote{$g_\Omega(.,z_0)$ is the unique harmonic function in 
$\Omega\setminus z_0$ such that $g_\Omega(z,z_0)+\log|z-z_0|$ is bounded in
a neighborhood of $z_0$ and $g_\Omega(z,z_0)\to0$ when $z\to\xi\in\partial\Omega$, see \cite[Sec. 4.4.]{Ransford}.}. 
\end{property}
When $\Omega=\D$, 
\eqref{Blaschkec} is equivalent to the classical Blaschke condition
$\sum_j(1-|\xi_j|)<\infty$.\\

{\it Proof.} If $f$ and $w$ are related by \eqref{ftow}, their $\log$-modulus
are comparable and they share 
the same zeros. Therefore it is enough to prove the result for $w$.
That $\log|\tr_{\partial\Omega}w|\in L^1_{\mathbb{R}}(\partial\Omega)$ unless $w\equiv0$ follows from 
\eqref{decompexpH}  and the corresponding result for holomorphic functions, 
{\it cf.} Lemma \ref{IntH}.
In another connection, \eqref{decompexpH} entails that the zeros of $w$ 
are those of the holomorphic function
$F$, hence they are isolated. Moreover,
since $F\in H^p(\Omega)$, it follows from the decomposition theorem 
\cite[Sec. 10.5]{duren}  and a classical result on the disk 
\cite[Thm 5.4]{gar} that the subharmonic function
$\log|F|$ has a harmonic majorant \cite[Sec. 2.6, Ex. 10]{duren}. Relation
\eqref{Blaschkec}  now follows from \cite[Thm 4.5.5]{Ransford}.
\hfill\boite
\begin{property}
\label{5}
Each $f\in H^p_\nu(\Omega)$ satisfies the maximum principle,
{\it i.e.} $|f|$ cannot assume a relative maximum in $\Omega$ unless 
it is constant. More generally, a non-constant function in $H^p_\nu(\Omega)$
is open and the preimage of
any value is discrete.
\end{property}
{\it Proof.} 
If we let $\nu_f(z):=\nu(z) \overline{\partial f}/\partial f(z)$ if 
$\partial f(z)\neq0$ and $\nu_f(z)=0$ otherwise, then $f$
is a pointwise a.e. solution in $\Omega$ of the \emph{classical} 
Beltrami equation
\begin{equation}
\label{clBel}
\overline\partial f  = \nu_f  \, \partial f, ~~~~~~|\nu_f|<\kappa<1.
\end{equation}
Moreover since $r>2$, it holds that $W^{1,r}_{loc}(\Omega)$ is an algebra
so that $w$ given by \eqref{decompexpH} hence also $f$ given by
\eqref{ftow} lies in $W^{1,r}_{loc}(\Omega)$.
It follows by Stoilov factorization \cite[Thm  11.1.2]{im2}
that $f=G(h(z))$, where $h$ is a quasi-conformal
topological map $\Omega\to\C$ satisfying \eqref{clBel} and $G$ a 
holomorphic function on $h(\Omega)$. The conclusion now follows at once 
from the corresponding properties of holomorphic functions.
\hfill\boite\\
\begin{property}
\label{6}
To any $p_1\in[p,2p)$ there is a constant $C$ depending only on 
$\Omega$, $\nu$ (resp. $\alpha$), and $p_1$ such that, for each
$f\in H^p_\nu(\Omega)$ (resp. $w\in G_\alpha^p(\Omega)$),
\begin{equation}
\label{impsum}
\|f\|_{L^{p_1}(\Omega)}\leq C \|f\|_p\,,
\qquad
\Bigl(\mbox{resp.}\ 
\|w\|_{L^{p_1}(\Omega)}\leq C
\|w\|_p
\Bigr).
\end{equation}
Moreover, 
to each open set ${\cal O}$ with $\overline{{\cal O}}\subset\Omega$,
there is a constant $c$ 
depending only on $\Omega$, ${\cal O}$, $\nu$ (resp. $\alpha$), 
$r$, and $p$ such that,
\begin{equation}
\label{estimint1}
\|f\|_{W^{2,r}({\cal O})}\leq c\, \|f\|_p\,,
\qquad
\Bigl(\mbox{resp.}\ 
\|w\|_{W^{1,r}({\cal O})}\leq c'\,
\|w\|_p
\Bigr).
\end{equation}
\end{property}
{\it Proof.} We may assume $\Omega$ is bounded.
By Proposition \ref{trick-hardy} and Remark \ref{remHolder},
inequality  \eqref{impsum} for $f$ follows from that on 
$w$. The latter is a consequence of \eqref{decompexpH}
and the corresponding property in $H^p(\Omega)$, {\it cf.} Lemma \ref{intsum}.

Inequality \eqref{estimint1} for $\|w\|_{W^{1,r}({\cal O})}$ 
follows at once from Property \ref{2} and
the Cauchy formula as applied to $F$ in \eqref{decompexpH}. Observe
from \eqref{ftow} that a similar inequality  holds for 
$\|f\|_{W^{1,r}({\cal O})}$ .

In addition, if we pick $\varepsilon>0$ so small that  
$\overline{{\cal O}}$ lies interior to $\widetilde{\Delta}_\varepsilon$
({\it cf.} \eqref{defDeltat}) and if, in the previous argument,
we apply the Cauchy formula to $F$ on each curve 
$\partial\widetilde{\Delta}_{t}$
for $t\in[\varepsilon,\varepsilon/2]$ and then integrate with respect to $t$,
we obtain an inequality of the form 
\begin{equation}
\label{Sobsur}
\|w\|_{W^{1,r}({\cal O})}\leq C
\|w\|_{L^r(\widetilde{\Delta}_{\varepsilon/2}\setminus
\widetilde{\Delta}_\varepsilon)}.
\end{equation}
 In view of \eqref{commew}, we may apply this to
$(1-\nu^2)^{1/2}\partial\,f$ with $\Omega$ replaced by the interior of
$\widetilde{\Delta}_{\varepsilon/3}$ and, since $(1-\nu^2)^{-1/2}\in W^{1,r}(\widetilde{\Delta}_{\varepsilon/2})$ which is an algebra,
we deduce an inequality of the form
$\|\partial f\|_{W^{1,r}({\cal O})}\leq C'
\|f\|_{W^{1,r}(\widetilde{\Delta}_{\varepsilon/2})}$. As  $f$ satisfies
\eqref{dbar} a similar inequality holds for $\overline{\partial} f$,
and since 
$\|f\|_{W^{1,r}(\widetilde{\Delta}_{\varepsilon/2})}\leq C''\|f\|_p$
as pointed out already, we get that part of \eqref{estimint1} dealing
with $\|f\|_{W^{2,r}({\cal O})}$. \hfill\boite

\begin{remarque}
If we pick $\sigma\in W^{1,2}(\D)$ satisfying \eqref{ellipticsigma} 
but having no nontangential limit a.e. on 
$\T$ \cite{Carleson}, then 
$\sigma^{1/2}+i\sigma^{-1/2}$ is a solution to \eqref{eq:w} meeting the 
Hardy condition \eqref{esssuplp} for all $p<\infty$ 
(since $W^{1/2,2}(\T_r)\subset VMO(\T_r)$ \cite{brezis2}), 
 but having no nontangential limit
a.e. on $\T$. Hence the assumption that $r>2$ is
necessary for Property \ref{1} to hold.
\end{remarque}

We conclude this section with a parameterization  of $G^p_\alpha(\Omega)$ by
$H^p(\Omega)$-functions which proceeds differently from Proposition
\ref{expsf}, and is fundamental to our approach of the Dirichlet problem. 
It was essentially obtained on the disk  in
\cite{BLRR} when $\sigma\in W^{1,\infty}(\Omega)$, and then carried over to 
the annulus in \cite{EfendievRuss}
under the assumption that $r>2$ and $p>r/(r-2)$.
It features the 
operator $T_\alpha$, defined for $h\in L^p(\Omega)$ by the formula
\begin{equation}
\label{defTalpha}
{T}_\alpha h(z)=\frac1{2\pi i}\iint_{\Omega}
\frac{\alpha(\xi)\overline{h}(\xi)}{\xi-z}d\xi\wedge
d\overline{\xi}\,,\  \ z\in \Omega.\,.
\end{equation}
Note that $T_\alpha h(z)$ is well-defined for a.e. $z$ when $p,r$ satisfy 
\eqref{hypothesesnup}, since 
$\alpha\overline{h}\in L^\gamma(\Omega)$ with $1/\gamma=1/p+1/r<1$. Also, 
$T_\alpha$ is linear when $L^p(\Omega)$ is viewed as a real Banach 
space.
\begin{proposition}
\label{paramwHp}
Assume that $\Omega \subset\overline{\C} $ is a bounded Dini-smooth
domain and that $\alpha\in L^r(\Omega)$ while $p$, $r$ satisfy 
\eqref{hypothesesnup}. Then $T_\alpha$ is compact from $L^p(\Omega)$
into itself, and $I-T_\alpha$ is invertible.
It holds that $G^p_\alpha(\Omega)=(I-T_\alpha)^{-1} H^p(\Omega)$,
and if $w\in G^p_\alpha(\Omega)$ then the unique $g\in H^p$ such 
that $w=(I-T_\alpha)^{-1}g$ is the Cauchy integral of $\tr_{\partial\Omega}w$:
\begin{equation}
\label{repCauchyg}
g(z)=\mathcal{C}(\tr_{\partial\Omega}w)=\frac{1}{2i\pi}\int_{\partial\Omega}
\frac{\tr_{\partial\Omega}w(\xi)}{\xi-z}\,d\xi, \qquad z\in\Omega.
\end{equation}
Moreover, we have that $\|w\|_p\leq C\|g\|_p$ where the constant $C$ depends
only of $\Omega$, $\alpha$, and $p$.
\end{proposition}
The proof of Proposition \ref{paramwHp} is given in Appendix \ref{app_proof_prop_3} and 
follows the lines of the proof of  \cite[Thm 4.4.1.1]{BLRR}, although 
technicalities arise  to handle the weaker assumption \eqref{hypothesesnup}.

\section{A decomposition theorem} 
\label{secGA}

Let $\Omega\subset\overline{\C}$ be Dini-smooth and write
$\Omega=\overline{\C}\setminus\cup_{j=0}^n K_j$ where the $K_j$ are disjoint 
compact sets in $\overline{\C}$.
We establish in this section a result which, loosely speaking, asserts 
that every  function in $H^p_\nu(\Omega)$  is a sum of 
members of $H^p_{\nu_j}(\overline{\C}\setminus K_j)$
for some appropriate extensions $\nu_j$  of $\nu$.
This result stands analogous to the decomposition theorem in $H^p(\Omega)$
\cite[Sec. 10.5]{duren} but, unlike Properties \ref{1}-\ref{4} in Section 
\ref{sechpnud}, it is not an immediate consequence of the latter 
{\it via}  Proposition \ref{expsf}.

We denote by  $H^{p,0}_{\nu}(\Omega)$ the subspace  of $H^p_{\nu}(\Omega)$ 
made of functions $f$  such that 
\begin{equation}
 \label{normim}
\int_{\partial\Omega} \mbox{Im}(\tr_{\partial \Omega} f(s))\,|ds| = 0.
\end{equation}
Moreover, we let $H^{p,00}_{\nu}(\Omega)$ be the subspace  of  
$H^{p,0}_{\nu}(\Omega)$ consisting of those $f$ for which
\begin{equation}
 \label{normal:fv}
\int_{\partial\Omega}\tr_{\partial \Omega} f(s) \ |ds| = 0.
\end{equation}
We record for later use a topological version of the decomposition theorem for
holomorphic Hardy spaces. For normalized cicular domains it is established 
{\it e.g.} in \cite[Lem. 2.1]{cp}\footnote{In this reference, Hardy spaces are defined through harmonic majorants, but we know 
this is equivalent to the definition based on \eqref{systema} for Dini-smooth domains.}, and 
it carries over immediately to Dini-smooth domains by conformal invariance.
\begin{lemma}
\label{rmq:Hardy-ann}
If $\Omega=\overline{\C}\setminus\cup_{j=0}^n K_j$ is a Dini-smooth domain
as above,
then additively 
\[H^p(\Omega) =  H^p(\overline{\C}\setminus K_0) \oplus {H^{p,00}}
(\overline{\C} \setminus K_1)\oplus\cdots\oplus {H^{p,00}}
(\overline{\C} \setminus K_n),
\] 
where the direct sum is topological. 
\end{lemma}

Recalling from \eqref{defcheck} the notation
$\check \nu$, we also have the following lemma which is established  in
\cite{flps} for $\nu \in W^{1,\infty}_\RR(\DD)$.
\begin{lemma}
\label{cor-oplus}
 Assume that $\Omega$ is a proper Dini-smooth
domain in $\overline{\C}$
and that $p$, $r$, $\nu$ satisfy \eqref{hypothesesnup}. 
Then the following topological decomposition holds:
\[L^p(\T) = \tr_\T H_{\nu}^p(\DD) \oplus \tr_\T {H^{p,00}_{ \check{\nu}}}(\overline{\CC} \setminus \overline{\DD}).\]
\end{lemma}
{\it Proof.} The proof of \cite[Cor. 3]{flps} applies without change to
$\nu\in W^{1,r}_\RR(\Omega)$, granted Lemma 
\ref{lemu} in Section \ref{secdir} to come. \hfill\boite\\

The main result in this section is the following generalization of 
Lemma \ref{rmq:Hardy-ann}.
\begin{theorem}
\label{rmq:Hardy-ann2}
Let $\Omega=\overline{\C}\setminus\cup_{j=0}^n K_j$ be a Dini-smooth domain,
the $K_j$ being disjoint compact sets in $\overline{\C}$,
and $p$, $r$, $\nu$ meet \eqref{hypothesesnup}. 
Then to each $j$, there is $\nu_j\in W_\RR^{1,r}(\overline{\C}\setminus K_j)$
with ${\nu_j}_{|\Omega}=\nu$ and ${\nu_j}_{|K_l}={\nu_k}_{|K_l}$ when
$l\neq j,k$, satisfying (\ref{kappa}) and such that
\begin{equation}
\label{decompHpnua}
H^p_\nu(\Omega) =  
H^p_{\nu_0}(\overline{\C}\setminus K_0) \oplus {H_{\nu_1}^{p,00}}
(\overline{\C} \setminus K_1)\oplus\cdots\oplus {H_{\nu_n}^{p,00}}
(\overline{\C} \setminus K_n),
\end{equation} 
where the direct sum is topological. 
\end{theorem}

{\it Proof.} It is clear that the right side of \eqref{decompHpnua} is
included in the left. To prove the converse, we may assume by conformal 
invariance (see Lemma \ref{conf}) that $\Omega$ is normalized circular.

 First, we consider the special case where $\Omega=\A_\varrho$ ($n=1$) for some 
$0<\varrho<1$. In this case, we put $\nu_i$ and $\nu_e$ for the seeked 
extensions of $\nu$ to $\D$ and $\overline{\C}\setminus\overline{\D_\varrho}$
(the subscripts respectively stand for ``interior'' and ``exterior'').
It is standard that there exists $\widetilde{\nu}\in W_\R^{1,r}(\D)$ such that 
$\widetilde{\nu}_{|\A_\varrho}=\nu$ \cite[Sec.  VI.3, Thm 5]{stein}.
Letting $\varepsilon>0$ be so small that 
$\|\nu\|_{L^\infty(\A_\varrho)}<\kappa-\varepsilon$, we set 
$\nu_i=\min(\widetilde{\nu},\kappa-\varepsilon)$, so that 
$\nu_i$ lies in $W^{1,r}(\D)$, extends $\nu$, and meets \eqref{kappa}.

Let $f \in H^p_{\nu}(\A_\varrho)$ so that $\tr_\T f \in L^p(\T)$. 
Lemma \ref{cor-oplus} is to the effect that
\begin{equation}
\label{eq:tr}
\tr_\T f = \tr_\T f_i + \tr_\T f_e \, ,
\end{equation}
for some $f_i \in H_{\nu_i}^p(\DD)$ and $f_e \in {H^{p,00}_{\check{\nu_i}}}(\overline{\CC} \setminus \overline{ \DD})$. 
Put $\mathcal{F}_i = f - f_{i|_{\A_\varrho}} \in H_{\nu}^p(\A_\varrho)$, 
and let $\mathcal{F}_e$ and $\nu_e$ be defined on 
$\overline{\CC} \setminus \varrho \overline{\DD}$ by:
\[
 \mathcal{F}_e  = \mathcal{F}_i \vee f_e = \left\{\begin{array}{l} 
\mathcal{F}_i \mbox{ on } \AA_\varrho \, ,\\
f_e \mbox{ on } \overline{\CC} \setminus  \overline{\DD} \, ,
\end{array}
\right.
\ \ \
\nu_e  = \nu \vee \check{\nu}_i
= \left\{\begin{array}{l}
\nu \mbox{ on } \AA_\varrho \, ,\\
\check{\nu}_i \mbox{ on } \overline{\CC} \setminus  \overline{\DD} \, .
\end{array}
\right.
\]
Since $\nu={\check \nu}_i$ on $\T$, it holds that 
$\nu_e \in W_\RR^{1, r}(\overline{\CC} \setminus \varrho \overline{\DD})$ 
(by absolutely continuity on lines in polar coordinates) and obviously 
it satisfies a condition similar to (\ref{kappa}).

{\it We claim} that $\mathcal{F}_e \in {H^p_{{\nu}_e}}(\overline{\C} \setminus \varrho \overline{\DD})$. 
By construction  it satisfies (\ref{systema}) on 
$\overline{\CC} \setminus \varrho \overline{\DD}$ and 
it is a solution to \eqref{dbar} (with $\nu_e$ instead of $\nu$) on
$\AA_\varrho \cup ({\overline \CC}
\setminus \overline{\DD})$. Thus, in order to establish that \eqref{dbar}
holds on 
the whole of $\overline{\CC} \setminus \varrho \overline{\DD}$, it is enough to prove that $\bar\d \mathcal{F}_e = \nu_e \ \overline{\partial \mathcal{F}_e}$, in the sense of distributions, on some annulus $\A_{r,R}$ 
with $\varrho < r < 1<R$. That is, we must show that
for all $\phi \in {\mathcal{D}}(\A_{r,R})$,
\begin{eqnarray}
\label{calFsol}
I_{r,R} (\phi) & = & \iint_{\A_{r,R}} (- \mathcal{F}_e \, \overline{\partial} \phi + \overline{\mathcal{F}_e} \, \overline{\partial} (\nu_e \, \phi)) \ dm(z) \nonumber \\
& = & \dfrac{i}{2} \iint_{\A_{r,R}} (- \mathcal{F}_e \, \overline{\partial} \phi + \overline{\mathcal{F}_e} \, \overline{\partial} (\nu_e \, \phi)) \ dz \wedge d\overline{z} = 0.
\end{eqnarray}
By Property \ref{6} in Section \ref{sechpnud}, we get that 
$f_i\in L^{p_1}(\D)$, $f_e\in L^{p_1}(\overline{\C}\setminus\overline{\D})$, 
and $f\in L^{p_1}(\A_{\varrho})$ hence ${\mathcal F}_i\in L^{p_1}(\A_\varrho)$
and ${\mathcal F}_e\in L^{p_1}(\overline{\C}\setminus\varrho\overline{\D})$, 
for some $p_1 >2$. Thus, in view of 
\eqref{hypothesesnup} and  H\"older's inequality,
the integrand in \eqref{calFsol} lies in $L^a(\A_{r,R})$ for some $a>1$ 
which justifies the limiting relation: 
\begin{eqnarray*}
-2i I_{r,R} (\phi) & = & \iint_{\A_r} (- \cF_i \, \overline{\partial} \phi + \overline{\cF_i} \, \overline{\partial} (\nu \, \phi)) \ dz \wedge d\overline{z}
+\iint_{\A_{1,R}} (- f_e \, \overline{\partial} \phi + \overline{f_e} \, \overline{\partial} (\check{\nu}_i \, \phi)) \ dz \wedge d\overline{z} \\
& = & \lim_{\epsilon \to 0} \underbrace{\iint_{\A_{r,1-\epsilon}} (- \cF_i \, \overline{\partial} \phi + \overline{\cF_i} \, \overline{\partial} (\nu \, \phi)) \ dz \wedge d\overline{z}}_{I_{r,\epsilon}(\phi)} \\
&&~~~~~~~~~~~~~~~~~~~~~~~~~~~~~~~+~\lim_{\epsilon \to 0} \underbrace{\iint_{\A_{1+\epsilon,R}} (- f_e \, \overline{\partial} \phi + \overline{f_e} \, \overline{\partial} (\check{\nu}_i \, \phi)) \ dz \wedge d\overline{z}}_{I_{R,\epsilon}(\phi)}.
\end{eqnarray*}

As ${\check\nu}_i\in W^{1,r}(\overline{\C}\setminus\overline{\D})$ and
$f_e \in W^{1,r}_{loc} (\overline{\C}\setminus\overline{\D})$ 
by Proposition \ref{expsf}, we can apply Stoke's theorem:
\begin{eqnarray*}
I_{R,\epsilon} (\phi) = 
<\overline{\partial} f_e - \check{\nu}_i \overline{\partial}  \overline{f_e},\phi>_{\A_{1+\epsilon,R}} + \int_{\partial \A_{1+\epsilon,R}} (f_e - \check{\nu}_i \overline{f_e}) \phi \ dz.
\end{eqnarray*}
Since $L^a(\overline{\C}\setminus\overline{\D})\ni\overline{\partial} f_e - \check{\nu}_i \overline{\partial}  \overline{f_e}=0$ a.e., we are left with
\begin{eqnarray*}
 I_{R,\epsilon} (\phi) & = & \int_{\partial \A_{1+\epsilon,R}} (f_e - \check{\nu}_i \overline{f_e}) \phi \ dz =  - \int_{\T_{1 + \epsilon}} (f_e - \check{\nu}_i \overline{f_e}) \phi \ dz
\end{eqnarray*}
since $\phi$ vanishes on $\TT_R$. 
Passing to the limit using Property \ref{1} in Section \ref{sechpnud} and the 
continuity of $\check{\nu}_i$ in $\overline{\C}\setminus\overline{\D}$ 
yields
\[
\lim_{\epsilon \to 0} 
I_{R,\epsilon} (\phi) = - \int_{\TT} (f_e - \check{\nu}_i \overline{f_e}) \phi \ dz.
\]
Likewise one can show that
\[
\lim_{\epsilon \to 0} 
I_{r,\epsilon} (\phi) = \int_{\TT} (\mathcal{F}_i - \nu \overline{\mathcal{F}_i}) \phi \ dz,
\]
and since $\tr_\T \cF_i = \tr_\T f_e$ by (\ref{eq:tr}) while 
$\nu_{|_\T} = {\nu}_{i|_\T} = \check{\nu}_{i|_\T}$, we finally conclude that
\[
-2i I_{r,R}(\phi) = \lim_{\epsilon \to 0} I_{\Omega_{r,\epsilon}} (\phi) +\lim_{\epsilon \to 0} I_{\Omega_{R,\epsilon}} (\phi) 
=0,
\]
which is  (\ref{calFsol}). This proves the claim that 
$\mathcal{F}_e \in {H^p_{{\nu}_e}}(\overline{\CC} \setminus \varrho \overline{\DD})$. 
Now, 
\begin{equation}
\label{decompf}
f = \Bigl( f_{i|_{\AA\varrho}}+\int_{\T_\varrho} \tr_{\T_\varrho} \mathcal{F}_e\Bigr)  + \Bigl(\mathcal{F}_{e|_{\AA_\varrho}}- \int_{\TT_\varrho} \tr_{\T_\varrho} \mathcal{F}_e\Bigr)
\end{equation}
is the decomposition we look for on $\A_\varrho$.

The (not yet topological) existence of  \eqref{decompHpnua} on 
$\Omega=\A_\varrho$
implies its existence on any Dini-smooth doubly connected domain by 
conformal invariance. 

Subsequently,
we get it over any normalized circular domain by induction on $n$: if
$\Omega=\D\setminus\cup_{j=1}^n{\overline{\D}_{a_j,r_j}}$ and $\varrho$
is close enough to 1 that 
$\D_\varrho\supset\cup_{j=1}^n{\overline{\D}_{a_j,r_j}}$,  we decompose
$f_{|\A_\varrho}=f_{0|\A_\varrho}+f_{e|\A_\varrho}$ with $f_0\in H^p_{\nu_0}(\D)$ and
$f_e\in H^p_{\nu_e}(\overline{\C}\setminus\varrho\overline{\D})$ 
for suitable extensions $\nu_0$, $\nu_e$ of $\nu$ to $\D$ and 
$\overline{\C}\setminus\varrho\overline{\D}$ respectively. As it coincides with
$f_e$ on $\A_\varrho$, the function $f-f_0$ 
lies in $H^p_{\nu_0\vee\nu_e}(\Omega')$ where
$\Omega'=\overline{\C}\setminus\cup_{j=1}^n{\overline{\D}_{a_j,r_j}}$ 
is $n-1$-connected, hence we can carry out the induction step.
Observe that in the latter  the $\nu_j$ will coincide with
$\nu_e$ on $\overline{\C}\setminus\varrho\overline{\D}$ for $1\leq j\leq n$, 
thereby proving the 
existence of decomposition \eqref{decompHpnua} in general.

To see that the sum is direct, write 
\[f=\sum_{j=0}^n f_j\qquad {\rm with}\ \ f_0\in H^p_{\nu_0}(\D)\ \ {\rm and}\ 
\  f_j\in H^{p,00}_{\nu_j}(\overline{\C}\setminus\overline{\D}_{a_j,r_j}),
\ 1\leq j\leq n.\]
Suppose that $f_0=-\sum_{j=1}^n f_j$ on $\Omega$. Then 
$h=f_0\vee(-\sum_{j=1}^n f_j)_{|\overline{\C}\setminus\overline{\D}}$
lies in $W^{1,p}(\overline{\C})$ and satisfies
$\overline\partial h  = \widetilde{\nu}  \, 
\overline{\partial h}$
with $\widetilde{\nu}=\nu_0\vee\nu_1\vee\cdots\vee\nu_n
\in W^{1,r}(\overline{\C})$ (remember two $\nu_j$ coincide wherever they are 
both defined). By the arguments in Remark  \ref{rmq:xcte}, we deduce that
$h$, thus also $f_0$ is a constant, say $C_0$. If we put $\widetilde{f}_1=f_1+C_0$,
then  $\widetilde{f}_1=-\sum_{j=2}^n f_j$ and arguing the same way we find
that $\widetilde{f}_1$, thus also $f_1$, is in turn a constant. 
Proceeding inductively
each $f_j$ is a constant $C_j$, and $C_j=0$ for $j\geq1$ by 
\eqref{normal:fv}. Therefore $f_j=0$ for $j\geq1$ and then $f_0=0$ as well.

Finally, having shown that the natural map
\[
H^p_{\nu_0}(\overline{\C}\setminus K_0) \oplus {H_\nu^{p,00}}
(\overline{\C} \setminus K_1)\oplus\cdots\oplus {H_{\nu_1}^{p,00}}
(\overline{\C} \setminus K_n)\longrightarrow H^p_{\nu_n}(\Omega) 
\] 
is injective and surjective, we observe from Property 2 in Section 
\ref{sechpnud} that it is continuous, hence a homeomorphism by the open mapping theorem.

\hfill\boite\\

\section{The Dirichlet problem}
\label{secdir}
For $\Omega$ a Dini-smooth domain,
we let  ${\cal U}^p(\Omega)$ consist of those functions $U$  
satisfying (\ref{div}) for which $\|U\|_p<\infty$ ({\it cf.} \eqref{systemn}).
In this section, we investigate the solvability of the Dirichlet
problem for  the class $\mathcal{U}^p(\Omega)$  with  boundary data in
$L^p(\partial\Omega)$. In other words,
the solution is  understood to meet condition \eqref{systema} for some 
admissible sequence 
$\Delta_n$, and to converge non-tangentially on $\partial\Omega$  to 
some prescribed member of $L^p(\partial\Omega)$.

The existence of non-tangential estimates of the form 
\eqref{ntestimates} for functions in $\mathcal{U}^p(\Omega)$ 
will make these requirements  equivalent, in the 
present context, to standard notions of solvability  \cite{FJR}. 
Note that $\sigma\in W^{1,r}(\Omega)$ 
with $p$, $r$ as in \eqref{hypothesesnup} 
is an assumption which is not covered by the Carleson 
condition set up in
\cite{DPP,KP}\footnote{The condition is that the {\it sup} on the Carleson 
domain 
$B(z,d(z,\partial\Omega)/2)$ of 
$d(.,\partial\Omega)|\nabla\sigma(.)|^2$, when viewed as a function of $z$,
 should be the density of a
Carleson (or vanishing  Carleson) measure.  Now, for $\chi_k$  the 
caracteristif function of $1-1/k-e^{-k^2}< |z|<1-1/k+e^{-k^2}$, $k\geq2$,
the radial function $\psi(r,\theta)=
\sum_{k=2}^\infty e^{k}\chi_k$ lies in $L^r(\D)$
for each $r\in(2,\infty)$. If we put 
$\sigma(r,\theta)=1+\int_0^r\psi(\rho)\,
d\rho$, the corresponding density is not even integrable on 
Carleson domains.}.


We study the Dirichlet problem for \eqref{div} in relation to the issue
of finding a function in $H^p_\nu(\Omega)$ with prescribed real part on 
$\partial\Omega$, $\nu$ being as in \eqref{dfnsigma}. 
Slightly abusing terminology, we call this issue  
the Dirichlet problem in $H^p_\nu(\Omega)$. Clearly, a solution to the 
Dirichlet problem for \eqref{div} is obtained from a solution to the 
Dirichlet problem in $H^p_\nu(\Omega)$ by taking the real part. 
However, we shall see that the Dirichlet problem in $H^p_\nu(\Omega)$ is not
always solvable  on multiply connected domains, whereas the
Dirichlet problem for \eqref{div} is  solvable.


On simply connected domains the two problems are equivalent as follows from 
the next lemma. When
$\sigma\in W_\R^{1,\infty}(\Omega)$,  
this is essentially proved in \cite{BLRR}
except for estimate \eqref{usobloc}.

\begin{lemma} 
\label{lemu}
Let $\Omega$ be a Dini-smooth simply connected domain and
$\sigma\in W^{1,r}(\Omega)$ satisfy \eqref{ellipticsigma} with 
$p$, $r$ as in  \eqref{hypothesesnup}. 
For every $u \in L^p_\RR(\partial\Omega)$, the following assertions hold.

$(i)$  
There exists a unique solution $U$ to (\ref{div}) in $\Omega$ such 
that $\|U\|_p<\infty$ ({\it cf.} \eqref{systemn}) and the non-tangential limit
of $U$ on $\partial\Omega$ is $u$. The function $U$ satisfies non-tangential 
estimates of the form \eqref{ntestimates}. Moreover,
to every open set $O$ with compact closure ${\overline{O}}\subset\Omega$, 
there is a constant 
$c = c(\Omega,O,\sigma,r,p)> 0$ such that:
\begin{equation}
\label{usobloc}
\|U\|_{W^{2,r}(O)} \leq c \, \|u\|_{L^{p}(\partial \Omega)} \, .
\end{equation}
$(ii)$ For $\nu$  defined by \eqref{dfnsigma}, there exists a unique 
$f \in H_{\nu}^{p,0}(\Omega)$ 
such that $\mbox{Re } \tr_{ \partial\Omega} f = u$. Moreover
$\|\tr_{\partial\Omega} f\|_{L^{p}(\partial\Omega)}\leq c\|u\|_{L^{p}(\partial\Omega)}$ for some constant $c=c(\Omega,p, \sigma) > 0$.
\end{lemma}
We prove Lemma \ref{lemu}  in Appendix \ref{app5}.


\subsection{Conjugate functions
}
\label{SEC_DIR:Hp}

A $\sigma$-harmonic
conjugate to $U\in{\cal U}^p(\Omega)$ is a real function $V$ 
such that $U+iV\in H^p_\nu(\Omega)$. If it exists, a $\sigma$-harmonic 
conjugate  is
unique up to an additive constant and is 
a solution to \eqref{system2}, see \eqref{system}.

We also say that $u\in L^p_\R(\partial\Omega)$ has
$\sigma$-harmonic conjugate $v\in L^p_\R(\partial\Omega)$ if
$u +iv =\tr_{\partial\Omega}f$ for some $f\in H^p_\nu(\Omega)$.
If $U\in{\cal U}^p(\Omega)$ has $\sigma$-harmonic conjugate $V$, then clearly
$\tr_{\partial\Omega} U$ has $\sigma$-harmonic conjugate $\tr_{\partial\Omega} V$.
Theorem \ref{Diru} further below asserts that 
each $u \in L^p_\R(\partial\Omega)$ is uniquely the trace of some 
$U\in{\cal U}^p(\Omega)$, so  the two notions of conjugacy
(in the domain and on the boundary)  will soon be proven equivalent.

Lemma \ref{lemu} entails that
if $\Omega$ is simply connected, then each $U\in{\cal U}^p(\Omega)$ 
(resp. $u \in L^p(\partial\Omega)$) has a
$\sigma$-harmonic conjugate. If $\Omega$ is multiply connected it is not so
as we will now see.
\begin{lemma}
\label{periodes}
Let $\Omega$ be a Dini-smooth domain and
$\sigma\in W^{1,r}(\Omega)$ satisfy \eqref{ellipticsigma} with
$p$, $r$ as in  \eqref{hypothesesnup}.
Each $U\in\mathcal{U}^p(\Omega)$ lies in $W^{2,r}_{loc}(\Omega)$,
and if $\Gamma\subset\Omega$  is a rectifiable Jordan curve then 
$\int_\Gamma\sigma\partial_n U$ depends only on the homotopy class of 
$\Gamma$. Moreover,
the function $U$ has a $\sigma$-harmonic conjugate if and only if
this integral is zero for all $\Gamma$.
\end{lemma}
{\it Proof.}
By Lemma \ref{conf} we may assume that 
$\Omega=\D\setminus\cup_{j=1}^n\overline{\D}_{a_j,r_j}$ is normalized circular.
Since $\mathcal{U}^p(\Omega)\subset L^p(\Omega)$, 
the proofs of Lemma \ref{lemu} and Corollary \ref{regint}
imply that  $U$ has a $\sigma$-harmonic conjugate in $W^{2,r}_{loc}(\Omega_1)$
on every Dini-smooth simply connected relatively compact subdomain 
 $\Omega_1\subset\Omega$. In particular $U\in W^{2,r}_{loc}(\Omega)$
so that $\nabla U$ is continuous, hence
$\int_\Gamma\sigma\partial_n U$ is well-defined.
For small $\varepsilon$, the circles 
$\T_{a_j,r_j+\varepsilon}$ for $1\leq j\leq n$
form a homotopy basis in $\Omega$. Hence in order to prove that
$\int_\Gamma\sigma\partial_n U$ depends only on the homotopy class of 
$\Gamma$,  it is enough to assume $\Gamma$  is homotopic to 
$\T_{a_j,r_j+\varepsilon}$ 
and to show that $\int_\Gamma\sigma\partial_n U=
\int_{\T_{a_j,r_j+\varepsilon}}\sigma\partial_n U$.
If $\Gamma$ is smooth and disjoint from $\T_{a_j,r_j+\varepsilon}$, or
intersects it transversally, the result follows immediately from the 
Green formula. The genericity of transversal intersections 
\cite{GuilleminPollack} now implies the result for every smooth $\Gamma$ 
by continuity.
If $\Gamma$ is merely rectifiable with parametrization $\gamma\in W^{1,1}(\T)$,we approximate the latter in the Sobolev sense 
by a smooth function on $\T$, which gives us the result.

Now, if $U$ has $\sigma$-harmonic conjugate $V$, then:
by \eqref{system} $\int_\Gamma\sigma\partial_n U=\int_\Gamma dV=0$ 
for every $\Gamma$.
Conversely, if $\int_\Gamma\sigma\partial_n U=0$ for every $\Gamma$,
then by continuation along any
path we can define globally in $\Omega$ a real-valued function $V$
such that $f=U+iV\in W^{2,r}_{\RR,loc}(\Omega)$ satisfies \eqref{dbar}. 
We have to prove that $f\in H^p_\nu(\Omega)$. 
By Proposition \ref{expsf}, it is equivalent to show that if
$w$ given by \eqref{ftow} gets factored as in \eqref{decompexpH},
then $F\in H^p(\Omega)$. The question localizes around each component of 
$\partial\Omega$, {\it i.e.} it is enough to establish for all $1\leq j\leq n$
that $\|F\|_{L^p(\T_{a_j,r_j+\eta})}$ is bounded independently of $\eta$ when 
the latter is small enough, and that
$\|F\|_{L^p(\T_{1-\eta})}$ is likewise bounded independently of $\eta$.
Consider this last case, the others being similar.
Since $U\in\mathcal{U}^p(\Omega)$ we know that
$\|\mbox{Re}\,f\|_{L^p(\T_{1-\eta})}$ is bounded independently of $\eta$, 
so by \eqref{ftow} the same is true of $\|\mbox{Re}\,w\|_{L^p(\T_{1-\eta})}$.
As
$\|\mbox{Im}\,F\|_{L^p(\T_{1-\eta})}\leq C_1+C_2 
\|\mbox{Re}\,F\|_{L^p(\T_{1-\eta})}$ by Lemma \ref{Rieszanneau},
we can argue as in the proof of Lemma \ref{lemu} (see \eqref{estimpu} and 
after)
to the effect that $\|F\|_{L^p(\T_{1-\eta})}$ is bounded independently of
$\eta$, as desired.
\hfill\boite\\

Next, we single out a particular subspace of $L^p(\partial\Omega)$ no
element of which 
has a conjugate except the zero function, and whose ``periods'' on 
a homotopy basis can be assigned arbitrarily. Namely,
if $\partial\Omega=\cup_{j=0}^n\Gamma_j$
where the $\Gamma_j$ are disjoint Jordan curves, we set 
\[
 \mathcal{S}_{\Omega} = \left\{(u_0, u_1,\cdots,u_n) \in \Pi_{j=0}^n
L^p_\R(\Gamma_j); \ u_j \equiv C_j\in\R~~\mbox{with}~\Sigma_{j=0}^n C_j=0\right\} \subset L^p_\R(\partial \Omega).
\]
\begin{proposition}
 \label{propS}
Let $\Omega$ be Dini-smooth and
$\sigma\in W^{1,r}(\Omega)$ satisfy \eqref{ellipticsigma} with
$p$, $r$ as in  \eqref{hypothesesnup}.
If $u \in \mathcal{S}_{\Omega}$  has a $\sigma$-harmonic
conjugate, then $u \equiv 0$. Each $u \in \mathcal{S}_{\Omega}$ is
uniquely the trace on $\partial\Omega$ of some $U\in W^{2,r}_\R(\Omega)$
satisfying \eqref{div}, and $\|U\|_{W^{2,r}(\Omega)}\leq C\|u\|_{L^p(\partial\Omega)}$ for some constant $C$ independent of $u$. 
If $\lambda_1,\cdots,\lambda_n$ are real numbers and
$\gamma_1,\cdots,\gamma_n$ is a homotopy basis for $\Omega$, 
then there exists a unique
$u\in \mathcal{S}_{\Omega}$ such that
$\int_{\gamma_j}\partial_n U=\lambda_j$.
\end{proposition}

{\it Proof.} 
By Lemma \ref{conf} we may assume that $\Omega = \mathbb{D} \setminus 
\cup_{j=1}^n \overline{\mathbb{D}}_{a_j,r_j}$ is normalized circular.
Set by convention $a_0=0$, $r_0=1$, so that $\mathbb{T} = \T_{a_0,r_0}$. 
Assume that $u \in {\cal S}_{\Omega}$ has a $\sigma$-harmonic
conjugate $v \in L^p_\R(\partial \Omega)$, {\it i.e.} 
$u+iv = \mbox{tr}_{\partial\Omega}f $ for some 
$f \in H^p_\nu(\Omega)$.
Denote by $\check{f}_j,\check{\nu}_j$ the reflections of $f$, $\nu$ across 
$\T_{a_j,r_j}$, $0 \leq j \leq n$,  that is
\begin{equation}
\label{defcheck2}
\check{f}_j(z-a_j) =2C_j-\overline{f\left(\tfrac{r_j^2}{\overline{z-a_j}}\right)} \, , \ 
\check{\nu}_j(z-a_j) = {\nu\left(\tfrac{r_j^2}{\overline{z-a_j}}\right)} \,  .
\end{equation}
By (\ref{equiv}), it holds that $\check{f_j} \in 
H^p_{\check{\nu}_j}(\check{\Omega}_j)$ where, for  $1\leq j\leq n$, we put 
$\check{\Omega}_j \subset \D_{a_j,r_j}$ 
(resp. $\check{\Omega}_0 \subset \overline{\C} \setminus \overline{\D}$) for 
the reflection of $\Omega$ across $\T_{a_j,r_j}$ (resp. across $\T$).
Put $\widetilde{f}=f\vee \check{f}_0 \vee \check{f}_1\cdots\vee \check{f}_n$,
$\widetilde{\nu}=\nu\vee \check{\nu}_0 \vee 
\check{\nu}_1\cdots\vee \check{\nu}_n$, and
$\check{\Omega}=
\overline{\Omega}\cup\check{\Omega}_0\cdots\cup\check{\Omega}_n\dag)$. Since
$\tr_{\T_{a_j,r_j}}\check{f}_j=\tr_{\T_{a_j,r_j}}f$ and
$\tr_{\T_{a_j,r_j}}\check{\nu}_j=\tr_{\T_{a_j,r_j}}\nu$, 
the argument leading to \eqref{calFsol}
in the proof of Theorem \ref{rmq:Hardy-ann2}
shows that
$\widetilde{f}\in H^p_{\check{\nu}}(\check{\Omega})$.
Thus $f=U+iV\in W^{2,r}(\Omega)$ by Property \ref{6}, 
hence $\sigma^{-1} \nabla V \in W^{1,r}_\mathbb{R}(\Omega)$ 
since the latter is an algebra for $r>2$. Moreover 
$\partial_n V \in W^{1-1/r,r}(\partial \Omega) \subset L^{r}(\partial \Omega)$,
which grants us enough smoothness, in view of \eqref{hypothesesnup},
to apply the divergence formula:
\begin{equation}
\label{eq:preuveS}
 0 = \iint_\Omega V \, \nabla \cdot (\sigma^{-1} \nabla V) \ dm = \int_{\partial \Omega} \sigma^{-1} \, v \, \partial_n V \ |d s|  - \iint_\Omega \sigma^{-1} |\nabla V|^2 \ dm \, ,
\end{equation}
where the first equality comes from \eqref{system2}.
Since  $u$ takes constant values a.e. on each $\T_{a_j,r_j}$,
$0 \leq j \leq n$, we have that $\partial_t u = 0$  on $\partial \Omega$ 
whence $\partial_n V = 0$ by 
(\ref{system}). Taking this into account in (\ref{eq:preuveS}), we obtain
\[
 \iint_\Omega \sigma^{-1} |\nabla V|^2 \ dm = 0 \, 
\]
implying by (\ref{ellipticsigma}) that $\nabla V =0$ a.e. in $\Omega$, 
hence $\nabla U=0$ by \eqref{system}. Thus $U$ is constant, in particular
all $C_j$ are equal, and since they add up to zero 
we obtain $u\equiv0$, as announced.

In another connection, since $\sigma\in W^{1,r}(\Omega)\subset VMO(\Omega)$,
it follows from elliptic regularity theory 
\cite{aq,Fazio} that any $u\in W^{1-1/l,l}(\partial\Omega)$ for some 
$l\in(1,\infty)$ is uniquely 
the trace of some $U\in W^{1,l}(\Omega)$ meeting \eqref{div}
with $\|U\|_{W^{1,l}(\Omega)}\leq
C(l,\Omega,\sigma)\|u\|_{W^{1-1/l,l}(\partial\Omega)}$. In particular,
since $u\in \mathcal{S}_\Omega$ is constant on each component of 
$\partial\Omega$, we have that  
\begin{equation}
\label{ellipticreg}
\|U\|_{W^{1,p}(\Omega)}\leq
C(p,\Omega,\sigma)\|u\|_{W^{1-1/p,p}(\partial\Omega)}=
C(p,\Omega,\sigma)\|u\|_{L^p(\partial\Omega)}.
\end{equation}
Moreover, if we let $C_j$ (resp. $C_0$) be the
constant value that $u$ assumes on $\T_{a_j,r_j}$ (resp. $\T$),  
the reflected function $\check{U}_j$  defined {\it via}
\eqref{defcheck2} (with $U$ in place of $f$ and $\sigma$ in place of $\nu$) 
allows us to define
$\check{U}=U\vee \check{U}_0\cdots \vee\check{U}_n\in W^{1,p}(\check{\Omega})$ 
meeting \eqref{div} in
$\check{\Omega}$ with conductivity
$\check{\nu}=\sigma \vee \check{\sigma}_0\vee\cdots\vee\check{\sigma}_n$.
Let $\Omega_1$ be a bounded Dini-smooth open set,
$\overline{\Omega}\subset\Omega_1\subset\overline{\Omega}_1\subset\check{\Omega}$. By compactness, we
can  cover $\overline{\Omega}$
with finitely many disks $D_k=\D_{\xi_k,\rho_k}$ such that
$\overline{\D}_{\xi_k,2\rho_k}\subset\Omega_1$. Using \eqref{usobloc}
and the trace theorem, we obtain
\[\|U\|_{W^{2,r}(D_k)}\leq c_k\|\tr_{\T_{\xi_k,2\rho_k}}U\|_{
L^p(\T_{\xi_k,2\rho_k})}\leq c_k\|\tr_{\T_{\xi_k,2\rho_k}}U\|_{
W^{1-1/p,p}(\T_{\xi_k,2\rho_k})}\leq 
c'_k\|U\|_{W^{1,p}(\D_{\xi_k,2\rho_k})}\]
so that $\|U\|_{W^{2,r}(\Omega)}\leq C\|U\|_{W^{1,p}(\Omega_1)}$. Moreover,
from the very form of \eqref{defcheck2}, it is easy to check that
$\|U\|_{W^{1,p}(\Omega_1)}\leq C\|U\|_{W^{1,p}(\Omega)}$. Therefore  by
\eqref{ellipticreg} $\|U\|_{W^{2,r}(\Omega)}\leq C\|u\|_{L^p(\Omega)}$,
as desired.

For $0\leq j\leq n$, let 
$\upsilon_j$ be equal to 1 on $\TT_{a_j,r_j}$
and to 0 on $\TT_{a_k,r_k}$, $k\neq j$. Each $u \in\mathcal{S}_\Omega$ 
decomposes uniquely as $u=\sum_j C_j\upsilon_j$ with $\sum_j C_j=0$. 
In addition, if 
$\Upsilon_j\in W^{2,r}(\Omega)$ is the solution to \eqref{div} such that
$\tr_{\partial\Omega} \Upsilon_j=\upsilon_j$, then
$U=\sum_j C_j\Upsilon_j$ is the solution to \eqref{div} with trace $u$
on $\partial\Omega$. Let $\gamma_k$ be a homotopy basis for $\Omega$,
$1\leq k\leq n$, and $\lambda_1,\cdots,\lambda_n$ be real constants.
The relations $\int_{\gamma_k}\partial_n U=\lambda_k$ are equivalent
to the linear system of equations:
\begin{equation}
\label{sysp}
\sum_{j=0}^n C_j\int_{\gamma_k}\partial \partial_n\Upsilon_j=\lambda_k,
\quad 1\leq k\leq n, \qquad\sum_{j=0}^n C_j=0.
\end{equation}
When all the $\lambda_k$  are zero, it follows from Lemma \ref{periodes}
and the first part of the proof that $C_j=0$ for all $j$ is the only solution.
Therefore, by elementary linear algebra, \eqref{sysp} has a unique solution
$C_0,C_1,\ldots, C_n$ for each $n$-tuple $\lambda_1,\cdots,\lambda_n$.
\hfill \boite
\begin{remarque}
\label{nteS}
If $u\in\mathcal{S}_\Omega$ and $U\in W^{2,r}(\Omega)$ is the solution to 
\eqref{div} such that $\tr_{\partial\Omega}U=u$, then $U$ certainly satisfies non-tangential estimates of the form \eqref{ntestimates}. Indeed, 
by the Sobolev embedding theorem and
Proposition \ref{propS}, we get that
\begin{equation}
\label{trmaj}
\|U\|_{L^\infty(\Omega)}\leq C\|U\|_{W^{1,r}(\Omega)}\leq C\|U\|_{W^{2,r}(\Omega)}\leq C'\|u\|_{L^p(\partial\Omega)}.\end{equation}
\end{remarque}

\subsection{Dirichlet problem for the conductivity equation} 
\label{secsoldiv}

The result below generalizes point $(i)$ of Lemma \ref{lemu} to multiply connected domains.
\begin{theorem} 
\label{Diru}
Let $\Omega$ be Dini-smooth and $\sigma\in W^{1,r}(\Omega)$ satisfy \eqref{ellipticsigma} with $p,r$ as in (\ref{hypothesesnup}). 
To each $u \in L^p_\RR(\partial \Omega)$, there is a unique 
$U \in{\cal U}^p(\Omega)$ whose non-tangential limit on 
$\partial \Omega$ is $u$. The function $U$ satisfies non tangential 
estimates of the form \eqref{ntestimates}. Moreover,
to every open set ${\cal O}$ with compact closure ${\overline{{\cal O}}}\subset\Omega$, 
there is a constant $c = c({\cal O},\Omega,\sigma,r,p)> 0$ such that 
\begin{equation}
\label{reginter}
\|U\|_{W^{2,r}({\cal O})} \leq c \, \|u\|_{L^{p}(\partial \Omega)} \, .
\end{equation}

\end{theorem}
{\it Proof.}
By Lemma \ref{conf} we may assume that 
$\Omega = \mathbb{D} \setminus \cup_{j=1}^n \overline{\mathbb{D}}_{a_j,r_j}$ 
is normalized circular, and we set by convention $a_0=0$, $r_0=1$.
Let $\nu$ be as in \eqref{dfnsigma}, so that $\nu$ meets 
$\eqref{kappa}$.
For $0\leq j\leq n$, we let $\nu_j \in W^{1,r}_\mathbb{R}(\overline{\mathbb{C}} \setminus \overline{\mathbb{D}}_{a_j,r_j})$ extend $\nu$ and satisfy
\eqref{kappa}, as in Theorem \ref{rmq:Hardy-ann2}. 
Subsequently, we let $\sigma_j \in W^{1,r}_\mathbb{R}(\overline{\mathbb{C}}
\setminus \overline{\mathbb{D}}_{a_j,r_j})$  be associated to $\nu_j$ 
through (\ref{dfnsigma}), so that it meets (\ref{ellipticsigma}).

Put $u = (u_0,u_1,\cdots,u_n) \in  \Pi_{j=0}^n \T_{a_j,r_j} = L^p_\RR(\partial \Omega)$. 
By Lemma \ref{lemu} point $(i)$, there exists a unique 
$U_0\in\mathcal{U}^p(\D)$
(with conductivity $\sigma_0$) such that $\tr_{\T} U_0 = u_0$.
Similarly,  for $1 \leq j \leq n$,  we let   $U_j\in \mathcal{U}^p(\CC \setminus \overline{\DD}_{a_j,r_j})$ (with conductivity $\sigma_j$)
be such that $\tr_{\T_{a_j,r_j}} U_j = u_j$.
Define operators 
$K_{i,j}: L^p_\R(\T_{a_i,r_i}) \to L^p_\R(\T_{a_j,r_j}), 0 \leq i\neq j \leq n$, by
\[
K_{i,j}(u_i) = U_{i_{|\mathbb{T}_{a_j,r_j}}}  \ , \quad i \neq j \ .
\]
Observe that $K_{i,j}$ is compact, for \eqref{usobloc} and the trace 
theorem imply that it maps continuously $L^p_\R(\T_{a_i,r_i})$ into
$W^{2-1/r,r}_\RR(\T_{a_j,r_j})\subset W^{1,r}_\RR(\T_{a_j,r_j})$ which is
compactly included in $L^l_\R(\T_{a_j,r_j})$ for all $l\in[1,\infty]$
by the Rellich-Kondratchov theorem.
 

Consider now the operator $\mathfrak{U}$ from $L^p(\partial\Omega)$ into itself
given by
\begin{equation}
\label{defmathfrakU}
\mathfrak{U}(u_0,\cdots,u_n) = \sum_{j=0}^n {U_j}_{|\partial\Omega}=
\Bigl(\sum_{j=0}^n {U_j}_{|\T}\,,\,\sum_{j=0}^n {U_j}_{|\T_{a_1,r_1}}\,,\cdots
\,,\,\sum_{j=0}^n {U_j}_{|\T_{a_n,r_n}}\Bigr).
\end{equation}
If we let $f_0 \in H^{p,0}_{\nu_0}(\mathbb{D})$ and
$f_j \in H^{p,0}_{\nu_j}(\overline{\mathbb{C}} \setminus \overline{\mathbb{D}_{a_j,r_j}}), 1 \leq j \leq n$ be the functions  granted by
Lemma \ref{lemu} point $(ii)$, with boundary condition 
$\mbox{Re }{f_j}_{|\T_{a_j,r_j}}=u_j$, we observe 
that $\mathfrak{U}(u) = \sum_{j=0}^n \mbox{Re } f_{j_{|\partial\Omega}} $ is 
the real part of $\tr_{\partial\Omega}\,f$,  where $f=\sum_{j=0}^n
{f_j}_{|\Omega}$ lies in $H_{\nu}^{p,0}(\Omega)$ 
since $\nu_{j_{|\Omega}} = \nu$ and we may use the characterization
by harmonic majorants in Theorem \ref{BS}. 
In particular 
\begin{equation}
\label{xx}
\mbox{Ran} \, \mathfrak{U} \subset \mbox{Re} \, \tr_{\partial \Omega} H_{\nu}^p(\Omega)
\subset \tr_{\partial \Omega}  \, {\cal U}^p(\Omega) \, .
\end{equation}
Moreover, we have the following matrix relation
\begin{equation}
 \label{truc1}
\mathfrak{U} = (I + K) 
\end{equation}
where $K:L^p(\partial\Omega)\to L^p(\partial\Omega)$ is given by
\begin{equation*}
\bigl(K(u_0,u_1,\cdots,u_n)\bigr)^T=
\begin{pmatrix}
 0   &   K_{1,0}   & K_{2,0} & \cdots & K_{n,0} \\
 K_{0,1} &    0    & K_{2,1} & \cdots & K_{n,1} \\
 K_{0,2} & K_{1,2} & 0 & \cdots & K_{n,2} \\
\vdots & & \ddots & \ddots & \vdots \\
K_{0,n} & \cdots & \cdots & K_{n-1,n} & 0
\end{pmatrix}
\begin{pmatrix}
u_0\\ u_1\\u_2\\\vdots\\u_n
\end{pmatrix}
\end{equation*}
Because the $K_{i,j}$ are compact, so is
$K$ from $L^p_\R(\partial \Omega)$ into itself.

Next, we prove that $\mbox{Ker} \, \mathfrak{U}=\mathcal{S}_{\Omega}$. 
Indeed let $u = (u_0,\cdots,u_n)\in L^p_\R(\partial \Omega)$ be such that
$\mathfrak{U}u= 0$. 
Using the notations introduced above, this means that
\begin{equation}
\label{42+}
\mbox{Re } \tr_{\partial \Omega} f=
\mbox{Re } \tr_{\partial \Omega}\sum_{j=0}^n f_{j_{|\Omega}} 
= 0 \, .
\end{equation}
Define further $\tilde{f} = f - i \int_{\mathbb{T}}\mbox{Im } \tr_{\mathbb{T}} f$ which lies in $H^p_\nu(\Omega)$ and set
\[
w = \frac{\tilde{f} - \nu \overline{\tilde{f}}}{\sqrt{1 - \nu^2}}
\in G^p_\alpha(\Omega),\qquad 
w= e^s F  \, ,
\] 
for the functions $w$, $s$, $F$ associated to 
$\tilde{f}$ through (\ref{ftow}) and (\ref{decompexpH}), 
with $\mbox{Im } s_{{|_\mathbb{T}}} = 0$ and, say,
$\mbox{Im } s_{{|_{\mathbb{T}_{a_j,r_j}}}} = \theta_j$ for  $1 \leq j \leq n$,
see Proposition \ref{expsf}. 
Clearly (\ref{defmathfrakU}) and (\ref{42+}) entail that  $\mbox{Re } \tr_{\mathbb{T}_{a_j,r_j}}w = 0$
for $0\leq j\leq n$. Thus, from the boundary conditions for  $s$,  
we see that $F$ has constant argument $-\theta_j+\pi/2$ modulo $\pi$ 
on each $\mathbb{T}_{a_j,r_j}$.
By Morera's theorem \cite[II, Ex. 12]{gar},
this allows one to reflect $F$  
across each $\mathbb{T}_{a_j,r_j}$ according 
to the rule
\begin{equation}
\label{reflecform}
\check{F}_j(z-a_j) =-e^{-2i\theta_j}\overline{F\left(\frac{r_j^2}
{\overline{z-a_j}}\right)},
\end{equation}
so that $F$ is in fact analytic on a neighborhood of $\overline{\Omega}$.
In particular, $F(\T_{a_j,r_j})$ is a segment $S_j$
on the line through the origin defined by 
$\{\arg z=-\theta_j+\pi/2\mbox{ mod }\pi\}$, and if $z_0\in\CC$ belongs to
none of the $S_j$ there is a single-valued branch of
$\log (F-z_0)$ on each $\T_{a_j,r_j}$, $0\leq j\leq n$. 
Thus, by the argument principle, $F$ cannot assume the value $z_0$ in $\Omega$,
but if $F$ is not constant $F(\Omega)$ is open, therefore it 
contains some $z_0\notin\cup_j S_j$. Hence $F$ is constant, all $\theta_j$ 
are equal to 0, and $F=ic$, $c\in\R$,
is a pure imaginary constant on $\T$.
Because $ \mbox{Im } \tr_{\mathbb{T}} \tilde{f}$ has vanishing mean 
by construction, and since 
\[
 \mbox{Im }  \tr_{\mathbb{T}} \tilde{f} =  \mbox{Im } \tr_{\mathbb{T}} \left(w \, \sqrt{\dfrac{1-\nu}{1+\nu}}\right) = c\, e^{\mbox{Re }s_{{|_\mathbb{T}}}} \, \sqrt{\dfrac{1-\nu}{1+\nu}}
\]
we must have $c = 0$ whence $F=w=\tilde{f} \equiv 0$. It follows that
\begin{equation}
\label{eq:iC}
\sum_{j=0}^n f_{j_{|\Omega}} = f = i C \ , \quad C \in \mathbb{R} \ .
\end{equation}
Now, since $f_j\in H_{\nu_j}^{p,0}(\overline{\C}\setminus\overline{\D}_{a_j,r_j})$ we find upon writing 
\[
 f = \underbrace{\left(f_{0_{|\Omega}} + \sum_{j=1}^n \int_{\T_{a_j,r_j}} 
\mbox{Re } f_{j_{|\Omega}}\right)}_{h_0}
+ \sum_{j=1}^n \underbrace{\left(f_{j_{|\Omega}} - \int_{\T_{a_j,r_j}} \mbox{Re } f_{j_{|\Omega}}\right)}_{h_j} \ ,
\]
that $f = h_0 + \sum_{j=1}^n h_j$ is the direct sum 
decomposition of $f$ furnished by Theorem \ref{rmq:Hardy-ann2}.
However, in view of  \eqref{eq:iC}, $f=iC+\sum_{j=1}^n 0$ is also such a 
decomposition, therefore by uniqueness $h_0 = iC$ and $h_j = 0$, 
$1\leq j\leq n$. Hence
$U_{j_{|\Omega}} = \mbox{Re } f_{j_{|\Omega}} = c_j \in \mathbb{R}$
for $j=0,1,\cdots,n$ and  $c_0 = -\sum_{j=1}^n c_j$ by (\ref{42+}), 
i.e. $u \in \mathcal{S}_{\Omega}$.

The reverse inclusion $ \mathcal{S}_{\Omega} \subset \mbox{Ker} \, 
\mathfrak{U}$ is immediate from \eqref{defmathfrakU} 
for $U_j$ is a constant
when $u_j$ is.

Finally, let us show that the Riesz number of the operator $\mathfrak{U}$ 
is equal to 1:
\begin{equation}
\label{eqRiesznb}
\mbox{Ker} \, \mathfrak{U} = \mbox{Ker} \, \mathfrak{U}^2 \, ,
\end{equation}
or equivalently $\mbox{Ker} \, \mathfrak{U}^2  \subset \mbox{Ker} \, 
\mathfrak{U}$.
Indeed, let $x \in \mbox{Ker} \, \mathfrak{U}^2$. Then 
\[
u = \mathfrak{U} \, x \in \mbox{Ker} \, \mathfrak{U} \cap \mbox{Ran} \, \mathfrak{U}=
\mathcal{S}_\Omega  \cap \mbox{Ran} \, \mathfrak{U} \subset \mathcal{S}_\Omega  \cap \mbox{Re} \, \tr_{\partial \Omega} H^{p}_{\nu}(\Omega) \, ,\]
in view of (\ref{xx}). By Proposition \ref{propS} it holds that
$\mathcal{S}_\Omega  \cap \mbox{Re} \, \tr_{\partial \Omega} H^{p}_{\nu}(\Omega) = \{0\}$,
therefore $u =0$ hence $x \in \mbox{Ker} \, \mathfrak{U}$, thereby establishing (\ref{eqRiesznb})

In view of what precedes, a theorem of F. Riesz\cite[Thm 1.16]{COLTKRESS}
implies the decomposition
\begin{equation}
\label{decompLp}
L^p_\RR(\partial \Omega)  = \mbox{Ker} \, \mathfrak{U} \oplus \mbox{Ran} \, 
\mathfrak{U}= \mathcal{S}_\Omega \oplus \mbox{Ran} \, \mathfrak{U}.
\end{equation} 
The existence of $U$ now follows from \eqref{decompLp}, \eqref{xx}, 
\eqref{ntestimates}, Theorem \ref{rmq:Hardy-ann2},  \eqref{usobloc},
Proposition \ref{propS} and Remark \ref{nteS}.

To prove uniqueness, suppose that $U\in\mathcal{U}^p(\Omega)$ satisfies
$\tr_{\partial\Omega}U=0$ and let $\gamma_j$ be a homotopy basis for $\Omega$,
$1\leq j\leq n$. By Proposition 
\ref{propS}, there is $\upsilon\in\mathcal{S}_\Omega$ and 
$\Upsilon\in W^{2,r}(\Omega)$  satisfying \eqref{div} such that
$\tr_{\partial\Omega}\Upsilon=\upsilon$ and 
$\int_{\gamma_j}\partial_n(U-\Upsilon)=0$ for all $j$. From
Lemma \ref{periodes}  it follows that 
$W=U-\Upsilon$ has a $\sigma$-harmonic conjugate $V$, {\it i.e.} 
$W+iV=f\in H^p_\nu(\Omega)$. Since $\mbox{Re}(\tr_{\partial\Omega}f)$ is equal
to some constant $C_j$ on $\TT_{a_j,r_j}$, the reflection formula 
\eqref{defcheck2} and the argument thereafter shows that 
$f\in W^{2,r}(\Omega)$. Then  $U\in W^{1,r}(\Omega)$ {\it a fortiori},
therefore $U\equiv0$ by uniqueness of $W^{1,r}(\Omega)$-solutions to the
Dirichlet problem \cite{aq,Fazio}.
\hfill \boite\\

Note that (\ref{decompLp}) and Proposition \ref{propS} immediately imply:
\begin{corollary}
\label{corS}
$\mbox{Ran} \, (\mathfrak{U}) = \mbox{Re} \, \tr_{\partial \Omega} H_{\nu}^p(\Omega)$.
\end{corollary}
\begin{remarque}
\label{majh}
From \eqref{decompLp}, Corollary \ref{corS}, Theorem \ref{BS} point $(ii)$, 
and Proposition \ref{propS}, we see that the condition 
``$U\in\mathcal{U}^p(\Omega)$'' may be replaced by `` $U$ solves \eqref{div} 
and $|U|^p$ has a harmonic majorant''.
\end{remarque}
It is standard in regularity theory that smoothness of the boundary
may be traded for smoothness of the coefficients. Here is an application of 
Theorem \ref{Diru} to the Dirichlet problem for equation \eqref{div} on 
non smooth domains. Given a weight $W\geq0$ on $\partial\Omega$,
we denote by $L^p(\partial\Omega,W)$ the weighted space of functions
$h$ for which $|h|^pW\in L^1(\partial\Omega)$.
\begin{corollary}
\label{Dirns}
Let $D\subset\overline{\C}$ be a finitely connected domain whose boundary is
a piecewise $C^{1,\lambda}$ polygon, $0<\lambda\leq1$, with $N$ vertices
$W_1,\cdots,W_N$. Let $\lambda_j\pi$ be the jump of the oriented tangent at
$W_j$, $-1\leq\lambda_j\leq1$, and assume that $\mu=\max\{\lambda_j\}<1$
({\it i.e.} there is no outward-pointing cusp). Define 
\[W(z)=\Pi_{j=1}^N|z-W_j|^{\lambda_j}.\]
If $\sigma\in W^{1,\infty}_\RR(D)$ meets \eqref{ellipticsigma} 
and $u\in L^p_\R(\partial D,W)$ for some
$p>2/(2-\max\{0,\mu\})$, there is a unique solution $U$ to \eqref{div} in $D$
such that $|U|^p$ has a harmonic majorant and $U$ has non tangential limit $u$ a.e. on $\partial D$.
\end{corollary}
{\sl Proof:} Let $\varphi$ map a circular domain $\Omega$ onto $D$. By Remark \ref{majh}, the statement is equivalent to the existence of a unique solution 
to \eqref{div} in $\mathcal{U}^p(\Omega)$ with $\sigma$ replaced by 
$\sigma\circ\varphi$ and
$u$ by $u\circ\varphi$. By \cite[Prop. 4.2]{BMSW}, $(W\circ\varphi)|\varphi'|$
extends continuously to $\overline{\Omega}$ and is never zero there\footnote{
This expresses that $W(z) |dz|$ is  comparable to  harmonic measure on
$\partial\Omega$.}.
Hence $\sigma\circ\varphi\in W^{1,r}(\Omega)$ for 
$r\in(2,2/\max\{0,\mu\})$, and the assumptions imply that
$p>r/(r-1)$ for $r$ close enough to $2/\max\{0,\mu\}$. Finally, the condition
on $p$ implies that $u\circ\varphi\in L^p(\partial\Omega)$, so that we can 
apply Theorem \ref{Diru}.
\hfill\boite\\

Although $\partial_n U$ needs not be a distribution on 
$\partial\Omega$ when $U\in\mathcal{U}^p(\Omega)$ (unless $p=\infty$), 
Theorem \ref{Diru} allows us to define 
$\sigma\partial_nU$.
Recall from (\ref{defDeltat}) the definition of
$\widetilde{\Delta}_\epsilon$.
\begin{corollary} 
\label{lemun}
Let $\Omega$ be a Dini-smooth domain and
$p$, $r$, $\sigma$  satisfy \eqref{dfnsigma} and \eqref{hypothesesnup}.
To each $U\in\mathcal{U}^p(\Omega)$, there is a unique distribution 
$\sigma\partial_n U  \in W^{-1,p}_\R(\partial \Omega)$
such that
\begin{equation}
\label{titi1}
\int_{\partial \Omega} \sigma \partial_n U  \, \varphi := 
\lim_{\epsilon \to 0} \iint_{\widetilde{\Delta}_\epsilon} \sigma \nabla U \cdot \nabla \varphi,\qquad \varphi \in \mathcal{D}(\R^2).
\end{equation} 
Further, there is a constant $C=C(\Omega, \sigma, r,p)$ such that
\begin{equation}
\label{titi2}
\|\sigma \partial_n U \|_{W^{-1,p}(\partial \Omega)} \leq C \|\mbox{tr}_{\partial \Omega} U\|_{L^p(\partial \Omega)} \, .
\end{equation}
\end{corollary}

{\sl Proof:} we may assume that $\Omega$ is bounded.
If $\mbox{tr}_{\partial \Omega} \,U\in\mathcal{S}_\Omega$, then
$U\in W^{2,r}(\Omega)$ by Proposition \ref{propS} and $\partial_n U$ 
may be defined a.e. on $\partial\Omega$ as the scalar product  of
$\tr_{\partial\Omega}\nabla U$ with the unit normal to $\partial\Omega$.
Thus, \eqref{titi1} follows 
from the divergence formula (the limit on the right is
equal to $\iint_\Omega\sigma\nabla U.\nabla\varphi$ by dominated convergence)
while \eqref{titi2} drops out from \eqref{trmaj} and the Sobolev
embedding theorem (because $\|\nabla U\|_{L^\infty(\Omega)}\leq C\|U\|_{W^{2,r}(\Omega)}$).
By \eqref{decompLp} and Corollary \ref{corS}, it remains to
handle the case where $U=\mbox{Re}f$ with $f=U+iV\in H^{p,0}_\nu(\Omega)$. 

In this case we know from \eqref{system} 
that $\sigma \partial_n U = \partial_t V$ on 
$\partial\widetilde{\Delta}_\epsilon$. Hence, using 
the Green formula and integrating by parts, we obtain for $\varphi\in\mathcal{D}(\RR^2)$:
\begin{equation}
\label{Greeninter}
  \iint_{\widetilde{\Delta}_\epsilon} \sigma \nabla U \cdot \nabla \varphi \,=
\int_{\partial\widetilde{\Delta}_\epsilon} \sigma \partial_n U  \, \varphi = - \int_{\partial\widetilde{\Delta}_\epsilon} \tr_{\partial\widetilde{\Delta}_\epsilon} V \partial_{t} \varphi.
\end{equation}
Consequently, by Property \ref{1} in section \ref{sechpnud} and Lemma \ref{conf}, we obtain
\begin{equation}
\label{defnormint} 
\lim_{\epsilon \to 0} \iint_{\widetilde{\Delta}_\epsilon} \sigma \nabla U \cdot \nabla \varphi = - \int_{\partial \Omega} \tr_{\partial \Omega}V  \, \partial_t \varphi    \, .
\end{equation}
As $\|\tr_{\partial\Omega}V \|_{L^p(\partial \Omega)}\leq c\|\tr_{\partial \Omega} U\|_{L^p(\partial \Omega)}$ by 
Theorem \ref{suite-thm1-lem1iiA}, the right hand side of \eqref{defnormint}
indeed defines a distribution  $\sigma\partial_nU\in W^{-1,p}_\RR(\partial \Omega)$ satisfying
our requirements.\hfill\boite

\subsection{The Dirichlet problem in $H_{\nu}^p(\Omega)$} 
\label{secHp}
%


We are now in position to solve the ``Dirichlet problem'' for equation 
(\ref{dbar}). Given $\Omega$ a Dini-smooth domain,
we put $\partial\Omega=\cup_{j=0}^n\Gamma_j$ where the $\Gamma_j$ are 
disjoint Dini-smooth Jordan curves.
Introduce for 
$U \in \mathcal{U}^p(\Omega)$ the compatibility condition 
\begin{center}
$(H_{\Omega,\sigma})$:~~~$\int_{\Gamma_j} \sigma \partial_n U = 0 \, , \ 0 \leq j \leq n$ 
\end{center}
where the normal derivative is understood in the sense of
Corollary \ref{lemun}.
Note that the relation $\int_{\partial \Omega} \sigma \partial_n U = 0$,
which follows from \eqref{titi1} when $\varphi=1$ on $\Omega$,
is to the effect that $(H_{\Omega,\sigma})$ holds as soon as 
$\int_{\Gamma_j} \sigma \partial_n U = 0$ for every $j$ {\sl but one}.

\begin{theorem}
\label{suite-thm1-lem1iiA}
Let $\Omega$ be a Dini-smooth domain, $\partial\Omega=\cup_{j=0}^n\Gamma_j$,
 and
$\sigma,p,r,\nu$ meet (\ref{dfnsigma}) and (\ref{hypothesesnup}). 
If $u \in L_\R^p(\partial \Omega)$ satisfies $(H_{\Omega,\sigma})$, and only in this case, there exists a unique $f \in H_{\nu}^{p,0}(\Omega)$ 
such that
$\mbox{Re } \tr_{\partial \Omega} f = u$ {a.e. on } $\partial \Omega$.
Moreover, there is a constant  $c=c_{p,\alpha,\nu} > 0$ such that 
$ \|\tr_{\partial \Omega} f\|_{L^{p}(\partial \Omega)}\leq c\|u\|_{L^{p}(\partial \Omega)} \ .$
\end{theorem}
{\it Proof.}  
Let $u\in L^p(\partial\Omega)$ and
$U \in \mathcal{U}^p(\Omega)$ be such that $u = \tr_{\partial \Omega} U$,
see Theorem \ref{Diru}.
From \eqref{defDeltat} and the 
definition of
$P_{\partial\Omega,\epsilon}$ (after equation \eqref{defntmax}), we see
that the connected components of $\partial\widetilde{\Delta}_\varepsilon$
form a system of smooth Jordan curves $\gamma_{j,\varepsilon}$, 
$0\leq j\leq n$, with $\gamma_{j,\varepsilon}$ homotopic to 
$\gamma_{j,\varepsilon'}$ for any two $\varepsilon,\varepsilon'$ small enough.
Hence $\int_{\gamma_{j,\varepsilon}}\partial_nU$ is independent of
$\varepsilon$ by Lemma \ref{periodes}, and letting $\varphi$ 
in equation \eqref{titi1} be 1 on a 
neighborhood of $\Gamma_j$ and $0$ on a neighborhood of $\Gamma_k$
for  $k\neq j$ , we deduce from 
Corollary \ref{lemun} (see \eqref{Greeninter}) 
that $\int_{\gamma_{j,\varepsilon}}\sigma\partial_nU=\int_{\Gamma_j}
\sigma\partial_nU$.
Since  the $\gamma_{j,\varepsilon}$ are a homotopy basis of
$\Omega$ for $1\leq j\leq n$, we conclude
from Lemma \ref{periodes} that
$U$ has a $\sigma$-harmonic conjugate $V$ if and only if 
$(H_{\Omega,\sigma})$ holds. Adding a constant to $V$ if necessary,
we can ensure that $f=U+iV\in H^{p,0}_\nu(\Omega)$. Uniqueness of $f$
follows from uniqueness of $\mathcal{U}$ and the fact that any two 
$\sigma$-harmonic conjugates differ by a constant. 
\hfill \boite

\begin{remarque}
Let $E(p,\sigma)\subset L^p_\RR(\partial\Omega)$ 
denote the closed subspace of functions
with zero mean meeting $(H_{\Omega,\sigma})$. 
Taking into account that $f\in H^p_\nu(\Omega)$ if, and only if 
$(if) \in H^p_{-\nu}(\Omega)$,
it follows from Theorem \ref{suite-thm1-lem1iiA} that the $\sigma$-conjugating
map $\mathcal{H}(u)=\mbox{\rm Im } tr_{\partial \Omega} f$ is an 
isomorphism
from $E(p,\sigma)$ onto $E(p,1/\sigma)$
satisfying $\mathcal{H}^2=-Id$.
\end{remarque}
\subsection{Neumann problem for the conductivity equation} 
\label{secsoldn}
Theorem \ref{suite-thm1-lem1iiA} allows us to solve a weighted 
Neumann problem for \eqref{div}, where data consist of the normal 
derivative of $u$ multiplied by the conductivity on $\partial\Omega$:
\begin{theorem} 
\label{Neumann}
Let $\Omega$ be a Dini-smooth domain and
$p$, $r$, $\sigma\in W^{1,r}(\Omega)$  
satisfy \eqref{ellipticsigma} and \eqref{hypothesesnup}.
To each $\phi\in W^{-1,p}_\R(\partial \Omega)$ such that
$\int_{\partial\Omega} \phi=0$, there is
$U\in\mathcal{U}^p(\Omega)$, unique up to an additive constant,
such that  $\sigma\partial_n U  =\phi$.
\end{theorem}

{\sl Proof:} by Proposition \ref{propS} and the proof of Theorem 
\ref{suite-thm1-lem1iiA}, there is $\Upsilon\in \mathcal{U}^p(\Omega)$
with $\tr_{\partial\Omega} U\in\mathcal{S}_\Omega$ such that
$\psi= \phi-\sigma\partial_n\Upsilon$  satisfies
$\int_{\Gamma_j}\psi=0$ for $0\leq j\leq n$.
Because $\psi_{|\Gamma_j}$ lies in $W^{-1,p}(\Gamma_j)$ 
and annihilates the constants, $\psi$  is of the form $\partial_t v$ for 
some $v\in L^p(\Omega)$. By Proposition \ref{propS} and Theorems
\ref{Diru}, \ref{suite-thm1-lem1iiA} applied with $1/\sigma$ and $-\nu$
rather than $\sigma$ and $\nu$, we can add to $v$ an element of 
$\mathcal{S}_\Omega$ 
(this does not change $\partial_t v$) so that 
$v=\mbox{Im}\,\tr_{\partial\Omega} f$ for some $f=W+iV\in H^{p,0}_\nu$.
From \eqref{Greeninter}, \eqref{titi1}, and Property \ref{1},
it follows that 
\[\int_{\partial \Omega} \sigma \partial_n W  \, \varphi
=\int_{\partial \Omega} \partial_t v  \, \varphi
=\int_{\partial \Omega} \psi  \, \varphi\,,\qquad \varphi\in\mathcal{D}(\R^2),
\]
hence $\sigma\partial_n W=\psi$. Then $U=W+\Upsilon$ satisfies our 
requirements. In another connection, if
$U\in\mathcal{U}^p(\Omega)$ is such that $\sigma\partial_n U=0$ then 
$U=\mbox{Re} f$ with $f=U+iV\in H^{p,0}_\nu(\Omega)$ by Theorem 
\ref{suite-thm1-lem1iiA}. Thus, we deduce from \eqref{defnormint} and
Property \ref{1} that $\int_{\partial\Omega} \tr_{\partial\Omega}
V\partial_t \varphi=0$ for $\varphi\in\mathcal{D}(\RR^2)$. 
Since in addition $\tr_{\partial\Omega}V$ has zero mean by contruction, it 
is the zero function so that $V\equiv0$ by the uniqueness part of
Theorem \ref{Diru} (applied with $1/\sigma$ rather than $\sigma$).
Finally, $U$ must be constant by \eqref{system}.
\hfill\boite

\section{Density of traces, approximation issues}
\label{sec:density}
On a rectifiable curve, we let $|E|$ indicate arclength of 
a measurable set $E$. 

We consider the following 

{\bf Conjecture}
\emph{Let $\Omega$ be a Dini-smooth domain and
$p$, $r$, $\nu$  satisfy  \eqref{hypothesesnup}.
If $E\subset\partial\Omega$ satisfies $|E|<|\partial\Omega|$, 
then $\bigl(\tr_{\partial\Omega}H_\nu^p(\Omega)\bigr)_{|E}$ is dense in 
$L^p(E)$.}\\

When  $\nu=0$\footnote{We deal then with holomorphic Hardy spaces
in which case we may take $r=\infty$ and $p\in(1,\infty)$.} and $E$ is closed, 
this is an easy consequence of Runge's theorem. Still when $\nu=0$, 
but this time
$E$ is arbitrary, it was proven to hold in \cite{bl} when $\Omega=\D$, 
hence it is true for all simply connected $\Omega$ by conformal invariance
of Hardy spaces and Lemma \ref{conf}. These results are
of key importance to the approach 
of bounded extremal problems developed in \cite{bl, cp,cps,flps,jmlp}.

When $\nu \in W^{1,\infty}(\Omega)$ and $\Omega$ is simply connected, 
the answer is again positive  as established on the disk in
\cite{BLRR}. Below, we prove that if $\Omega$ is multiply connected but 
$E$ is contained in a single connected component of $\partial\Omega$, then the
statement is indeed correct:

\begin{theorem} 
\label{densityA}
Let $\Omega$ be a multiply connected
Dini-smooth domain, $\partial\Omega=\cup_{j=0}^n\Gamma_j$, $n>0$,
and $\nu,p,r,\nu$ meet  (\ref{hypothesesnup}). 
Let $E \subset \Gamma_{j_0}$ for some $j_0\in\{0,\cdots,n\}$.
If $|E|<|\Gamma_{j_0}|$, then restrictions to $E$
of traces of $H_\nu^p(\Omega)$-functions are dense in $L^{p}(E)$. 
\end{theorem}

To prove Theorem \ref{densityA}, we establish in Appendix \ref{app6}
the following result which may be of independent interest.
Note that the arguments of proof are in fact  to the effect that the result holds if there is at least one connected component of the boundary that does not intersect $E$. 
\begin{proposition}
\label{casinvHp}
Let $\Omega=\D$   and assume that 
$\nu$, $p$, $r$ satisfy \eqref{hypothesesnup}. 
Let moreover $\alpha\in L^r(\D)$ and
$\varrho \in (0,1)$. 

Then $\Bigl(H_{\nu}^p(\DD)\Bigr)_{|\D_\varrho}$ (resp. 
$\Bigl(G_{\alpha}^p(\DD)\Bigr)_{|\D_\varrho}$)
is dense in  $H_{\nu_{|\D_\varrho}}^p(\DD_\varrho)$ (resp. 
$G_{\alpha_{|\D_\varrho}}^p(\DD_\varrho)$).
\end{proposition} 
{\it Proof of Theorem \ref{densityA}.} 
Assume without loss of generality that $\Omega=\D\setminus\cup_{j=1}^n 
\overline{\D}_{\xi_j,r_j}$ and that $j_0=1$. Extending if necessary
$\alpha$  by zero to each $\D_{\xi_j,r_j}$ with $2\leq j  \leq n$, we 
are back to the case $n=1$ and then to
$\Omega=\A_\varrho$.
Clearly, it is enough to consider $E=\T_\varrho$.

Let $\nu_e\in W^{1,r}(\overline{\C}\setminus\overline{\D}_\varrho)$ extend 
$\nu$ and satisfy \eqref{kappa} (see proof of Theorem \ref{rmq:Hardy-ann2}),
and set $\nu_i\in W^{1,r}(\overline{\D})$ to
be $\nu\vee\check{\nu}_e$ where 
$\check{\nu}_e(z)=\nu_e(\varrho^2/\overline{z})$ for $z\in\D_\varrho$.
By Lemma \ref{cor-oplus}, any $\psi\in L^p(\T_\varrho)$ can be written as 
$\tr_{\T_\varrho}\psi_1+\tr_{\T_\varrho}\psi_2$ where 
$\psi_1\in H^p_{\check{\nu}_e}(\D_\varrho)$ and
$\psi_2\in  H^{p,00}_{\nu_e}(\overline{\C}\setminus\overline{\D}_\varrho)$.
By Proposition \ref{casinvHp}, to each $\varepsilon>0$ there is $\psi_3\in
H^p_{\nu_i}(\D)$ such that $\|\tr_{\T_\varrho}\psi_1-\tr_{\T_\varrho}\psi_3\|_{L^p(\T_\varrho)}<\varepsilon$.
Since $(\psi_2+\psi_3)_{|\A_\varrho}\in H^p_\nu(\A_\varrho)$, the result follows.
\hfill\boite\\




We shall illustrate the use of Theorem \ref{densityA} 
in  certain bounded extremal problems (BEP) in $H_{\nu}^p({\A_\varrho})$
which play an important role in the works \cite{TheseYannick,flps,fl,FMP}
where inverse boundary problems for equation \eqref{div} are considered.
 
We let now $\Omega=\A_\varrho$ for some $\varrho\in(0,1)$ and
we assume that $p,r,\nu$ meet (\ref{hypothesesnup}). 
Fix $I\subset \T_\varrho$ with $|I|>0$, and define 
$J=\partial\A_\varrho\setminus I$. 
To $M>0$ and $\phi \in L^p_\R(J)$, we associate 
the following subset of 
$L^p(I)$.

\[ 
{\mathcal B}_p^{\AA_\varrho}=\left\{g_{|I}:\ g\in \tr_{\partial \A}  H^{p}_{\nu}(\AA_\varrho) ;\ \left\Vert \mbox{Re} \,  g - \phi\right\Vert_{L^{p}(J)}\leq M\right\}\subset L^{p}(I) \, . 
\] 
Note that a function $g\in H^p_\nu(\A_\varrho)$ is completely determined by 
$g_{|I}$ in view of Property \ref{3}, Section \ref{sechpnud}.

The theorem below extends to annular geometry with weaker smoothness 
assumptions a result obtained on the disk with Lipschitz-continuous
$\nu$ in \cite[Thm 3]{flps} (see also \cite{bl, cp, cps,jmlp} for the case 
$\nu = 0$).

\begin{theorem} \label{BEP2} 
Let notations and assumptions be as above.
Then, to every
$F_d\in L^{p}(I)$, there exists a unique function $g_*\in {\mathcal B}_p^{\AA_\varrho}$ such that 
\begin{equation} \label{pb:bep} \tag{BEP} 
\left\Vert F_d-g_*\right\Vert_{L^{p}(I)}=\min_{g\in {\mathcal B}_p^{\AA_\varrho}} \left\Vert F_d-g\right\Vert_{L^{p}(I)}. 
\end{equation}
Moreover, if $F_d\notin \bigl({\mathcal B}_{p}^{\AA_\varrho}\bigr)_{|I}$, then $\left\Vert \mbox{Re} \, g_*- \phi\right\Vert_{L^{p}(J)}=M$. 
\end{theorem} 
{\it Proof.}
Since ${\mathcal B}^{\AA_\varrho}_p$ is a convex subset of the uniformly convex Banach space 
$L^{p}(I)$, 
it is enough, in order to prove existence and uniqueness of 
$g_*$, to check that ${\mathcal B}^{\AA_\varrho}_p$ is closed in
$L^{p}(I)$  \cite[Prop. 5]{beauzany}. Let $\varphi_k\in H^p_\nu(\A_\varrho)$
be such that $\varphi_{k|_I} \in {\mathcal B}^{\AA_\varrho}_p$ converges 
to some function $\varphi_I$ in $L^p(I)$ as $k \to \infty$. 
Then $\tr_{\partial\A_\varrho}\varphi_k$ is bounded in 
$L^p(\partial \A_\varrho)$
by definition of ${\mathcal B}^{\AA_\varrho}_p$.  
Hence, extracting a subsequence if necessary, we may assume that
$(\tr_{\partial\A_\varrho}\varphi_k)$ converges weakly in 
$L^p(\partial \A_\varrho)$ to 
$(\tr_{\partial\A_\varrho}\psi)$ for some
$\psi \in H^p_\nu(\partial \AA_\varrho)$
by Properties \ref{2}-\ref{3} in Section \ref{sechpnud} and
Mazur's theorem, \cite{brezis}. Because $(\tr_{J}\varphi_k)$ 
{\it a fortiori} converges weakly to $\tr_{J}\psi)$
in  $L^p(J)$, we deduce from the weak-* compactness of balls that
$\psi_{I}\in {\mathcal B}^{\AA_\varrho}_p$.
Moreover, we must have $\psi_{|_I} = \varphi_I$ by the strong convergence of 
$(\varphi_k)_{|I}$ in $L^p(I)$, hence ${\mathcal B}^{\AA_\varrho}_p$ is 
indeed closed.


\medskip
To prove that $\left\Vert \mbox{Re } g_* - \phi\right\Vert_{L^{p}(J)}=M$
when  $F_d\notin {\mathcal B}^{\AA_{\varrho}}_p$
assume for a 
contradiction that $\left\Vert \mbox{Re } g_* - \phi\right\Vert_{L^{p}(J)}<M$. 
By Theorem \ref{densityA}, there is a function $h\in \tr_{\partial \A_\varrho} H^p_\nu(\A)$ such that 
\[
\|F_d-g_*- h\|_{L^p(I)}<\|F_d-g_*\|_{L^p(I)}
 \, ,
\]
and by the triangle inequality we have
\[
\|F_d-g_*-\lambda h\|_{L^p(I)}<\|F_d-g_*\|_{L^p(I)}
\]
for all $0 < \lambda < 1$.
Now for $\lambda>0$ sufficiently small we have
$\|\mbox{Re } (g_* + \lambda  h )- \phi \|_{L^p(J)} \le M$, contradicting the optimality of 
$g_*$. \hfill\boite

\section{Conclusion}
\label{sec:conclu}

We developped in this paper a theory of Hardy spaces and conjugate functions on Dini smooth domains for the conjugate Beltrami equation that runs parallel to the holomorphic case. We conjecture the assumptions 
$\nu \in W^{1,r}$, $r>2$ and  $p>r/(r-1)$ are best possible for the above mentioned results to hold. We applied our results to Dirichlet and Neumann problems for the conductivity equation with $L^p$ and $W^{-1,p}$ data. Whether those continue to hold in higher dimension \cite{BEP3d} and for matrix-valued conductivity coefficients is an interesting open question.\\


\noindent{\small{{\bf Acknowledgements.} }} 
The research of the authors was partially supported by grant
AHPI (ANR-07-BLAN-0247) and the ``Region PACA''.



\newpage
{\Large{\textbf{Appendix}}}

\bigskip
\appendix
\numberwithin{equation}{section}
\section{Conformal maps of Dini-smooth annular domains}
\label{app1}
As is well-known \cite[Thm 3.5]{Pommerenke},
a conformal map between Dini-smooth simply connected domains extends
to a homeomorphism of their closures, and the derivative extends
continuously to the closure of the initial domain in such a way that it is 
never zero. This remains true in the multiply connected case,
but the authors could not 
locate the result in the literature which is why we provide a proof.
\begin{lemma}
\label{conf}
Let $\varphi$ conformally map a Dini-smooth domain $\Omega$ onto a Dini-smooth
domain $\Omega'$. Then $\varphi$ extends to a homeomorphism
from $\overline{\Omega}$ onto $\overline{\Omega'}$ whose derivative also 
extends continuously to $\overline{\Omega}$ and is never zero there.
\end{lemma}
{\it Proof.} That $\varphi$ extends to a homeomorphism
from $\overline{\Omega}$ onto $\overline{\Omega'}$ can be proved as in the 
simply connected case (see {\it e.g.}  \cite[Thm 14.18]{Rudin}), 
granted that each boundary point of $\Omega$ is 
accessible, by Dini-smoothness of $\partial\Omega$,
and that every bounded analytic function on $\Omega$ 
has nontangential limits
at almost every boundary point \cite[Thms 10.3, 10.12]{duren}.
We  are thus left to show that $\phi'$ extends in a continuous nonvanishing 
manner to $\overline{\Omega}$.

For this, observe that it is enough to consider the doubly connected 
case. For if $J$ is one of the Jordan curves composing $\partial\Omega$
and $J'$ is another Dini-smooth  Jordan curve contained in 
$\Omega$, disjoint from $J$, such that
the annular region $A(J,J')$ between $J$ and $J'$ lies entirely in $\Omega$,
then $\varphi$ conformally maps $A(J,J')$ onto some annular region in 
$\Omega'$ whose boundary  consists of two Dini-smooth Jordan curves,
one of which is a connected component of $\partial\Omega'$ (by what precedes).
If $\varphi'$ continuously extends to $J$ in a nonvanishing manner,
we will be done since $J$ was an arbitrary 
connected component of $\partial\Omega$.

Now, let $\Omega$  be doubly connected and lie between two Dini-continuous 
Jordan curves
$\Gamma_1$, $\Gamma_2$, the latter being interior to the former. 
Let $\psi_1$ map the interior $\Omega_1$ of $\Gamma_1$ onto the unit disk 
$\D$. 
Because $\Gamma_1$ is Dini-smooth, $\psi_1$  extends to a homeomorphism 
from $\overline{\Omega_1}$ onto 
$\overline{\D}$ and the derivative $\psi_1'$ extends continuously to
 $\overline{\Omega_1}$ and is never zero there. Clearly $\psi_1(\Gamma_2)$ is
a Dini-smooth Jordan curve. Let  $\Omega_2$ indicate the interior of 
$\Gamma_2$  and $\psi_2$ conformally map 
$\overline{\C}\setminus \psi_1(\overline{\Omega_2})$ onto
$\overline{\C}\setminus \overline{\D}$. Then $\psi_3:=\psi_2\circ\psi_1$ 
maps $\Omega$ 
onto an annular region $\Omega_3$ bounded by analytic Jordan curves
(namely a cicle and an analytic image of a circle), $\psi_3$ 
extends to a homeomorphism of the closures, and by the chain rule 
$\psi_3'$ extends 
continuously to $\overline{\Omega}$ where it is never zero.
Let $\varrho\in(0,1)$ be such that $\psi_4$ conformally maps 
$\Omega_3$ onto $\A_\varrho$.
Then $\psi_4$ extends to a homeomorphisn 
from $\overline{\Omega_3}$  onto $\overline{\A_\varrho}$,
and since $\partial\Omega_3$ consists of analytic curves 
it follows from the reflexion principle that 
$\psi_4$ extends analytically and locally injectively to a neighborhood of 
$\overline{\Omega_3}$. Altogether, we constructed a conformal map
from $\Omega$ onto $\A_\varrho$, namely  $\psi_4\circ\psi_3$,
that extends continuously from
$\overline{\Omega}$ onto $\overline{\A_\varrho}$, and whose derivative extends 
continuously to $\overline{\Omega}$ where it is never zero.
Because self-conformal maps of  $\A_\varrho$ must be M\"obius transforms, 
similar properties hold for any conformal map from $\Omega$ onto $\A_\varrho$.
The same is true of $\Omega'$ which is conformally equivalent to $\Omega$ 
and therefore to the same $\A_\varrho$. Factoring $\varphi$ into a 
conformal map from $\Omega$ onto  $\A_\varrho$ followed by a conformal map
from $\A_\varrho$ onto $\Omega'$ ({\it e.g.} $(\psi_4\circ\psi_3)\circ
(\psi_4\circ\psi_3)^{-1}\varphi$), we get the desired result.
\hfill \boite

\section{Proof of Proposition \ref{expsf}}
\label{app2}
We may assume that $\Omega$ is bounded.
Set by convention $\overline{w(\xi)}/w(\xi)=0$ if $w(\xi)=0$ and define
\[\lambda(z)=
\frac{1}{2i\pi}\iint_\Omega \frac{\overline{w(\xi)}}{w(\xi)}\frac{\alpha(\xi)}
{\xi-z}\,d\xi\wedge
d\overline{\xi},\qquad z\in\Omega.
\]
As in \eqref{defA}, we find that $\lambda\in W^{1,r}(\Omega)$
with $\overline{\partial}\lambda=\alpha\overline{w}/w$ and 
$\|\partial \lambda\|_{L^r(\Omega)}\leq C_1
\|\alpha\overline{w}/w\|_{L^r(\Omega)}\leq C_1\|\alpha\|_{L^r(\Omega)}$
for some constant $C_1=C_1(r)$, see \cite[Ch. 1, (1.7)-(1.9)]{im2}.
Thus, if we set $s(z)=\lambda(z)-\iint_\Omega \lambda\,dm/m(\Omega)$,
we obtain by Poincar\'e's inequality that 
\begin{equation}
\label{majua}
\|s\|_{W^{1,r}(\Omega)}\leq C_2 (\|\partial s\|_{L^r(\Omega)}+
\|\overline{\partial} s\|_{L^r(\Omega)})=C_2
(\|\partial \lambda\|_{L^r(\Omega)}+
\|\overline{\partial} \lambda\|_{L^r(\Omega)})\leq 2C_2C_1
\|\alpha\|_{L^r(\Omega)}
\end{equation}
for some constant $C_2=C_2(r,\Omega)$.
Now \eqref{regs} follows from the Sobolev embedding theorem, since $r>2$.

Next, we show that $F=e^{-s}w\in L^p_{loc}(\Omega)$ is in fact holomorphic.
By Weyl's lemma, it is enough to check 
that $\bar\d F=0$ as a distribution. 
Let $\psi\in {\cal D}(\Omega)$ and $\psi_n$ a sequence in 
${\cal D}(\R^2)_{|_\Omega}$ converging to $s$ in $W^{1,r}(\Omega)$. As $r>2$,
$\psi_n$ converges uniformly to $s$ on $\Omega$ by Sobolev's 
embedding theorem,
hence by dominated convergence and since 
$\alpha\overline{w}\in L^1_{loc}(\Omega)$ 
$$
\langle \bar\d F,\psi\rangle=
-\langle e^{-s}w,\bar\d\psi\rangle=
-\lim_n\langle w, e^{-\psi_n} \bar\d\psi\rangle=
-\lim_n\langle w,\bar\d(e^{-\psi_n} \psi)+\psi e^{-\psi_n}\bar\d\psi_n\rangle$$
$$
=\lim_n\langle \alpha\overline{w}, e^{-\psi_n}\psi\rangle
-\lim_n\langle w,\psi e^{-\psi_n}\bar\d\psi_n\rangle
=\langle e^{-s}(\alpha\overline{w}-w\bar\d s),\psi\rangle=0
$$
since $w\bar\d s=\alpha \overline w$, where we used in the fourth 
equality that $e^{-\psi_n}\psi\in{\cal D}(\Omega)$.
This proves \eqref{decompexpH}. Because 
$s\in W^{1,r}(\Omega)$ is bounded we have that $e^s\in W^{1,r}(\Omega)$, and 
as  $F$ is locally smooth we get that $w\in W_{loc}^{1,r}(\Omega)$, 
as announced.

Clearly $F$ satisfies \eqref{systema} if and only if $w$ does by \eqref{regs},
{\it i.e.} $F\in H^p(\Omega)$ if and only if $w\in G_\alpha^p(\Omega)$.

As for the normalization, let 
$u\in W^{1,r}_\RR(\Omega)$ be harmonic in $\Omega$ with 
$u_{|\partial\Omega}={\rm Im}\,s_{|\partial\Omega}\in W^{1-1/r,r}_\RR
(\partial\Omega)$. Such a function uniquely exists with 
$\|u\|_{W^{1,r}(\Omega)}\leq
C_3\|{\rm Im }\,s\|_{W^{1-1/r,r}(\partial\Omega)}$, where $C_3=C_3(r,\Omega)$
\cite{campanato}. Thus, by \eqref{majua} and continuity of the trace, it 
holds that 
\[
\|u\|_{W^{1,r}(\Omega)}\leq C_4\|\alpha\|_{L^r(\Omega)},\ \ 
\mbox{\rm with}\ C_4=C_4(r,\Omega).
\]
Set $a_j=\int_{\Gamma_j} \partial_n u$, $0\leq j\leq n$, where 
$\partial_n u\in W^{-1/r,r}_\RR(\partial\Omega)$.
Note that $\sum_ja_j=0$ by \eqref{Greenn} (applied with $\sigma=g\equiv1$).
We can find a function $\omega$, harmonic on $\Omega$ and $C^1$-smooth on  
$\overline{\Omega}$, which is constant on each $\Gamma_j$  and such that
$\int_{\Gamma_j} \partial_n \omega=a_j$,
see \cite[Sec. 6.5.1]{ahlfors00} and Lemma \ref{conf}.
By construction 
\[\|\omega\|_{W^{1,r}(\Omega)}\leq C_5\|\partial_n u\|_{W^{-1/r,r}(\partial\Omega)}\leq C_5\|u\|_{W^{1,r}(\Omega)}\leq C_6
\|\alpha\|_{L^r(\Omega)}
\ \ \mbox{\rm with}\ C_6=C_6(r,\Omega).\]
The harmonic function $v=u-\omega$ lies in $W^{1,r}(\Omega)$ and its conjugate
differential $d^*v=-\partial v/\partial y dx+\partial v/\partial x dx$ 
is exact\footnote{Indeed, since $d^* v=\partial_n v|dz|$ along any curve,
its integral over a cycle $\gamma\subset\Omega$ is 
zero by Green's formula \eqref{Greenn} (applied with $\sigma=g\equiv1$ on
the domain bounded by $\gamma$ and all the $\Gamma_j$ located inside $\gamma$)
because $\int_{\Gamma_j}\partial_n v=0$ for all $j$ by construction,
see \cite[Sec. 4.6]{ahlfors00}.}. Thus, there is a harmonic conjugate 
$\tilde{v}$ in $\Omega$, unique up to an additive constant, such that
$G=v+i\tilde{v}$ is holomorphic; if we normalize $\tilde{v}$ so that 
$\int_\Omega \tilde{v}dm=0$, it is immediate from the Cauchy-Riemann 
equation and Poincar\'e's inequality that 
$\|\tilde{v}\|_{W^{1,r}(\Omega)}\leq C_7 \|v\|_{W^{1,r}(\Omega)}$ with
$C_7=C_7(r,\Omega)$. Altogether,
$\|G\|_{W^{1,r}(\Omega)}\leq C_8\|\alpha\|_{L^r(\Omega)}$, and  since $r>2$
we see that $G$ is bounded by the Sobolev embedding theorem.
Finally, setting $\tilde{s}=s-iG$ and $\tilde{F}=e^{iG}F$, 
we find that $w=e^{\tilde{s}}\tilde{F}$ is a factorization of the 
form \eqref{decompexpH} in which ${\rm Im}\,\tilde{s}$ is constant on each 
$\Gamma_j$. Clearly, we may impose the value of this constant on any given 
$\Gamma_j$ upon renormalizing $\tilde{F}$. \hfill\boite.

\section{Proof of Lemma \ref{holoHp}}
\label{app3}
We must show that if $f$ is holomorphic in $\Omega$ and
\eqref{systema} holds for some sequence of admissible
compact sets $\Delta_n$, then it holds for $\widetilde{\Delta}_n$ 
defined in \eqref{defDeltat} as well.
When $\Omega$ is simply connected, this a well-known consequence of 
Carath\'eodory's 
kernel convergence theorem, see \cite[Thm 10.1]{duren}. 

Assume next that 
$\Omega$ is $m$-connected.
By Lemma \ref{conf} and the change of variable formula, it is enough to prove 
the result when $\Omega$ is a normalized circular domain
(so that $\Omega=\Omega'$ and $\varphi$ is the identity map in definition 
\eqref{defDeltat}). Let $\T_{a_j,r_j}$, $1\leq j\leq m$ 
denote the connected components of  $\partial\Omega$ lying inside $\D$.
By the decomposition theorem \cite[Sec. 10.5]{duren}, we can write
$f=f_1+\cdots+f_{m+1}$ with $f_{m+1}\in H^p(\D)$ and  $f_j\in H^p(\overline{\C}\setminus\overline{\D_{a_j,r_j}})$ for $1\leq j\leq m$. The result just quoted
in the simply connected case inplies that 
\[\sup_{n\in\NN}\|f_j\|_{L^p(\T_{a_j,r_j+\delta_\Omega/n})}<\infty,
\]
and since the $\T_{a_k,r_k+\delta_\Omega/n}$ are compactly embedded in 
$\C\setminus\overline{\D_{a_j,r_j}})$ when $k\neq j$, the inequality
\[\sup_{n\in\NN}\|f_j\|_{L^p(\T_{a_k,r_k+\delta_\Omega/n})}<\infty,
\]
follows from H\"older's inequality and  the Cauchy representation formula 
for $f_j$ from its values on
$\T_{a_j,r_j+\delta_\Omega/n_0}$ for some fixed $n_0$. Likewise
$\|f_{m+1}\|_{L^p(\T_{1-\delta_\Omega/n})}$, 
$\|f_{m+1}\|_{L^p(\T_{a_j,r_j+\delta_\Omega/n})}$ are uniformly bounded,
so that $\|f\|_p <\infty$ as desired.
\hfill\boite
\section{Traces of holomorphic functions}
\begin{lemma}
\label{IntH}
Let $\Omega$ be a Dini-smooth domain, and $g\in H^p(\Omega)$. Then
$g$ has a non-tangential limit a.e. on $\partial \Omega$
defining a trace function  $\tr_{\partial \Omega} g\in L^p(\partial\Omega)$.
In fact
\begin{equation}
\label{cvLpH}
\lim_{\varepsilon\to0}\|\tr_{\partial\Omega}g-g\circ P^{-1}_{\partial\Omega,\varepsilon}
\|_{L^p(\partial\Omega)}=0,
\end{equation}
and if $g$ is not identically zero then 
$\log|\tr_{\partial\Omega}g|\in L^1(\partial\Omega)$.

The quantity $\|\tr_{\partial\Omega}g\|_{L^p(\partial\Omega)}$
defines a norm on $H^p(\Omega)$ which is equivalent to $\|g\|_p$. 
Moreover
\begin{equation}
\label{MFLp}
\|{\mathcal M}_g\|_{L^p(\partial\Omega)}\leq C
\|\tr_{\partial\Omega}f\|_{L^p(\partial\Omega)},
\end{equation}
where $C$ depends on $\Omega$, $p$ and the 
aperture $\beta$ used in the definition of the maximal function.
\end{lemma}
{\it Proof.} 
It is well-known that functions in $H^p(\Omega)$ (recall from 
Theorem \ref{BS} that it coincides both with the Hardy and Smirnov class) 
have non-tangential limit in $L^p(\partial\Omega)$ of which $g$ is the Cauchy
integral \cite[Thm 10.4, Sec. 10.5]{duren}.
By Lemma \ref{conf}, we may assume that $\Omega$ is a
normalized circular domain. When $\Omega=\D$ all properties stated 
are standard, see  \cite[Thms 1.6,  2.2, 2.6]{duren}, \cite[Thm 3.1]{gar} 
and the remarks thereafter. By reflection, they also
hold for Hardy spaces of the complement of a disk. 
Next, assume that
$\Omega=\D\setminus\cup_{j=1}^n \overline{\D}_{a_j,r_j}$, with
$a_j\in\D$ and $0<r_j<1-|a_j|$. 
The decomposition theorem \cite[Sec. 10.5]{duren} 
tells us that $g=\sum_{j=0}^n g_j$ with
$g_0\in H^p(\D)$ and $g_j\in H^p(\overline{\C}\setminus\overline{\D_{a_j,r_j}})$, $1\leq j\leq n$.
From the known result in the simply connected case and the smoothness of 
holomorphic functions on their domain of analyticity, we thus obtain
\[
\lim_{\varepsilon\to0}\left\|\left(\tr_{\T}g_0+\sum_{j=1}^n{g_j}_{|\T}\right)
-\left(\sum_{j=0}^n{g_j}\right)\circ P^{-1}_{\T,\varepsilon}
\right\|_{L^p(\T)}=0,
\]
and by a similar argument we also get for $1\leq j\leq n$ that
\[
\lim_{\varepsilon\to0}\left\|\left(\tr_{\T_{a_j,r_j}}g_j+\left(\sum_{k\neq j}^n{g_k}\right)_{|\T_{a_j,r_j}}\right)
-\left(\sum_{j=0}^n{g_j}\right)\circ P^{-1}_{\T_{a_j,r_j},\varepsilon}
\right\|_{L^p(\T_{a_j,r_j})}=0,
\]
from which  \eqref{cvLpH} follows.

Next,  observe from
Hardy's convexity theorem (see \cite[Thms 1.5, 1.6]{duren} and the remark 
thereafter) that $\log\|g\|_{L^p(\T_r)}$ is a convex function of $\log r$
for $r\in(1-\delta_\Omega,1)$, hence by \eqref{cvLpH}
\begin{equation}
\label{majHc1}
\sup_{1-\delta_\Omega\leq r<1}\|g\|_{L^p(\T_r)}\leq \max\{\|g\|_{L^p(\T_{1-\delta_\Omega})},\|\tr_{\T}g\|_{L^p(\T)}\}.
\end{equation}
Likewise, for $j=1,\cdots,n$,
\begin{equation}
\label{majHc2}
\sup_{r_j+\delta_\Omega\geq r>r_j}\|g\|_{L^p(\T_{a_j,r_j})}\leq \max\{\|g\|_{L^p(\T_{a_j,r_j+\delta_\Omega})},\|\tr_{\T_{a_j,r_j}}g\|_{L^p(\T_{a_j,r_j})}\}.
\end{equation}
From \eqref{majHc1}-\eqref{majHc2} and H\"older's inequality applied
to the representation of $g$ as the Cauchy integral of $\tr_{\partial\Omega}g$,
we deduce that $\|g\|_p\leq C\|\tr_{\partial\Omega}g\|_{L^p(\partial\Omega)}$
where the constant $C$ depends only on $p$ and $\Omega$. In the other 
direction, the inequality $\|\tr_{\partial\Omega}g\|_{L^p(\partial\Omega)}
\leq \|g\|_p$ follows from the Fatou lemma applied to \eqref{cvLpH},
hence $\|\tr_{\partial\Omega}g\|_{L^p(\partial\Omega)}$ is equivalent to 
$\|g\|_p$. 

In the same vein, \eqref{MFLp} is easily obtained from the known simply 
connected case, the decomposition theorem and 
H\"older's inequality applied
to the representation of $g$ as the Cauchy integral of $\tr_{\partial\Omega}g$.

To prove that $\log|\tr_{\partial\Omega}g|\in L^1(\partial\Omega)$
if $g$ is not identically zero, we observe that $\Omega$ can be decomposed as
a finite union of Dini-smooth simply connected domains $\Omega_l$ such that
$\partial\Omega\subset\cup_l\Omega_l$. By Theorem \ref{BS} 
$g_{|\Omega_l}\in H^p(\Omega_l)$ since $|g|^p$ has a harmonic majorant on 
$\Omega$, {\it a fortiori} on $\Omega_l$. The result now follows from
the one in the simply connected case.
\hfill\boite
\begin{lemma}
\label{intsum}
Let $\Omega$ be a Dini-smooth domain, and $g\in H^p(\Omega)$. To each
$p_1\in [p,2p)$, there is a constant $c$ depending only on $\Omega$ and $p_1$
such that
$\|g\|_{L^{p_1}(\Omega)}\leq c\|g\|_p$.
\end{lemma}
{\it Proof.} 
When $\Omega=\D$, this is established in
\cite[appendix, proof of lem. 5.2.1]{BLRR}. So, by conformal mapping, we get it
for simply connected $\Omega$\footnote{Remember summability is understood 
with respect to area measure on the sphere.}. In the multiply connected case,
the result follows from its simply connected version and 
Lemma \ref{rmq:Hardy-ann}.
\hfill\boite
\begin{lemma}
\label{Rieszanneau}
Let $f$ be holomorphic in the annulus $\A_\varrho$, $0<\varrho<1$.
To each $p\in(1,\infty)$, there are constants $C_1,C_2$, depending on $f$,
$\varrho$  and  $p$, such that
\[\|\mbox{Im}\,f\|_{L^{p}(\T_r)}\leq C_1\|\mbox{Re}\,\|_{L^{p}(\T_r)}+
C_2,\quad \varrho<r<1.\]
\end{lemma}
{\it Proof.} 
Set $r_1=(1-\varrho)/3$, and pick 
$(1-\varrho)/2<r_2<1$. For $(1-\varrho)/2\leq|z|<r_2$, 
we get by the Cauchy formula
\[f(z)=F_2(z)-F_1(z),\qquad F_2(z)=\frac{1}{2i\pi}\int_{\T_{r_2}}\frac{f(\xi)}
{\xi-z}\,d\xi,\quad
F_1(z)=\frac{1}{2i\pi}\int_{\T_{r_1}}\frac{f(\xi)}
{\xi-z}\,d\xi.\]
Thus, setting $r=|z|$, it holds that
\[\|\mbox{Im}\,f\|_{L^{p}(\T_r)}\leq C(\rho,f,p)+
\|\mbox{Im}F_2\|_{L^{p}(\T_r)},\]
and since $F_2$ is holomorphic in $\D_{r_2}$ we obtain
from the M. Riesz theorem
\[\|\mbox{Im}\,f\|_{L^{p}(\T_r)}\leq  
C(\rho,f,p)+C(p)\|\mbox{Re}F_2\|_{L^{p}(\T_r)}\]
\[\leq
C(\rho,f,p)+C(p)\bigl(\|\mbox{Re}f\|_{L^{p}(\T_r)}+
\|\mbox{Re}F_1\|_{L^{p}(\T_r)}\bigr)
\leq C'(\rho,f,p)+C(p)\|\mbox{Re}f\|_{L^{p}(\T_r)}.
\]
A similar estimate holds for $\varrho<|z|\leq (1-\varrho)/2$ 
upon swaping the role of $F_1$ and $F_2$.
\hfill\boite

\section{A lemma on Sobolev functions}
\begin{lemma}
\label{estSobb}
Let $\Omega$ be a bounded Dini-smooth domain and assume that $p$, $r$ satisfy
\eqref{hypothesesnup}. Let $p_1\in[p,2p)$ be such that $2/p_1-1/p<1-2/r$ and 
set $1/\beta=1/p_1+1/r$.  Then 
\begin{itemize}
\item[(i)]$W^{1-1/\beta,\beta}(\partial\Omega)$ is
compactly included in $L^p(\partial\Omega)$;
\item[(ii)] $\widetilde{\Delta}_{1/n}$ being  as in \eqref{defDeltat},
there is a constant $C$ depending only of $\Omega$, $p$, and $\beta$ such that
for each $h\in W^{1,\beta}(\Omega)$
\begin{equation}
\label{SobHar}
\sup_{n\in\NN}\|\tr_{\partial\widetilde{\Delta}_{1/n}}h\|_
{L^p(\partial\widetilde{\Delta}_{1/n})}< C\|h\|_{W^{1,\beta}(\Omega)}.
\end{equation}
\end{itemize} 
\end{lemma}
{\it Proof:} let $\varphi$ conformally map $\Omega$ onto a normalized 
circular domain $\Omega'$. By Lemma \ref{conf} it is clear that 
$\|h\|_{W^{1,\beta}(\Omega)}$ and 
$\|h\circ\varphi^{-1}\|_{W^{1,\beta}(\Omega')}$ are comparable.
Likewise (see \eqref{defDeltat}), for any $l\in(1,\infty)$ and any 
smooth $\Phi$, 
$\|\tr_{\partial\widetilde{\Delta}_{1/n}}\Phi\|_
{L^l(\partial\widetilde{\Delta}_{1/n})}$ and
$\|\tr_{\partial{K}_{1/n}}\Phi\circ\varphi^{-1}\|_
{L^l(\partial{K}_{1/n})}$ on the one hand, 
$\|\tr_{\partial\widetilde{\Delta}_{1/n}}\Phi\|_
{W^{1,l}(\partial\widetilde{\Delta}_{1/n})}$ and
$\|\tr_{\partial{K}_{1/n}}\Phi\circ\varphi^{-1}\|_
{W^{1,l}(\partial{K}_{1/n})}$ on the other hand are comparable.
Hence, by interpolation,
$\|\tr_{\partial\widetilde{\Delta}_{1/n}}\Phi\|_
{W^{1-1/l,l}(\partial\widetilde{\Delta}_{1/n})}$ and
$\|\tr_{\partial{K}_{1/n}}\Phi\circ\varphi^{-1}\|_
{W^{1-1/l,l}(\partial{K}_{1/n})}$ are also comparable.
Altogether, we may assume for the proof that 
$\Omega$ is normalized 
circular. Moreover, in view of the 
extension theorem for Sobolev functions \cite[Sec. VI.3, Thm 5]{stein},
we may proceed componentwise on the boundary so it is enough to
consider the case where $\Omega=\D$.

From \cite[Thm 4.54]{Demangel}, we know 
if $\beta \geq 2$ that the inclusion 
$W^{1-1/\beta,\beta}(\T) \subset L^l(\T)$ is compact
for all $l \in (1, \infty)$, while if $\beta < 2$
the inclusion $W^{1-1/\beta,\beta}(\T) \subset L^{\beta/(2 - \beta)}(\T)$
is compact. One can check that $\beta/(2 -\beta) > p$
when $2/p_1-1/p<1-2/r$, thereby proving $(i)$.

From $(i)$ and the trace theorem, there is a constant $c=c(p,\beta)$ such that
\[\|\tr_\T h\|_{L^p(\T)}< c \|h\|_{W^{1,\beta}(\D)}, \qquad 
h\in W^{1,\beta}(\D).\]
Picking $r\in(0,1)$ and applying the above inequality to $h_r(z)=h(rz)$, we 
obtain (remember arclength is normalized)
\[\|\tr_{\T_r} h\|_{L^p(\T_r)}< c \|h_r\|_{W^{1,\beta}(\D)}
\leq \frac{c}{r^{2/\beta}} \|h\|_{W^{1,\beta}(\D)}.\]
Since $\partial{K}_{1/n}=\T_{1-1/2n}$ in the present case, assertion $(ii)$
follows. \hfill\boite

\section{Proof of Proposition \ref{paramwHp}}
\label{app_proof_prop_3}
By H\"older's inequality and standard properties of the Cauchy and 
Beurling transforms \cite[Ch. 1, (1.7)-(1.9)]{im2}, $T_\alpha$ maps
$L^p(\Omega)$ into $W^{1,rp/(r+p)}(\Omega)$. 
By the Rellich-Kondratchov theorem,
either $rp/(r+p)\geq2$ in which case
$W^{1,rp/(r+p)}(\Omega)$ is compactly embedded in $L^\lambda(\Omega)$,
$1\leq\lambda<\infty$, or else $rp/(r+p)<2$ and then 
$W^{1,rp/(r+p)}(\Omega)$ is compactly embedded in
every $L^\lambda(\Omega)$ with $1\leq\lambda<2rp/(2(r+p)-rp)$. Since
$2rp/(2(r+p)-rp)>p$ when $r>2$, this proves that $T_\alpha$ is compact from
$L^p(\Omega)$ into itself.

Next, we show that $I-T_\alpha$ is injective.
Indeed, if $h=T_\alpha h$, we get from what precedes that
$h\in L^\lambda(\Omega)$ for $1\geq1/\lambda>\max(0,1/p+1/r-1/2)$.
Thus, by H\"older's inequality, 
$\alpha \overline{h}\in L^t(\Omega)$ for every $t$ such that
$1\geq1/t>\max(1/r,1/p+2/r-1/2)$, and in turn
$h\in W^{1,t}(\Omega)$ for all such $t$. 
Therefore, by the Sobolev embedding theorem,
$h\in L^\lambda(\Omega)$ for each $\lambda$
such that $1\geq1/\lambda>\max(0,1/p+2(1/r-1/2))$. Iterating, we find that
$h\in L^\lambda(\Omega)$ whenever
$1\geq1/\lambda>\max(0,1/p+k(1/r-1/2))$ for some $k\geq1$, and that
$h\in W^{1,t}(\Omega)$ for $1\geq1/t>\max(1/r,1/p+1/r+k(1/r-1/2))$. Since
$r>2$, we deduce that $h\in W^{1,t}(\Omega)$ as soon as
$1\leq t<r$, in particular we may pick $t>2$.
The rest of the argument proceeds as in \cite[App. A]{BLRR}: we put
\[
H(z)=\frac1{2\pi i}\iint_{\Omega}
\frac{\alpha(\xi)\overline{h}(\xi)}{\xi-z}d\xi\wedge
d\overline{\xi}\,,\  \ z\in \C\,,
\]
noting  that $h=H_{|\Omega}$ and $H\in W^{1,t}_{loc}(\C)$ by what precedes.
Clearly $\overline{\partial}H=(\alpha\vee0)H$ on $\C$, and
since $t>2$ while $H$ vanishes at infinity we can apply the extended 
Liouville theorem  \cite[Prop. 3.3]{ap}
to the effect that $H\equiv0$ hence $h\equiv0$, as desired.

It now follows from a theorem of F. Riesz  
\cite[Thm 1.16]{COLTKRESS} that
$I-T_\alpha$ is an isomorphism of $L^p(\Omega)$. 

In another connection, let $w\in G^p_\alpha(\Omega)$ and set $g=w-T_\alpha w$.
Then $\overline{\partial}g=0$  because 
$\overline{ \partial}T_\alpha w=\alpha\overline{w}$, hence
$g$ is holomorphic in $\Omega$. Moreover, we know from Property \ref{2} that
$w \in L^{p_1}(\Omega)$ for $p \leq p_1 < 2 p$, hence for such $p_1$ it holds
that $T_\alpha w \in W^{1,\beta}(\Omega)$ with $1/\beta = 1/{p_1} + 1/r$.
Choosing $p_1$ such that $2/p_1-1/p<1-2/r$,
Lemma \ref{estSobb} point $(ii)$ implies that $T_\alpha w$ satisfies 
\eqref{systema}, hence so does $g$, that is, $g\in H^p$.

Conversely, for $g\in H^p$, let us put $w=(I-T_\alpha)^{-1}g$.
Since $g\in L^{p_1}(\Omega)$ for 
$p \leq p_1 < 2 p$  and \eqref{hypothesesnup} continues to
hold with $p$ replaced by 
$p_1$, we get from the previous part of the proof that
$w$ lies in $L^{p_1}(\Omega)$ and consequently that $T_\alpha w$ satisfies 
\eqref{systema}. Hence $w$ in turn meets \eqref{systema},
and since $\overline{\partial}(I-T_\alpha)w=0$ it is a solution to 
\eqref{eq:w}, hence a member of $G_\alpha^p(\Omega)$.

Equation \eqref{repCauchyg} simply means that 
${\mathcal C}(\tr_{\partial\Omega} T_\alpha w)(z)=0$ for $z\in\Omega$.
To see this, let 
\[F(z)=
\frac1{2\pi i}\iint_{\Omega}
\frac{\alpha(\xi)\overline{w}(\xi)}{\xi-z}d\xi\wedge
d\overline{\xi}\,,\  \ z\in \C.\,.
\]
By what precedes  
$F(z)\in W^{1,\beta}_{loc}(\C)$, $F_{|\Omega}=T_\alpha w$, and clearly 
$F$ is holomorphic
in $\overline{\C}\setminus\overline{\Omega}$ with $F(\infty)=0$.
Lemma \ref{estSobb} point $(ii)$, applied with $h=F$ on
$(\overline{\C}\setminus\overline{\Omega})\cap \D_R$ where $R$ is a 
large positive number, shows that 
$F\in H^p(\overline{\C}\setminus\overline{\Omega})$, and we have that
$\tr_{\partial\Omega} F=\tr_{\partial\Omega} T_\alpha w$.
Consequently, by Cauchy's theorem, it holds for any system of 
rectifiable Jordan curves $\Gamma$ homotopic to $\partial\Omega$ in 
$\overline{\C}\setminus\Omega)$ that
\[{\mathcal C}(\tr_{\partial\Omega} T_\alpha w)(z)=
\frac1{2\pi i}\int_\Gamma
\frac{F(\xi)}{\xi-z}d\xi
\,,\  \ z\in \Omega.\,.
\]
Deforming the inner components of $\partial\Omega$ to a point and the 
outer component to $\infty$, we see that the above integral is zero, 
as desired.

Finally, since $\int_{\partial\Omega} T_\alpha w=0$ by what we just said,
we get from the Poincar\'e inequality, H\"older's inequality, the continuity
of $(I-T_\alpha)^{-1}$ and
Lemma \ref{intsum} that
\[\|T_\alpha w\|_{W^{1,\beta}(\Omega)}\leq c_1 \|\alpha\overline{w}\|_{L^\beta(\Omega)}\leq c_2\|\alpha\|_{L^r(\Omega)}\|w\|_{L^{p_1}(\Omega)}\leq
c_3\|g\|_{L^{p_1}(\Omega)}\leq c_4  \|g\|_p.\]
Thus, by the trace theorem and Lemma \ref{estSobb} point $(i)$, we obtain
$\|\tr_{\partial\Omega}T_\alpha w\|_{L^p(\partial\Omega)}\leq c_5 \|g\|_p$,
hence $\|\tr_{\partial\Omega} w\|_{L^p(\partial\Omega)}\leq c_6\|g\|_p$
since $w=g+T_\alpha w$. Property \ref{2} now yields
$\|w\|_p\leq C\|g\|_p$, as announced.
\hfill\boite

\section{Proof of Lemma \ref{lemu}}
\label{app5}
By Lemma \ref{conf}, we may assume that $\Omega=\D$.
Let ${\bf P}_{+}(h)=\tr_{\TT}{\mathcal C}(h)$ denote
the Riesz projection  which discards Fourier coefficients of
non-positive index. It is continuous from $L^p(\TT)$ onto 
$\tr_{\TT} \ H^p(\DD)\subset L^p$  \cite[Sec. III.1]{gar}. Moreover,
to each $u\in L^p(\T)$ there uniquely exists $\widetilde{u}\in L^p(\T)$ 
such that $u+i\widetilde{u}\in {\bf P}_{+} (H^p)$ and $\int_\T\widetilde{u}=0$.

For $u\in L^p(\T)$ and $c\in\R$, let $w_{u,c}\in G^p_\alpha(\D)$ 
satisfy  $ (u+i(\widetilde{u}+c)) = P_+ w_{u,c}$. 
Such a function uniquely exists 
by Proposition \ref{paramwHp} and depends continuously on $u$ and $c$. Define
\[A(u,c):=\left(\mbox{Re}\,\bigl(\tr_\T\, w_{u,c}\bigr)\,\,,\,
\mbox{Im}\int_\T 
\tr_\T\, \sigma^{1/2}
w_{u,c}\right)\in L^p_\R(\T)\times\R.
\]
Since $(I-T_\alpha)w_u=g$
where $g\in H^p$ satisfies $\tr_\T\, g=u+i(\widetilde{u}+c)$, we can decompose
the operator $A$ as $A(u,c)=(u,c)+B(u,v)$ where
\[B(u,c):=\left(
\mbox{Re}\,\bigl(\tr_\T \, T_\alpha( w_{u,c} )\bigr)\,\,,
\,\, 
\mbox{Im}\int_\T 
\sigma^{1/2}\tr\,  w_{u,c} \,-c\right).\]
From the proof of Proposition \ref{paramwHp}, we know  that 
$(u,c)\mapsto T_\alpha w_{u,c}$ 
is continuous from
$L^p(\T)\times \R$ into $W^{1,\beta}(\D)$ when $1/\beta=1/p_1+1/r$ for some
$p_1\in[p,2p)$, hence $B$ is compact from $L^p(\T)\times \R$ into itself 
by Lemma \ref{estSobb} point $(i)$ and the trace theorem.

In another connection, if $w\in G^p_\alpha(\D)$ is such that
$\mbox{Re}\,(\tr_\T\, w)=0$ and $\mbox{Im}\int_\T \sigma^{1/2}
\tr_\T\, w=0$, then $w=0$. Indeed,
normalizing $\tr_\T\,s$  to be real 
in \eqref{decompexpH} we find that $F\in H^p$ has zero real part on $\T$,
hence it is an imaginary constant, and in fact $F=0$ as the mean on $\T$ 
of $\sigma^{1/2} F$ must vanish. 
Consequently $A$ is injective, hence a homeomorphism of 
$L^p(\T)\times\R$ by Riesz's theorem. This shows one can impose uniquely in
$L^p(\T)$ the real part of $w\in G^p_\alpha(\D)$ on $\T$ together with
the mean of 
$\sigma^{1/2}$ times its imaginary part there, and that
\[\|\tr_\T\,w\|_{L^p(\T)}\leq c_1\|\mbox{Re}\tr_\T\,w\|_{L^p(\T)}+c_2\left|
\int_\T\sigma^{1/2}\mbox{Im}\tr_\T\,w\right|.\]
From this and \eqref{ftow}, it follows
easily that one can impose uniquely
the real part $u\in L^p(\T)$ of $f\in H^p_\alpha(\D)$ on $\T$ and the mean of 
its imaginary part there. In addition, if the latter is taken to be zero, 
there is an inequality of the form
$\|\tr_{\partial\Omega} f\|_{L^{p}(\partial\Omega)}\leq c\|u\|_{L^{p}(\partial\Omega)}$ which proves assertion $(ii)$.




Now, in view of Properties \ref{1},\ref{2}, and \ref{6} in
section \ref{sechpnud}, taking real parts in assertion 
$(ii)$  yields assertion $(i)$ except for the uniqueness part.
To establish the latter, assume $u=0$ and let us prove that $U=0$.
As $u$ satisfies \eqref{div}, there is a distribution $V$ such that
\eqref{system} holds. Since $\|U\|_p<\infty$, hence {\it a fortiori} 
$U\in L^p(\D)$, we observe much as in the proof of 
\cite[Thm 4.4.2.2]{BLRR} that
$V\in L^p(\D)$; the only difference is that, in order to obtain 
equation (63) {\it loc. cit.}, one must know whenever $\Phi$ is smooth with 
compact support that 
$\|\Phi\nabla\sigma\|_{L^q(\D)}\leq c\|\nabla\Phi\|_{L^q(\D)}$,
$1/p+1/q=1$, which follows easily from the H\"older and the 
Sobolev inequalities for $r>2$. 
Then $f=U+iV$ satisfies
\eqref{dbar}, so that $w$ given by \eqref{ftow} satisfies \eqref{eq:w}.
Since $U={\rm Re }f$ satisfies \eqref{esssuplp} by assumption, so does
${\rm Re }\,w$. By Proposition \ref{expsf} $w$ assumes the form
\eqref{decompexpH} where $\tr_\T\,{\rm Im}\,s=0$ and, say,
$F=a+ib$ is holomorphic. Assume for a contradiction that 
$\|a\|_{L^p(\T_\rho)}$, $0\leq\rho<1$  is unbounded.
Then it must tend to $+\infty$ as it increases with $\rho$
\cite[Thms 1.5, 1.6]{duren}. 
By the continuity of $s$
({\it cf.} Remark \ref{remHolder}), to each
$\varepsilon\in (0,1)$ there is $\rho_0 \in (0,1)$ such that
$|\mbox{Im} \, \exp(s(z))|<\varepsilon |\exp (s(z))|$ as soon as 
$\rho_0<|z|\leq1$. For such $z$, we deduce from \eqref{decompexpH} that
\begin{equation}
\label{estimpu}
|\mbox{Re} \, w(z)|\geq e^{-\|s\|_{L^\infty(\D)}}
((1-\varepsilon^2)^{1/2}|a(z)|-\varepsilon|b(z)|).
\end{equation}
By a theorem of  M. Riesz $\|b-b(0)\|_{L^p(\T_\rho)}\leq C \|a\|_{L^p(\T_\rho)}$
with $C=C(p)$, uniformly with respect to $\rho\in(\rho_0,1]$
\cite[Thm 4.1]{duren}. Hence integrating
\eqref{estimpu}, we get
\[\|\mbox{Re}\, w\|_{L^p(\T_\rho)} \geq e^{-\|s\|_{L^\infty(\D)}}
\Bigl((1-\varepsilon^2)^{1/2}-|b(0)|/\|a\|_{L^p(\T_\rho)}-\varepsilon C
\Bigr)\|a\|_{L^p(\T_\rho)}\]
and taking $\varepsilon$ small enough we find this absurd since the
left hand-side is bounded while the right hand side goes to infinity
when $\rho\to1$. Hence $\|a\|_{L^p(\T_\rho)}$ is bounded and so is
$\|b\|_{L^p(\T_\rho)}$, by the M. Riesz theorem again,
in other words $F\in H^p(\D)$. 
Therefore 
 $w\in G_\alpha^p(\D)$ by Proposition \ref{expsf}, thus $f\in H^p_\nu(\D)$.
Moreover, since $\mbox{Re}\,w$ has 
nontangential limit $0$  on $\T$, so does $\mbox{Re}\,F$ since $\tr_\T\,e^s>0$,
thus $F$ is an imaginary constant.
However $V$ was defined up to an additive constant only, and since
we have just shown that $\tr_\T V\in L^p(\T)$ we can pick this constant so that
$\int_\T V \sqrt{(1+\nu)/(1-\nu)} =0$. 
Then $\tr_\T {\rm Im}\,w$ has 
zero mean by \eqref{ftow}, consequently $F=0$, hence $w=f=0$.
In particular $U=0$, as desired.
\hfill\boite

\section{Proof of Proposition \ref{casinvHp}}
\label{app6}
By Proposition \ref{trick-hardy}, the conclusion for $H^p_\nu(\D)$ 
follows from  the result for  $G_\alpha^p(\D)$ which we now prove.

Define a Hermitian duality pairing on 
$L^p(\DD) \times L^q(\DD)$,  $1/p+1/q = 1$, by the formula
\[<h,g>_\DD=\frac1{2\pi i}\iint_{\DD} h(z)\overline{g(z)}
\,dz\wedge d\overline{z}.\]
If  $A:L^p(\D)\to L^p(\D)$ is antilinear ({\it i.e.} real linear such that
$A(\lambda h)=\bar{\lambda} A(h)$), then $A(h)=B(h)+iB(ih)$ where
$B=\mbox{Re} A$. We let $A^\sharp$ 
designate the antilinear operator on $L^q(\D)$ 
such that $<Ah,g>=<A^\sharp g,h>$.
It is easy to check that $A^\sharp(g)=B^*(g)+iB^*(ig)$, where 
$B^*$ is the adjoint of $B$ when $L^p(\D)$, $L^q(\D)$
are viewed as real vector spaces endowed with the pairing $\mbox{Re}<.\,,.>$.

For $\alpha\in L^r(\D)$ and $h\in L^p(\D)$, define functions $T_\alpha(h)$ and 
$\mathcal{T}_\alpha(h)$ on $\D$ and $\C$ by
\begin{equation}
\label{defT}
\frac1{2\pi i}\iint_{\DD} \frac{\alpha(\xi)\overline{h(\xi)}}
{\xi-z}d\xi\wedge 
d\overline{\xi} = 
\left\{
\begin{array}{ll}
\mathcal{T}_\alpha h(z) & \mbox{ for } z \in \CC \, , \\
 & \\
T_\alpha h(z) & \mbox{ for } z \in \DD \,  : \, 
T_\alpha h = (\mathcal{T}_\alpha h)_{|_\DD} \, .
\end{array}
\right.
\end{equation}


From Proposition \ref{paramwHp}, we get that $T_\alpha$
is compact from $L^p(\D)$ into itself and clearly it is antilinear.
Moreover $I-T_\alpha$ is an isomorphism of $L^p(\DD)$ and the restriction
map $I-T_{\alpha}: G_{\alpha}^p(\DD) \to H^p(\DD)$  is an isomorphism
which coincides with the analytic (Cauchy) projection, see \eqref{repCauchyg}.
%

By a theorem of Schauder ($\mbox{Re}\,T_\alpha)^*$ is compact \cite[Thm. VI.4]{brezis}, hence also $T_\alpha^\sharp$.
In addition $I-T_\alpha^\sharp$ is injective, for if $g=T^\sharp_\alpha(g)$
we get from the definition of $A^\sharp$ 
\[\langle (I-T_\alpha)h,g\rangle=\langle h,g\rangle-
\langle g,h\rangle,\qquad h\in L^p(\D),
\]
which is absurd if $g\neq 0$ since the right hand side is pure imaginary 
whereas $I-T_\alpha$ is surjective.
Hence $I-T_\alpha^\sharp$ is an isomorphism of $L^q(\D)$ by Riesz's theorem.
By Fubini's theorem, we obtain for $h \in L^p(\DD)$ and $g \in L^q(\DD)$
that
\begin{eqnarray*}
<T^\sharp_\alpha g,h>_\DD = <T_\alpha h,g>_\DD &=&
-\frac1{4\pi^2}\iint_{\DD} \left(
\iint_{\DD} \frac{\alpha(\xi)\overline{h(\xi)}}{\xi-z}d\xi\wedge 
d\overline{\xi}\right)
\overline{g(z)}\,dz\wedge 
d\overline{z}\\
&=&
-\frac1{4\pi^2}\iint_{\DD} \left(
\iint_{\DD} \frac{\overline{g(z)}}{\xi-z}dz\wedge 
d\overline{z}\right)\alpha(\xi)\overline{h(\xi)}
\,d\xi\wedge 
d\overline{\xi}\\
&=&
-\frac1{4\pi^2}\iint_{\DD} \left(-\alpha(\xi)
\iint_{\DD} \frac{\overline{g(z)}}{z-\xi}dz\wedge 
d\overline{z}\right)\overline{h(\xi)}
\,d\xi\wedge 
d\overline{\xi}
\end{eqnarray*}
so that 
\begin{equation}
\label{adj}
T^\sharp_\alpha g=-\alpha T_{\chi_\D}g
\end{equation}
where $\chi_\D$ is the characteristic function of $\D$.

Let $\alpha_\varrho$ be  $0$ on $\DD_\varrho $ and $\alpha$ on $\AA_\varrho$.
If $p_1\geq p$, then 
$T_{\alpha_\varrho}$ maps $L^{p_1}(\D)$ into $W^{1,\beta}(\D)$
with $1/\beta=1/r+1/p_1$, see proof of Proposition \ref{paramwHp}.
Besides, $\left(\mbox{Ran} T_{\alpha_\varrho}\right)_{|\D_\varrho}$ 
consists of holomorphic functions. 
Therefore, for $p_1$  as in Lemma \ref{estSobb}, we get from Property \ref{6} and \eqref{SobHar} that 
\begin{equation}
\label{restrick}
 \Bigl(T_{\alpha_\varrho}(I-T_\alpha)^{-1}H^{p}(\D)\Bigr)_{|\D_\varrho}= 
 \Bigl(T_{\alpha_\varrho}G^{p}_\alpha(\D)\Bigr)_{|\D_\varrho}\subset 
 \Bigl(T_{\alpha_\varrho}L^{p_1}(\D)\Bigr)_{|\D_\varrho}\subset 
H^p(\D_\varrho).
\end{equation}    
Moreover, the operator, 
$A_\varrho=I+T_{\alpha_\varrho}(I-T_\alpha)^{-1}$ maps $H^p(\DD)$ into
$L^{p_1}(\D)$. Let us introduce the operator 
$B_\varrho=J\circ A_\varrho$, where 
$J:L^{p_1}(\DD)\to L^{p_1}( \DD_\varrho)$ is the natural restriction. 
In view of \eqref{restrick}, $B_\varrho$ maps continuously
$H^p(\DD)$ into $H^p( \DD_\varrho)$.
\begin{lemma}
\label{denselisse}
The operator $B_\varrho:H^p(\DD)\to H^p(\DD_\varrho )$ has dense range.
\end{lemma}
{\it Proof:} It is equivalent to prove that if 
$\Psi\in H^{q,00}(\CC\setminus  \overline{\DD}_\varrho)$ satisfies
\begin{equation}
\label{HB}
\frac1{2i\pi}\int_{\TT_\varrho }\Psi(z)\,\tr_{\T_\varrho}B_\varrho g(z)\,dz=0,~~~~\forall g\in H^p(\DD),
\end{equation}
then $\Psi=0$; indeed, the line integral of the product
over $\T_\varrho$ identifies 
$H^{q,00}(\CC\setminus \overline{\DD}_\varrho )$ (non isometrically) with the 
dual of 
$H^p(\DD_\varrho )$ \cite[Theorem 7.3]{duren}. 
For small $\varepsilon>0$, since $\Psi$ and $g$ are  smooth on 
$\overline{\A}_{\varrho+\varepsilon,1-\varepsilon}$ while 
$A_\varrho g-g\in W^{1,\beta}(\D)$, 
we  get from Stokes' theorem  
\[
\dfrac{1}{2i\pi}\int_{\partial \AA_{\varrho+\varepsilon, 1-\varepsilon}} 
\!\!\Psi(z)\,\tr
A_\varrho g(z)   \ dz =- \dfrac{1}{2i\pi}\iint_{\AA_{\varrho+\varepsilon,1-\varepsilon}} \overline{\partial} \bigl(\Psi(z)A_\varrho g(z)\bigr) \ dz \wedge d\overline{z}.
\]
where the trace of $A_\varrho g$ in the first integral is on
$\partial \AA_{\varrho+\varepsilon, 1-\varepsilon}$.
Because $\Psi$, $g$ are holomorphic in $\AA_{\varrho}$ and
$\overline{\partial}T_{\alpha_\varrho} h =  \alpha_\varrho \overline{h}$ for 
$h \in L^p(\DD)$ by standard properties of the Cauchy transform,
we may compute the surface integral using the definition of $A_\varrho$ 
to obtain
\[0=\dfrac{1}{2i\pi}\int_{\partial \AA_{\varrho+\varepsilon, 1-\varepsilon}} 
\!\!\Psi(z) \,\tr A_\varrho g(z) \ dz + \dfrac{1}{2i\pi}\iint_{\AA_{\varrho+\varepsilon,1-\varepsilon}}\Psi(z)  \alpha(z)\overline{(I-T_\alpha)^{-1}g(z)}
\ dz \wedge d\overline{z}.
\]
Now, $\Psi_{|\A_\varrho}\in H^q(\A_\varrho)$ and 
$A_\varrho g\in H^p(\D)+W^{1,\beta}(\D)$, while
$(I-T_\alpha)^{-1}g\in G^p_\alpha$. 
Let $\Psi_\varrho$ be 0 on $\D_\varrho$ and $\Psi$ elsewhere.
Pick $q_1\in(2,2q)$ 
with $2/q_1-1/q<1-2/r$ and recall from Lemma \ref{intsum} that
$\Psi_{|\A_\varrho}\in L^{q_1}(\A_\varrho)$, so that
$\alpha\Psi_\varrho\in
L^{\delta}(\D)$ where $1/\delta=1/r+1/q_1<(p_1-1)/p_1$.
Thus, letting $\varepsilon\to0$, we get from Property \ref{1}, \ref{6},
Lemma \ref{estSobb}, and H\"older's inequality   that
\begin{equation}
\label{CGA}
\frac1{2i\pi}\!\!\int_{ \TT_\varrho} \!\Psi(z)\, \tr_{\T_\varrho}B_\varrho g(z) \, dz=
\frac1{2i\pi}\!\!\int_{\TT}\Psi(z)\, \tr_{\T}A_\varrho g(z)  dz +
<\alpha\Psi_\varrho,(I-T_\alpha)^{-1}g>_\DD,
\end{equation}
where we took into account that
$\tr_{\T_\varrho}B_\varrho g =\tr_{\T_\varrho}A_\varrho g$ on $\TT_\varrho $.
Put
$b=(I-T_\alpha^\sharp)^{-1}(\alpha\Psi_\varrho)\in L^{p_1/(p_1-1)}(\D)$. 
Then, by definition of $A_\varrho^\sharp$, it holds that
\begin{eqnarray}
<\alpha\Psi_\varrho,(I-T_\alpha)^{-1}g>_\DD&=&
<(I-T_\alpha^\sharp) b, (I-T_\alpha)^{-1}g>_\DD\label{transfPA}\\
&=&<b,(I-T_\alpha)^{-1}g>_\DD-<T_\alpha(I-T_\alpha)^{-1}g,b>\nonumber\\
&=&<g, (I-T^\sharp_\alpha)^{-1}(\alpha\Psi_\varrho)>_\DD+
2i\mbox{Im}<b,(I-T_\alpha)^{-1}g>_\DD.\nonumber
 \end{eqnarray}
Representing $g\in H^p(\D)$ by the  Cauchy integral of $\tr_\T g$,
we further have that
\begin{eqnarray*}
<g, (I-T^\sharp_\alpha)^{-1}(\alpha\Psi_\varrho)>_{_\mathbb{D}} & = & \dfrac{1}{2i\pi} \iint_\mathbb{D} \left( \int_\mathbb{T} \dfrac{\tr_\T g(z)}{z-\xi} \ dz \right) \overline{(I-T^\sharp_\alpha)^{-1}(\alpha\Psi_\varrho)(\xi)} \ d\xi \wedge d\overline{\xi} \\
\\
& = & \dfrac{1}{2i\pi}  \int_\mathbb{T} \tr_\T g(z) \left( \iint_\mathbb{D} \dfrac{\overline{(I-T^\sharp_\alpha)^{-1}(\alpha\Psi_\varrho)(\xi)}}{z - \xi} \ d\xi \wedge d\overline{\xi} \right) \ dz \ ,
\end{eqnarray*}
and letting
$\mathfrak{I}=2i\mbox{Im}<b,(I-T_\alpha)^{-1}g>_\DD$ 
we get in view of \eqref{transfPA}
\begin{equation}
 \label{eq:CGA2}
<\alpha\Psi_\varrho,(I-T_\alpha)^{-1}g>_\DD=
-\dfrac{1}{2i\pi}\int_\TT \tr_\T g(z)\tr_\T T_{\chi_\D}(I-T^\sharp_\alpha)^{-1}(\alpha\Psi_\varrho)(z) \ dz + \mathfrak{I}\ .
\end{equation}
Next, by \eqref{defT}, the function $\mathcal{T}_{\alpha_{\varrho}} h\in
W^{1,\beta}_{loc}(\C)$ has the same trace on $\T$ as $T_{\alpha_{\varrho}} h$
for all $h\in L^{p_1}(\D)$.
Hence,  the first integral in the right hand side of
(\ref{CGA}) can be rewritten as
\begin{equation}
 \label{decCauchy}
\dfrac{1}{2i\pi}\int_\TT\Psi(z) \,\tr_\T g(z) \ dz
+\dfrac{1}{2i\pi}\int_\TT\Psi(z)\, \tr_\T\mathcal{T}_{\alpha_\varrho}(I-T_\alpha)^{-1}g(z) \ dz \ .
\end{equation} 

Moreover,  an argument  similar to the one that led us to \eqref{restrick} 
easily yields that
\begin{equation}
\label{simdelta}
\left(\mathcal{T}_{\alpha_\varrho}L^{p_1}(\D)\right)_{|\overline{\C}\setminus\overline{\D}}
\subset H^{p,00}(\overline{\C}\setminus\overline{\D}),
\end{equation} 
therefore $\Psi_{|\T}\tr_\T \mathcal{T}_{\alpha_\varrho}(I-T_\alpha)^{-1}g$
is the trace on $\T$ of a function
in $H^{1,00}(\overline{\C}\setminus\overline{\D)}$ and
the second integral in (\ref{decCauchy}) is zero by Cauchy's theorem
(deform $\T$ to infinity).
Thus, 
\begin{equation}
 \label{decCauchy2}
 \dfrac{1}{2i\pi}\int_{\TT} \Psi(z) \,\tr_\T A_\varrho g(z) \ dz = \dfrac{1}{2i\pi}\int_\TT\Psi(z)\,\tr_\T  g(z) \ dz \ ,
\end{equation}
and in view of (\ref{CGA}), (\ref{eq:CGA2}), and (\ref{decCauchy2}), we
conclude if \eqref{HB} holds that
\begin{equation}
 \label{eq:preuve_densite_A_2_bis}
 \frac1{2i\pi}\int_\TT \tr_\T g(z)\,\left(\Psi(z)-\tr_\T T_{\chi_\D}
(I-T^\sharp_\alpha)^{-1}(\alpha_\varrho\Psi)(z)\right) \ dz=-\mathfrak{I} \ , \quad g\in H^p(\DD).
\end{equation}
Then, observe from \eqref{defT} that $\mathcal{T}_{\chi_\D} h\in
W^{1,\delta}_{loc}(\C)$ has the same trace on $\T$ as $T_{\chi_\D} h$
for all $h\in L^\delta(\D)$. In addition, since either $\delta>2$ or
$\delta/(2-\delta)> q$, we get from \cite[Thm 4.54]{Demangel} as in the 
proof of Lemma  \ref{estSobb} point $(i)$ that
\begin{equation}
\label{simq}
\left(\mathcal{T}_{\chi_\D}L^{q_1}(\D)\right)_{|\overline{\C}\setminus\overline{\D}}
\subset H^{q,00}(\overline{\C}\setminus\overline{\D})
\end{equation} 
(compare  \eqref{simdelta}).
Hence  \eqref{eq:preuve_densite_A_2_bis} can be rewritten as
\begin{equation}
 \label{formetest}
 \frac1{2i\pi}\int_\TT \tr_\T g(z)\,\left(\Psi(z)-\tr_\T \mathcal{T}_{\chi_\D}
(I-T^\sharp_\alpha)^{-1}(\alpha\Psi_\varrho)(z)\right) \ dz=-\mathfrak{I} \ , \quad g\in H^p(\DD),
\end{equation}
and the left hand side of \eqref{formetest} is
a complex linear form $\mathcal{L}(g)$ on $H^p(\D)$
while the right hand side is always
pure imaginary. Therefore $\mathcal{L}$ is the zero form,
which means that 
\[
\Psi-\mathcal{T}_{\chi_\D}(I-T^\sharp_\alpha)^{-1}(\alpha\Psi_\varrho)
\in H^{q,00}(\overline{\C}\setminus\overline{\D})(\alpha\Psi_\varrho)
\]
must be the zero function. Translating back to $\D$,
this amounts to say that the function
\begin{equation}
\label{defG}
G(z)=\Psi_\varrho(z)-T_{\chi_\D}(I-T^\sharp_\alpha)^{-1}\left(\alpha\Psi_\varrho\right),
\end{equation}
which lies in $\in \Bigl(H^q(\A_\varrho)\vee 0_{|\D_\varrho}\Bigr)+
W^{1,\delta}(\D)$, satisfies $\tr_\T G=0$. 
A short calculation using \eqref{adj} and the identity $(A-B)^{-1}=A^{-1}+A^{-1}B(A-B)^{-1}$ shows that 
\begin{equation}
\label{eqG}
\alpha G=\alpha\Psi_\varrho+T_\alpha^\sharp(I-T^\sharp_\alpha)^{-1}(\alpha\Psi_\varrho)=
(I-T^\sharp_\alpha)^{-1}(\alpha\Psi_\varrho),
\end{equation}
hence 
\begin{equation}
\label{eqdG}
G=\Psi_\varrho-T_{\chi_\D}\left(\alpha G\right)=
\Psi_\varrho-T_{\bar{\alpha}}\left(G\right)
\end{equation}
by \eqref{defG} and \eqref{eqG}. Note, since $q_1>2$ and either $\delta>2$ or $2\delta/(2-\delta)>2q$,
 that $G\in L^\lambda(\D)$ for some $\lambda>{r/(r-1)}$
by Lemma \ref{intsum} and the Sobolev embedding theorem. 
Set for simplicity $G_1=G_{\A_\varrho}$. Applying $\overline\partial$ to
\eqref{eqdG}, we find that $\overline{\partial} G_1=
-\overline{\alpha} \,\overline{G_1}$.
Thus, by Proposition \ref{expsf}, it holds that $G_{|\A_\varrho}=e^sF$ 
where $s$ is bounded
and $F$ is holomorphic in $\A_\varrho$. Moreover, as we noticed already that
either $\delta>2$ or $\delta/(2-\delta)>q$, we get as in the 
proof of Lemma  \ref{estSobb} point $(i)$ that $\|G\|_q<\infty$ on 
$\A_\varrho$. In particular 
$F\in H^q(\A_\varrho)$, and since $\tr_\T F=0$ we have that $F=0$ 
so that $G_1=0$. Plugging this in \eqref{eqdG}, we get that 
$(I+T_{\bar{\alpha}})(0_{\A_\varrho}\vee G_{\D_\varrho})=0$
hence $G_{\D_\varrho}=0$ since $(I+T_{\bar{\alpha}})$ is injective 
on $L^\lambda(\D)$  by Proposition \ref{paramwHp}. Altogether $G=0$, hence 
$\Psi_\varrho=0$, and finally $\Psi=0$, by analytic continuation.
\hfill\boite
\medskip
{\it Proof of Proposition \ref{casinvHp}:} for $0<r\leq1$, we
let for simplicity $T_{\alpha,r}=\T_{\alpha_{|\D_r}}: L^p(\DD_r)\to
L^p(\DD_r)$, so that $\T_{\alpha,1}=T_\alpha$.
so that $T_{\alpha,1}=T_\alpha$. Recall also from \eqref{repCauchyg} the 
notation $\mathcal{C}$ for Cauchy integrals.
By Proposition \ref{paramwHp}, a function 
$w_\varrho\in G^p_{\alpha_{\D_\varrho}}(\D_\varrho)$ lies in 
$\bigl(G^p_\alpha(\D)\bigr)_{|\D_\varrho}$ if, and only if there is 
$g\in H^p(\D)$ such that
\begin{equation}
\label{carres}
(I-T_{\alpha,\varrho})^{-1}\mathcal{C}(\tr_{\T_\varrho}w_\varrho)=
\left((I-T_\alpha)^{-1}g\right)_{|\DD_\varrho}.
\end{equation}
Define  $\widetilde{\alpha}_\varrho=\chi_{\DD_\varrho}\alpha$. 
Since 
$T_{\alpha,\varrho}(h_{|\DD_\varrho})=(T_{\widetilde{\alpha}_\varrho}h)_{|\DD_\varrho}$ for $h\in L^p(\DD)$, equation \eqref{carres} means that
\begin{equation}
\label{calcP+}
\mathcal{C}(\tr_{\T_\varrho}w_\varrho)=
\Bigl((I-T_{\widetilde{\alpha}_\varrho})(I-T_\alpha)^{-1}g\Bigr)_{|\DD_\varrho}.
\end{equation}
Observe that
$T_\alpha=T_{\alpha_\varrho}+T_{\widetilde{\alpha}_\varrho}$
where the notation $\alpha_\varrho=\chi_{|\A_\varrho}\alpha$ was introduced in the proof of Lemma \ref{denselisse}. Hence, using the identity 
$(A-B)^{-1}=A^{-1}+A^{-1}B(A-B)^{-1}$, we obtain
\[(I-T_\alpha)^{-1}=(I-T_{\widetilde{\alpha}_r})^{-1}+
(I-T_{\widetilde{\alpha}_r})^{-1}T_{\alpha_r}(I-T_\alpha)^{-1}.
\]
Substituting in \eqref{calcP+}, we get
\[\mathcal{C}(\tr_{\T_\varrho}w_\varrho)=B_\varrho g,
\]
and we conclude from Lemma \ref{denselisse} that 
$\bigl(G^p_\alpha(\D)\bigr)_{|\D_\varrho}=(I-T_{\alpha,\varrho})^{-1}\mbox{Ran}B_\varrho$ is dense in $G^p_{\alpha_{\D_\varrho}}(\D_\varrho)$, as desired.
\hfill\boite

\end{document}